\documentclass[12pt,reqno]{amsart}
\usepackage{mathrsfs}
\usepackage{cases}
\usepackage{epic}
\usepackage{amsfonts}
\usepackage{graphicx}
\usepackage{amsmath}
\usepackage{amssymb, upgreek}
\usepackage{bm}
\usepackage{latexsym,todonotes}
\usepackage{pdflscape}
\usepackage[all]{xypic}
\usepackage[all]{xy}
\usepackage{color}
\usepackage{colordvi}
\usepackage{multicol}
\usepackage[normalem]{ulem}
\newcommand{\arxiv}[1]{\href{http://arxiv.org/abs/#1}{\tt arXiv:\nolinkurl{#1}}}

\numberwithin{equation}{section}
\usepackage[linktocpage=true]{hyperref}
\hypersetup{colorlinks,linkcolor=blue,urlcolor=cyan,citecolor=blue}

\newtheorem{innercustomthm}{{\bf Theorem}}
\newenvironment{customthm}[1]
{\renewcommand\theinnercustomthm{#1}\innercustomthm}
{\endinnercustomthm}
\def \haH{\widehat{H}}
\def \haT{\widehat{\Theta}}

\def \Ui {\mathbf{U}^\imath}
\def \bt{\mathbf t}
\def \bn{\mathbf n}
\def \bh{\mathbf h}
\def \U{\mathbf{U}}
\def \Lg{L\fg}
\def \II{\mathbb{I}_0}
\def \haTh{\widehat{\Theta}}
\def \tK{K}
\def\SS{\mathbb{S}}
\def \RR{\mathbb{R}}
\def \bla{\boldsymbol{\lambda}}
\def \blx{x}
\def \G{\mathbb G}
\def \PL{\mathbb{P}^1_{\bfk}}
\def \scrM{\mathscr M}
\def \scrf{\mathscr F}
\def \scrt{\mathscr T}
\def \TH{\Theta}
\newcommand{\tUD}{{}^{\text{Dr}}\tU}

\def \haB{\widehat{B}}
\def \haX{\widehat{x}}

\allowdisplaybreaks
\begin{document}
	\input diagxy
	\xyoption{all}
	\newcommand{\ad}{\operatorname{ad}\nolimits}
	
	\newcommand{\iadd}{\operatorname{iadd}\nolimits}
	\newcommand{\Sym}{\operatorname{Sym}\nolimits}
	
	\renewcommand{\mod}{\operatorname{mod}\nolimits}
	\newcommand{\fproj}{\operatorname{f.proj}\nolimits}
	\newcommand{\Fac}{\operatorname{Fac}\nolimits}
	\newcommand{\aut}{\operatorname{Aut}\nolimits}
	\newcommand{\Cone}{\operatorname{cone}\nolimits}
	
	\newcommand{\proj}{\operatorname{proj}\nolimits}
	\newcommand{\inj}{\operatorname{inj.}\nolimits}
	\newcommand{\rad}{\operatorname{rad}\nolimits}
	\newcommand{\Span}{\operatorname{Span}\nolimits}
	\newcommand{\soc}{\operatorname{soc}\nolimits}
	\newcommand{\ind}{\operatorname{inj.dim}\nolimits}
	\newcommand{\Ginj}{\operatorname{Ginj}\nolimits}
	\newcommand{\res}{\operatorname{res}\nolimits}
	\newcommand{\np}{\operatorname{np}\nolimits}
	\newcommand{\nat}{\operatorname{nat.}\nolimits}
	\newcommand{\Aut}{\operatorname{Aut}\nolimits}
	\newcommand{\Char}{\operatorname{char}\nolimits}
	
	\newcommand{\Mod}{\operatorname{Mod}\nolimits}
	\newcommand{\R}{\operatorname{\mathbb R}\nolimits}
	\newcommand{\End}{\operatorname{End}\nolimits}
	\newcommand{\lf}{\operatorname{l.f.}\nolimits}
	\newcommand{\Iso}{\operatorname{Iso}\nolimits}
	
	\newcommand{\colim}{\operatorname{colim}\nolimits}
	\newcommand{\gldim}{\operatorname{gl.dim}\nolimits}
	\newcommand{\cone}{\operatorname{cone}\nolimits}
	\newcommand{\rep}{\operatorname{rep}\nolimits}
	\newcommand{\Ext}{\operatorname{Ext}\nolimits}
	\newcommand{\Tor}{\operatorname{Tor}\nolimits}
	\newcommand{\Hom}{\operatorname{Hom}\nolimits}
	\newcommand{\Top}{\operatorname{top}\nolimits}
	\newcommand{\Coker}{\operatorname{Coker}\nolimits}
	\newcommand{\thick}{\operatorname{thick}\nolimits}
	\newcommand{\rank}{\operatorname{rank}\nolimits}
	\newcommand{\Gproj}{\operatorname{Gproj}\nolimits}
	\newcommand{\Len}{\operatorname{Length}\nolimits}
	\newcommand{\RHom}{\operatorname{RHom}\nolimits}
	\renewcommand{\deg}{\operatorname{deg}\nolimits}
	\renewcommand{\Im}{\operatorname{Im}\nolimits}
	\newcommand{\Ker}{\operatorname{Ker}\nolimits}
	\newcommand{\Coh}{\operatorname{coh}\nolimits}
	\newcommand{\Id}{\operatorname{Id}\nolimits}
	\newcommand{\coh}{\operatorname{coh}\nolimits}
	\newcommand{\CM}{\operatorname{CM}\nolimits}
	\newcommand{\sgn}{\operatorname{sgn}\nolimits}
	\newcommand{\tMH}{\operatorname{\cs\cd\widetilde{\ch}(\cc_{\varrho}(\ca))}\nolimits}
	\newcommand{\utMH}{\operatorname{\cs\cd\ch(\cc_{\varrho}(\ca))}\nolimits}
	\newcommand{\For}{\operatorname{{\bf F}or}\nolimits}
	\newcommand{\coker}{\operatorname{Coker}\nolimits}
	\renewcommand{\dim}{\operatorname{dim}\nolimits}
	\newcommand{\rankv}{\operatorname{\underline{rank}}\nolimits}
	\newcommand{\dimv}{{\operatorname{\underline{dim}}\nolimits}}
	\newcommand{\diag}{{\operatorname{diag}\nolimits}}
	\newcommand{\tor}{{\operatorname{tor}\nolimits}}
	\renewcommand{\Vec}{{\operatorname{Vec}\nolimits}}
	\newcommand{\pd}{\operatorname{proj.dim}\nolimits}
	\newcommand{\gr}{\operatorname{gr}\nolimits}
	\newcommand{\id}{\operatorname{id}\nolimits}
	\newcommand{\Cha}{\operatorname{char}\nolimits}
	\newcommand{\swa}{\operatorname{swap}\nolimits}
	\def \ov{\overline}
	\def \cI{\mathcal{I}}
	\def \iH{{}^{\imath}\widetilde{\ch}}
	\newcommand{\pdim}{\operatorname{proj.dim}\nolimits}
	\newcommand{\idim}{\operatorname{inj.dim}\nolimits}
	\newcommand{\Gd}{\operatorname{G.dim}\nolimits}
	\newcommand{\Ind}{\operatorname{Ind}\nolimits}
	\newcommand{\add}{\operatorname{add}\nolimits}
	\newcommand{\pr}{\operatorname{pr}\nolimits}
	\newcommand{\oR}{\operatorname{R}\nolimits}
	\newcommand{\oL}{\operatorname{L}\nolimits}
	\newcommand{\Perf}{{\mathfrak Perf}}
	\newcommand{\cc}{{\mathcal C}}
	\newcommand{\gc}{{\mathcal GC}}
	\newcommand{\ce}{{\mathcal E}}
	\newcommand{\cs}{{\mathcal S}}
	\newcommand{\cf}{{\mathcal F}}
	\newcommand{\cx}{{\mathcal X}}
	\newcommand{\cy}{{\mathcal Y}}
	\newcommand{\ct}{{\mathcal T}}
	\newcommand{\cu}{{\mathcal U}}
	\newcommand{\cv}{{\mathcal V}}
	\newcommand{\cn}{{\mathcal N}}
	\newcommand{\mcr}{{\mathcal R}}
	\newcommand{\ch}{{\mathcal H}}
	\newcommand{\ca}{{\mathcal A}}
	\newcommand{\cb}{{\mathcal B}}
	\newcommand{\ci}{{\mathcal I}}
	\newcommand{\cj}{{\mathcal J}}
	\newcommand{\cm}{{\mathcal M}}
	\newcommand{\cp}{{\mathcal P}}
	\newcommand{\cg}{{\mathcal G}}
	\newcommand{\cw}{{\mathcal W}}
	\newcommand{\co}{{\mathcal O}}
	\newcommand{\cq}{{\mathcal Q}}
	\newcommand{\cd}{{\mathcal D}}
	\newcommand{\ck}{{\mathcal K}}
	\newcommand{\cz}{{\mathcal Z}}
	
	\def \bl{\mathbf l}
	\def \bs{\mathbf s}
	
	\newcommand{\calr}{{\mathcal R}}
	\newcommand{\ol}{\overline}
	\renewcommand{\ul}{\underline}
	\newcommand{\st}{[1]}
	\newcommand{\ow}{\widetilde}
	\renewcommand{\P}{\mathbf{P}}
	\newcommand{\pic}{\operatorname{Pic}\nolimits}
	\newcommand{\Spec}{\operatorname{Spec}\nolimits}
	\newtheorem{theorem}{Theorem}[section]
	\newtheorem{acknowledgement}[theorem]{Acknowledgement}
	\newtheorem{algorithm}[theorem]{Algorithm}
	\newtheorem{axiom}[theorem]{Axiom}
	\newtheorem{case}[theorem]{Case}
	\newtheorem{claim}[theorem]{Claim}
	\newtheorem{conclusion}[theorem]{Conclusion}
	\newtheorem{condition}[theorem]{Condition}
	\newtheorem{conjecture}[theorem]{Conjecture}
	\newtheorem{construction}[theorem]{Construction}
	\newtheorem{corollary}[theorem]{Corollary}
	\newtheorem{criterion}[theorem]{Criterion}
	\newtheorem{definition}[theorem]{Definition}
	\newtheorem{example}[theorem]{Example}
	\newtheorem{exercise}[theorem]{Exercise}
	\newtheorem{lemma}[theorem]{Lemma}
	\newtheorem{notation}[theorem]{Notation}
	\newtheorem{problem}[theorem]{Problem}
	\newtheorem{proposition}[theorem]{Proposition}
	\newtheorem{remark}[theorem]{Remark}
	\newtheorem{solution}[theorem]{Solution}
	\newtheorem{summary}[theorem]{Summary}
	\newtheorem*{thm}{Theorem}
	
	\newcommand{\qbinom}[2]{\begin{bmatrix} #1\\#2 \end{bmatrix} }
	
	\def \LaC{\Lambda_{\texttt{Can}}}
	\newcommand{\tUiD}{{}^{\text{Dr}}\tUi}
	\def \nua{a}
	\def \BF{\digamma}
	
	\newtheorem{innercustomprop}{{\bf Proposition}}
	\newenvironment{customprop}[1]
	{\renewcommand\theinnercustomprop{#1}\innercustomprop}
	{\endinnercustomthm}
	\newcommand{\UU}{{\mathbf U}\otimes {\mathbf U}}
	\newcommand{\UUi}{(\UU)^\imath}
	\newcommand{\tUU}{{\tU}\otimes {\tU}}
	\newcommand{\tUUi}{(\tUU)^\imath}
	\newcommand{\tUi}{\widetilde{{\mathbf U}}^\imath}
	\newcommand{\tU}{\widetilde{{\mathbf U}}}
	
	\def \haP{\widehat{P}}
	\def \cR{\mathcal R}
	\def \bfk{\mathbf k}
	\def \bp{{\mathbf p}}
	\def \bA{{\mathbf A}}
	\def \bL{{\mathbf L}}
	\def \bF{{\mathbf F}}
	\def \bh{{\mathbf h}}
	\def \bS{{\mathbf S}}
	\def \bC{{\mathbf C}}
	\def \bU{{\mathbf U}}
	\def \bm{{\mathbf m}}
	
	\def \fp{{\mathfrak p}}
	\def \fg{\mathfrak{g}}
	\def \sqq{\mathbf{v}}
	\def \scrM{{\mathscr M}}
	\def \scrf{{\mathscr F}}
	\def \scrt{{\mathscr T}}
	\def \fpr{\mathcal{P}^{<\infty}}
	\def \cp{{\mathcal{P}}}
	\def \Z{{\Bbb Z}}
	\def \F{{\Bbb F}}
	\def \D{{\Bbb D}}
	\def \C{{\Bbb C}}
	\def \N{{\Bbb N}}
	\def \Q{{\Bbb Q}}
	\def \G{{\Bbb G}}
	\def \P{{\Bbb P}}
	\def \K{{\Bbb K}}
	\def \E{{\Bbb E}}
	\def \I{{\Bbb I}}
	\def \X{{\Bbb X}}
	\def \A{{\Bbb A}}
	\def \L{{\Bbb L}}
	\def \BH{{\Bbb H}}
	\def \T{{\Bbb T}}
	\def \de{{\delta}}
	\newcommand{\UD}{{}^{\text{Dr}}\U}
	\def \TT{\mathbf T}
	\newcommand{\browntext}[1]{\textcolor{brown}{#1}}
	\newcommand{\greentext}[1]{\textcolor{green}{#1}}
	\newcommand{\redtext}[1]{\textcolor{red}{#1}}
	\newcommand{\bluetext}[1]{\textcolor{blue}{#1}}
	\newcommand{\brown}[1]{\browntext{ #1}}
	\newcommand{\green}[1]{\greentext{ #1}}
	\newcommand{\red}[1]{\redtext{ #1}}
	\newcommand{\blue}[1]{\bluetext{ #1}}
	
	\title[$\imath$Hall algebras of weighted projective lines and quantum symmetric pairs III]
	{$\imath$Hall algebras of weighted projective lines and quantum symmetric pairs III: quasi-split type}
	
	\author[Ming Lu]{Ming Lu}
	\address{Department of Mathematics, Sichuan University, Chengdu 610064, P.R.China}
	\email{luming@scu.edu.cn}

	\author[Shiquan Ruan]{Shiquan Ruan}
	\address{ School of Mathematical Sciences,
		Xiamen University, Xiamen 361005, P.R.China}
	\email{sqruan@xmu.edu.cn}
	
	
	\begin{abstract}
		From a category $\mathcal{A}$ with an involution $\varrho$, we introduce $\varrho$-complexes, which are  a generalization of (bounded) complexes, periodic complexes and modules of $\imath$quiver algebras.  The homological properties of the category $\mathcal{C}_\varrho(\mathcal{A})$ of $\varrho$-complexes are given to make the machinery of semi-derived Ringel-Hall algebras applicable. The $\imath$Hall algebra of the weighted projective line $\mathbb{X}$ is the twisted semi-derived Ringel-Hall algebra of $\mathcal{C}_\varrho({\rm coh}(\mathbb{X}))$, where $\varrho$ is an involution of ${\rm coh}(\mathbb{X})$. This $\imath$Hall algebra is used to realize the quasi-split $\imath$quantum loop algebra, which is a generalization of the $\imath$quantum group arising from the quantum symmetric pair of quasi-split affine type ADE in its Drinfeld type presentation. 
	\end{abstract}
	
	\maketitle
	\setcounter{tocdepth}{1}
	
	\tableofcontents
	
	\section{Introduction}
	
	\subsection{Hall algebras and quantum groups}
	
	Ringel \cite{Rin90} used the Hall algebra associated with a Dynkin quiver $Q$ over a finite field $\bfk=\mathbb F_q$ to  realize the  positive part $\U^+ =\U^+_v(\fg)$ of the Drinfeld-Jimbo quantum group $\U$; see Green \cite{Gr95} for an extension to acyclic quivers. 
	Building on these works, Bridgeland \cite{Br} in 2013 used a Hall algebra of $2$-periodic complexes (also called $\Z_2$-graded complexes) to realize
	the entire Drinfeld double $\widetilde \U$, a variant of $\U$ with the Cartan subalgebra doubled. 
	The semi-derived Ringel-Hall algebras formulated in \cite{LP16} is a further generalization of Bridgeland's construction to arbitrary hereditary abelian categories; compare with the semi-derived Hall algebras defined in \cite{Gor18}.

	
	There has been a current realization (also called Drinfeld type presentation) of the affine quantum groups formulated by Drinfeld \cite{Dr88, Be94, Da15}, which plays a fundamental role in  (algebraic and geometric) representation theory and mathematical physics. The quantum loop algebra $\U_v(\Lg)$ was defined as a generalization of the Drinfeld's current realization of the quantum affine algebra of any Kac-Moody algebra $\fg$.
	
	The Hall algebra of a weighted projective line was developed in \cite{Sch04} to realize the ``half part" of a quantum loop algebra, building on the realization of the  current half of quantum affine $\mathfrak{sl}_2$ via 
	the Hall algebra of the projective line \cite{Ka97,BKa01}. This realization was then upgraded to the whole quantum loop algebra via the Drinfeld double technique  \cite{DJX12, BS12, BS13} and the semi-derived Ringel-Hall algebra of the weighted projective line \cite{LP16}.
	

	\subsection{$\imath$Hall algebras and $\imath$quantum groups}
	
	The $\imath$quantum groups $\Ui$ (namely, quantum symmetric pair coideal subalgebras) arising from quantum symmetric pairs $(\U,\Ui)$ introduced by Letzter \cite{Le99, Let02} (also cf. \cite{Ko14}) can be viewed as a vast generalization of Drinfeld-Jimbo quantum groups; see the survey \cite{W22} and references therein. Just as the quantum groups $\U$ are associated to Dynkin diagrams, the $\imath$quantum groups $\Ui$ are associated to Satake diagrams (a bicolored diagram with a diagram involution $\tau$). The quasi-split $\imath$quantum groups, which correspond to Satake diagrams with no black nodes, already form a rich family of algebras. In this case, the Satake diagrams are just Dynkin diagrams with diagram involutions $\tau$ (possibly $\tau=\Id$, and called split type for this special case).

	In the framework of semi-derived Ringel-Hall algebras \cite{LP16},  Lu-Wang \cite{LW19a,LW20a} have developed $\imath$Hall algebras of $\imath$quivers to
	realize the universal quasi-split $\imath$quantum groups $\tUi$ of Kac-Moody type. The universal $\imath$quantum group $\tUi$ is a coideal subalgebra of the Drinfeld double quantum group $\tU$, and the
	$\imath$quantum group $\Ui$ \cite{Le99,Let02,Ko14} (with parameters) is obtained by a central reduction of $\tUi$. In this paper, we shall work with universal quantum symmetric pairs $(\tU,\tUi)$ following \cite{LW19a}, as this allows us to formulate the braid group action conceptually \cite{LW19b,LW22,WZ22}.
	
	In \cite{LWZ23}, a Drinfeld type presentation of affine $\imath$quantum groups $\tUi$ of quasi-split ADE type is obtained, building on the Drinfeld type presentation of $\tUi$ of split ADE type given in \cite{LW21b}; see an extension to split BCFG types in \cite{Z22}.   
	
	In \cite{LRW20a}, the authors used the $\imath$Hall algebra of the projective line to give a geometric realization of the $q$-Onsager algebra (i.e., $\imath$quantum group of split affine $\mathfrak{sl}_2$) in its Drinfeld type presentation. 
	In \cite{LR21,LR23}, the $\imath$Hall algebras of weighted projective lines are used to give a geometric realization of split $\imath$quantum loop algebras associated with star-shaped graphs (a generalization of $\imath$quantum groups arising from quantum symmetric pairs of split affine ADE type \cite{LW21b}). 
	(Here the $\imath$Hall algebra refers to the
	twisted semi-derived Ringel-Hall algebra of the category of 1-periodic complexes
	of coherent sheaves.)

	
	

	\subsection{Goal} 
	This paper is a sequel to \cite{LRW20a,LR21}. We use $\imath$Hall algebras of weighted projective lines to give a geometric realization of the quasi-split $\imath$quantum loop algebras of star-shaped graphs (a generalization of $\imath$quantum groups of quasi-split affine ADE type in Drinfeld type presentation). The key point is to find  suitable categories and use suitable Hall algebras.
	We introduce the $\varrho$-complexes for a category $\ca$ with automorphism $\varrho$, generalizing (bounded) complexes, $m$-periodic complexes and representations of $\imath$quiver algebras \cite{LW19a}. The $\imath$Hall algebra here is the (twisted) semi-derived Ringel-Hall algebra of the category $\cc_\varrho(\ca)$ of $\varrho$-complexes  for a hereditary abelian category $\ca$ with involution $\varrho$, by using the general machinery of semi-derived Ringel-Hall algebras in \cite[Appendix A]{LW19a}. Then we consider the category $\ca=\coh(\X_\bfk)$ of coherent sheaves over a weighted projective line $\X_\bfk$ appended with an involution $\varrho$. Then we use the $\imath$Hall algebra $\iH(\X_\bfk,\varrho)$ to realize the quasi-split $\imath$quantum loop algebra $\tUiD(\Lg)$.

	\subsection{Main results}
	
	\subsubsection{$\imath$-Categories and $\imath$Hall algebras }
	
	Let $\bfk=\F_q$, the field of $q$ elements. Let $\ca$ be a hereditary $\bfk$-linear abelian  category with an involution $\varrho:\ca\rightarrow \ca$. The $\varrho$-complex is a pair $(M,d)$ with $M\in\ca$ and $d:M\rightarrow \varrho(M)$ satisfying $\varrho(d)\circ d=0$. Let $\cc_\varrho(\ca)$ be the category of $\varrho$-complexes, which is called the $\imath$-category of $(\ca,\varrho)$. Typical examples are categories of $m$-periodic complexes for $m=1,2$, and categories of modules of $\imath$quiver algebras \cite{LW19a}.

	
	Let $\cd^b(\ca)$ be the derived category of $\ca$ with suspension functor $\Sigma$. The involution $\varrho$ induces a triangulated automorphism $\widehat{\varrho}: \cd^b(\ca)\rightarrow \cd^b(\ca)$. Let $\cd_\varrho(\ca)$ be the relative derived category of $\cc_{\varrho}(\ca)$; see \S\ref{subsec:ca-auto}.   
	
	\begin{customthm}{{\bf A}}
		[Theorems 
		\ref{prop: orbit}]
		Let $\ca$ be a hereditary abelian category with an involution $\varrho$.
		Then  $\cd_\varrho(\ca)\simeq \cd^b(\ca)/\Sigma\circ \widehat{\varrho}$. 
	\end{customthm}
	
	For the homological dimension of $\varrho$-complex $M^\bullet$, we have the following result, which is similar to $1$-Gorenstein algebras.
	
	\begin{customprop}{{\bf B}}
		[Proposition \ref{cor:projfinite}]
		For any $M^\bullet\in\cc_\varrho(\ca)$, the following statements are equivalent.
		\begin{enumerate}
			\item $\pd_{\cc_{\varrho}(\ca)} M^\bullet<\infty$,
			\item $\pd_{\cc_{\varrho}(\ca)} M^\bullet\leq1$,
			\item $\ind_{\cc_{\varrho}(\ca)} M^\bullet<\infty$,
			\item $\ind_{\cc_{\varrho}(\ca)} M^\bullet\leq1$,
			\item $M^\bullet$ is acyclic.
		\end{enumerate}
	\end{customprop}
	
	Using Proposition {\bf B}, for $K^\bullet,M^\bullet\in\cc_\varrho(\ca)$ with $K^\bullet$ acyclic,  we define the Euler form $\langle K^\bullet,M^\bullet\rangle$ and $\langle M^\bullet,K^\bullet\rangle$ in \eqref{left Euler form}--\eqref{right Euler form}. This Euler form is deeply related to the Euler form $\langle\cdot,\cdot\rangle$ of $\ca$ by using the restriction functor $\res:\cc_\varrho(\ca)\rightarrow \ca$. Recall from \eqref{eq:KX} the acyclic $\varrho$-complex $K_X$ for any $X\in\ca$.
	
	\begin{customthm}{{\bf C}}
		[Theorem \ref{prop:Euler}]
		For $K^\bullet=(K,d)$ and $M^\bullet \in\cc_{\varrho}(\ca)$ with $K^\bullet$ acyclic, we have 
		\begin{align}
			\label{eq:KMeuler1}
			\langle K^\bullet,M^\bullet\rangle=\langle K_{\Im(d)},M^\bullet\rangle=\langle \Im(d),\res(M^\bullet)\rangle ,
			\\
			\label{eq:MKeuler2}
			\langle M^\bullet,K^\bullet\rangle=\langle M^\bullet,K_{\Im(d)}\rangle =\langle \res(M^\bullet),\varrho(\Im(d))\rangle.
		\end{align}
		In particular, if $M^\bullet$ is also acyclic, then 
		\begin{align}
			\label{eq:MKeuler3}
			\langle K^\bullet,M^\bullet\rangle=\frac{1}{2}\langle \res (K^\bullet),\res (M^\bullet)\rangle.
		\end{align}
	\end{customthm}
	
	Let $\ch(\cc_\varrho(\ca))$ be the Ringel-Hall algebra of $\cc_\varrho(\ca)$. With the help of Theorem {\bf A,C} and Proposition {\bf B}, we can use the general machinery established in \cite[Appendix A]{LW19a} to define the semi-derived Ringel-Hall algebra $\cs\cd\ch(\cc_\varrho(\ca))$ to be the localization $(\ch(\cc_{\varrho}(\ca))/\cI)[\cs^{-1}]$, where $\cs$ is defined in \eqref{eq:Sca}, and the ideal $\cI$ is given in \eqref{eq:ideal}. 
	
	The $\imath$Hall algebra $\iH(\ca,\varrho)$ (also denoted by $\cs\cd\widetilde{\ch}(\cc_\varrho(\ca))$) is defined to be the semi-derived Ringel-Hall algebra with the multiplication twisted by the Euler form of $\ca$; see \eqref{eqn:twsited multiplication}. An $\imath$Hall basis of $\iH(\ca,\varrho)$ is given in Proposition \ref{prop:hallbasis}, and an $\imath$Hall multiplication formula can be found in Proposition \ref{prop:iHallmult}; cf. \cite{LW19a,LW20a,LinP}.

	\subsubsection{Involutions of weighted projective lines}
	
	Let $\X:=\X_{\bp,\ul{\bla}}$ be the weighted projective line over $\bfk$ of weight type $(\bp,\ul{\bla})$, where $\bp=(p_1,p_2,\dots,p_\bt)$, and $\ul{\bla}=\{\bla_1,\dots,\bla_\bt\}$ is a collection of distinguished closed points of degree one. Let $\coh(\X)$ be the category of coherent sheaves over $\X$, and $\co$ be its structure sheaf.

	Let $\Gamma=\T_{p_1,p_2,\dots,p_\bt}$ be the star-shaped graph displayed in  \eqref{star-shaped}. The automorphism group $\Aut(\X)$ of $\coh(\X)$ (fixing $\co$) is studied in \cite{LM}, which is a subgroup of $\Aut(\P_\bfk^1)$. In this paper, we focus on involutions of $\coh(\X)$.
	
	\begin{customprop}{{\bf D}}
		[Proposition \ref{prop:involution}]
		Let $\bp=(p_1,\dots,p_\bt)\in\Z_+^\bt$, and $\varrho$ be a nontrivial involution of the star-shaped graph $\Gamma$. Then there exists a lifting of  $\varrho$ to $\coh(\X)$ for some weighted projective line $\X$ if and only if the number of $\varrho$-fixed branches $\leq 2$. Moreover, any involution $\varrho$ of the star-shaped graph $\Gamma$ of finite and affine type can be lifted to $\coh(\X)$. 
	\end{customprop}
	
	

	\subsubsection{$\imath$Quantum loop algebras and $\imath$Hall algebras}
	
	Let $\Gamma=\T_{p_1,\dots,p_\bt}$ be a star-shaped graph with its vertex set $\II$ as in \eqref{star-shaped}, $\fg$ be its associated Lie algebra. Let $\varrho$ be an involution of $\Gamma$, and
	$\tUiD$ be its associated quasi-split $\imath$quantum loop algebra \cite{LWZ23} defined as in Definition \ref{def:iDRA1}. If $\Gamma$ is of ADE type, then $\tUiD$ is the Drinfeld type presentation of the $\imath$quantum group arising from quantum symmetric pairs of quasi-split affine ADE type; see \cite{LWZ23}. Recall that $\tUiD$ is generated by $B_{\star,l},  B_{[i,j],l}, H_{\star,m},H_{[i,j],m}, \Theta_{\star,m},\Theta_{[i,j],m}$ and central elements $\K_\star^{\pm1}, \K_{[i,j]}^{\pm}, C^{\pm1}$ for $l\in\Z,m>0$, and $1\leq i\leq \bt$, $1\leq j\leq p_i-1$, subject to \eqref{qsiDR0}--\eqref{qsiDR9}.

	Let $\X_\bfk$ be the weighted projective line over a finite field $\bfk=\F_q$ associated to $\Gamma$. Assume $\varrho$ is an involution of $\coh(\X)$.  Then we have the $\imath$Hall algebra $\iH(\X_\bfk,\varrho)$. 
	
	For the special case $\P_\bfk^1$, the projective line over $\bfk$, we consider the $q$-Onsager algebra $\tUi(\widehat{\mathfrak{sl}}_2)$ (i.e., the split $\imath$quantum group of $\widehat{\mathfrak{sl}}_2$). Its Drinfeld type presentation $\tUiD(\widehat{\mathfrak{sl}}_2)$ is generated by $\K_1^{\pm1}$, $C^{\pm1}$, $B_{1,r}$, $\Theta_m$, $H_m$ for $r\in\Z$, $m\geq1$. 
	Then the following result is a generalization of \cite[Theorem {\bf A}]{LRW20a}.
	
	\begin{customthm}{{\bf F}}
		[Propositions \ref{prop:P1-Onsager}, \ref{prop:HaH}]
		Let $\varrho$ be an involution of $\coh(\P_\bfk^1)$. Then there exists a $\Q(\sqq)$-algebra homomorphism
		\begin{align}
			\label{eq:phi}
			\Omega: \tUiD_{{\sqq}}(\widehat{\mathfrak{sl}}_2)\longrightarrow \iH(\P^1_\bfk,\varrho)
		\end{align}
		which sends, for all $r\in \Z$ and $m \ge 1$,
		\begin{align*}
			\K_1\mapsto [K_\co],  \quad
			C\mapsto [K_\de],  \quad
			B_{1,r} \mapsto -\frac{1}{q-1}[\co(r)],  \quad
			\Theta_{m}\mapsto  \haT_m,\quad H_m\mapsto \widehat{H}_m.
		\end{align*}
	\end{customthm}
	Here $\haT_m$ is defined to 
	\begin{align*}
		\haT_{m}= \frac{1}{(q-1)^2\sqq^{m-1}}\sum_{0\neq f:\co\rightarrow \co(m) } [\coker f],\qquad \forall m\geq1,
	\end{align*}
	and $\widehat{H}_m$ is given by Proposition \ref{prop:HaH}.

	Let $\Gamma$ be the star-shaped graph with an involution  $\varrho$. Let $\tUiD$ be the $\imath$quantum loop algebra of $(\Gamma,\varrho)$. Let 
	$$\I_\varrho:=\{\text{the chosen representatives of $\varrho$-orbits of nodes in $\Gamma$} \}.$$
	Assume that $\varrho$ can be lifted to an involution $\varrho$ of $\coh(\X)$ for a weighted projective line $\X$ associated to $\Gamma$; cf. Proposition {\bf D}. 
	Then we have the following theorem. 
	
	\begin{customthm}{{\bf G}}
		[Theorem \ref{thm:morphi}]
		There exists a $\Q(\sqq)$-algebra homomorphism
		\begin{align}
			\Omega: \tUiD_\sqq\longrightarrow \iH(\X_\bfk,\varrho),
		\end{align}
		which sends
		\begin{align}
			&\K_{\star}\mapsto [K_{\co}], \qquad \K_{[i,j]}\mapsto [K_{S_{ij}}], \qquad C\mapsto [K_\de];&
			\\
			&{B_{\star,l}\mapsto \frac{-1}{q-1}[\co(l\vec{c})]},\qquad\Theta_{\star,r} \mapsto {\widehat{\Theta}_{\star,r}}, \qquad H_{\star,r} \mapsto {\widehat{H}_{\star,r}};\\
			& \Theta_{[i,j],r}\mapsto\widehat{\Theta}_{[i,j],r}, \quad H_{[i,j],r}\mapsto \widehat{H}_{[i,j],r}
			,\qquad B_{[i,j],l}\mapsto \begin{cases}{\frac{-1}{q-1}}\haB_{[i,j],l}, \text{ if }[i,j]\in\I_\varrho,\\
				\frac{\sqq}{q-1}B_{[i,j],l},\text{ if }[i,j]\notin\I_\varrho;\end{cases}
		\end{align}
		for any $[i,j]\in\II-\{\star\}$, $l\in\Z$, $r>0$.
	\end{customthm}

	Let us explain the counterparts of the Drinfeld type generators in $\iH(\X_\bfk,\varrho)$. First, $\widehat{\Theta}_{\star,r}$ and $\haH_{\star,r}$ are defined by using the corresponding ones $\widehat{\Theta}_r$, $\haH_{r}$ defined in $\iH(\P^1_\bfk,\varrho)$ (see Theorem {\bf F}) via the natural embedding (see \eqref{the embedding functor F on algebra})
	$$F_{\X,\P^1}:\iH(\P^1_\bfk,\varrho) \longrightarrow \iH(\X_\bfk,\varrho).$$
	
	Second, 
	the subcategory $\scrt_{\bla_i}$ of $\coh(\X_\bfk)$, consisting of torsion sheaves supported at the distinguished point $\bla_i$, is isomorphic to $\rep^{\rm nil}_\bfk( C_{p_i})$, the category of finite-dimensional nilpotent representations of  the cyclic quiver $C_{p_i}$ with $p_i$ vertices.  For $1\leq i\leq \bt$, if the $i$-th branch of $\Gamma$ is fixed by $\varrho$ (i.e., $\varrho(\bla_i)=\bla_i$), then there is an algebra embedding $$\widetilde{\psi}_{C_{p_i}}: \tUi_\sqq(\widehat{\mathfrak{sl}}_{p_i})\longrightarrow \iH(\bfk C_{p_i},\Id)$$ by \cite[Theorem 9.6]{LW20a}. As a consequence, each branch of $\Gamma$  corresponds to a subalgebra of $\iH(\X_\bfk,\varrho)$ which is isomorphic to
	$\tUiD_v(\widehat{\mathfrak{sl}}_{p_i})$ by the composition $\Omega_{C_{p_i}}$ of $\widetilde{\psi}_{C_{p_i}}$ and the isomorphism of two presentations \cite{LW21b} (see \eqref{the map Psi A}) $$\Phi:\tUiD_v(\widehat{\mathfrak{sl}}_{p_i})\stackrel{\cong}{\longrightarrow} \tUi_v(\widehat{\mathfrak{sl}}_{p_i}).$$
	The equivalence $\scrt_{\bla_i}\simeq \rep^{\rm nil}_\bfk( C_{p_i})$ induces an algebra embedding (see \eqref{eq:embeddingx})
	$$\iota_i: \iH(\bfk C_{p_i},\Id)\longrightarrow\iH(\X_\bfk,\varrho).$$
	So we define $\haB_{[i,j],l},\widehat{\Theta}_{[i,j],r},\widehat{H}_{[i,j],r}$ to be the images of the Drinfeld generators of $\tUiD_v(\widehat{\mathfrak{sl}}_{p_i})$ under the composition $\iota_i\circ \Omega_{C_{p_i}}$; see \eqref{def:haBThH}.
	
	Otherwise, assume $\bla_{\varrho(i)}=\varrho(\bla_i)\neq\bla_i$. First, there is an algebra embedding \begin{align}
		\widetilde{\psi}_{C_{p_i}}:\tU_\sqq(\widehat{\mathfrak{sl}}_{p_i})&\longrightarrow \cs\cd\widetilde{\ch}_{\Z_2}( \bfk C_{p_i}),
	\end{align}
	where $\cs\cd\widetilde{\ch}_{\Z_2}( \bfk C_{p_i})$ is the twisted semi-derived Ringel-Hall algebra of $\cc_{\Z_2}(\rep^{\rm nil}_\bfk(C_{p_i}))$; see Proposition \ref{lem:UslnSDH}. Let $\tUD_v(\widehat{\mathfrak{sl}}_{p_i})$ be the Drinfeld presentation of $\tU_v(\widehat{\mathfrak{sl}}_{p_i})$, which is generated by $x_{i k}^{\pm}$, $h_{i l}$ and the invertible elements $K_i$, $K_i'$, $C$, $\widetilde{C}$ for $i\in\II$, $k\in\Z$, $l\in\Z\backslash\{0\}$; see \S\ref{sec:QLA}. Composing with the Drinfeld-Beck isomorphism $\Phi:\tUD_v(\widehat{\mathfrak{sl}}_{p_i})\rightarrow\tU_v(\widehat{\mathfrak{sl}}_{p_i})$ (see \cite{Be94}), we obtain
	\begin{align}
		\Omega_{C_n}:=\widetilde{\psi}_{C_n}\circ \Phi:\tUD_{\sqq}(\widehat{\mathfrak{sl}}_{p_i})\longrightarrow \cs\cd\widetilde{\ch}_{\Z_2}(\bfk C_n).
	\end{align}
	On the other hand, the restriction of $\varrho$ to  $\scrt_{\bla_i}\times \ct_{\bla_{\varrho i}}$ gives a subcategory $\cc_\varrho(\scrt_{\bla_i}\times \ct_{\bla_{\varrho i}})$ of $\cc_\varrho(\coh(\X))$, which is identified with $\cc_{\Z_2}(\rep^{\rm nil}_\bfk(C_{p_i}))$ by using the equivalence
	$\scrt_{\bla_i}\simeq \rep^{\rm nil}_\bfk( C_{p_i})$. This induces a morphism of $\imath$Hall algebras:
	\begin{align}
		\iota_i:\cs\cd\widetilde{\ch}_{\Z_2}(\bfk C_{p_i})\longrightarrow \iH(\X_\bfk,\varrho).
	\end{align}
	So we define $\haB_{[i,j],l},\widehat{\Theta}_{[i,j],r},\widehat{H}_{[i,j],r}$ to be the images of the Drinfeld generators of $\tUD_v(\widehat{\mathfrak{sl}}_{p_i})$ under the composition $\iota_i\circ \Omega_{C_{p_i}}$; see \eqref{def:haBXno}--\eqref{def:haphino}.
	
	To show that $\Omega: \tUiD_{\sqq} \rightarrow \iH(\X_\bfk,\varrho)$ is a homomorphism, we must verify  the Drinfeld type relations \eqref{qsiDR0}--\eqref{qsiDR9} for the $\imath$quantum loop algebra $\tUiD$.
	The counterparts in $\iH(\X_\bfk,\varrho)$ of the relations at the $\star$ point  follow from Theorem {\bf F}, and the counterparts in $\iH(\X_\bfk,\varrho)$ of the relations involving all the vertices $[i,j]$ follow from the above definitions and the fact any two torsion sheaves supported at distinct points have zero Hom and $\Ext^1$-spaces.
	
	In order to check the remaining Drinfeld type relations in $\iH(\X_\bfk,\varrho)$, the key part is to verify the relations between $\star$ and $[i,1]$ for $1\leq i\leq \bt$, especially the ones between $\star$ and $[i,1]$, which are given in \S\ref{subsec:relationsstarjneq1}--\S\ref{sec:Relationsstari1 II}. 
	
	\subsection{Comparison with previous works}
	
	In \cite{LRW20a,LR21}, the $\imath$Hall algebra is the twisted semi-derived Ringel-Hall algebra of the category of $1$-periodic complexes. The category $\cc_{\Z_1}(\ca)$ of $1$-periodic complexes of $\ca$ is well known, and studied by many experts before (see e.g., \cite{RZ17,St17} and the references therein). The semi-derived Ringel-Hall algebra of this category $\cc_{\Z_1}(\ca)$ can be established easily by using \cite{LP16} (see also \cite{LinP}). However, it is much difficult to verify the Drinfeld type relations of the $\imath$quantum loop algerbra of split type in the $\imath$Hall algebra of the weighted projective line in \cite{LR21}.
	
	Compared with \cite{LRW20a,LR21}, we introduce $\varrho$-complexes for a category $\ca$ with involution $\varrho$. This kind of complexes is new, which is a nontrivial generalization of $1$-periodic complexes, as $1$-periodic complexes are $\varrho$-complexes with $\varrho=\Id$. In order to make the general machinery of semi-derived Ringel-Hall algebras applicable for the category of $\varrho$-complexes, we study its homological properties carefully in Section \ref{sec:cat-inv}, and many proofs in \cite{LP16,LW19a} can not work in this general setting (see e.g. Propositions \ref{prop:Extconincide}, \ref{cor:projfinite},  Theorems \ref{prop: orbit}, \ref{prop:Euler}).
	
	In order to realize the quasi-split $\imath$quantum loop algebras, 
	we should generalize some results in \cite{LRW20a,LR21,LRW21}; see e.g. Propositions \ref{prop:realroot1}, \ref{prop:imageroot}, \ref{prop:P1-Onsager}, \ref{prop:HaH}. Then the route of verification of Drinfeld type relations  of $\imath$quantum loop algebra in the $\imath$Hall algebra is similar to \cite{LR21}, with some new computations and some proofs there modified.
	
	\subsection{Further works}
	
	The category of $\varrho$-complexes introduced in this paper is a new kind of categories, we shall study it further. 
	We shall also study the homological properties of $\varrho$-complexes for arbitrary abelian categories (or even exact categories) which are not hereditary.
	
	We expect the morphism $\Omega: \tUiD\rightarrow \iH(\X_\bfk,\varrho)$ to be injective for arbitrary Kac-Moody algebra $\fg$; cf. \cite{Sch04,DJX12, LR21, LR23}. In fact, following the arguments of \cite{LR23}, we can prove that 
	$\Omega$ is injective for $\fg$ of finite or affine type.

	In a sequel, we shall describe the composition subalgebra of the $\imath$Hall algebra of a weighted projective line, and give a PBW basis for the $\imath$quantum loop algebra via coherent sheaves.  
	
	The semi-derived Ringel-Hall algebras of Jordan quiver and cyclic quivers should be considered in depth, for example, their connections to Laurent symmetric functions and Schur-Weyl duality; cf. \cite{DDF12}.
	For a cyclic quiver $C_n$, its composition  algebra of $\cs\cd\ch_{\Z_2}(\bfk C_n)$  is used to realize the quantum group $\tU$ of affine ${\mathfrak{sl}}_{n}$, while the entire $\cs\cd\ch_{\Z_2}(\bfk C_n)$ shall be used to realize the quantum group $\tU$ of affine ${\mathfrak{gl}}_{n}$; cf. \cite{Sch02,DDF12}. 
	
	The $\imath$Hall algebra $\iH(\bfk Q,\varrho)$ of the quiver algebra $\bfk Q$ of affine type A was known earlier to realize the same algebra in its Serre presentation.
	The Geigle-Lenzing's derived equivalence induces an isomorphism of these two $\imath$Hall algebras $\iH(\X_\bfk,\varrho)\stackrel{\cong}{\longrightarrow}\iH(\bfk Q,\varrho)$, explaining the isomorphism of the quasi-split $\imath$quantum group of affine type A
	under the Serre and Drinfeld type presentations \cite{LWZ23}.
	
	It is also interesting to study the $\imath$Hall algebras of higher genus curves, in particular, of elliptic curves.

	\subsection{Organization}
	
	This paper is organized as follows. We 
	define $\varrho$-complexes and study their homological properties  in Section \ref{sec:cat-inv}.
	In Section \ref{sec:Semi-derived}, we 
	construct $\imath$Hall algebras for the category of $\varrho$-complexes, and in Section \ref{sec:WPL} we review the weighted projective lines and coherent sheaves. 
	In Section \ref{sec:QG-iQG}, we review the materials on quantum groups, quantum loop algebras,  $\imath$quantum group and $\imath$quantum loop algebras. 
	
	In Section \ref{sec:cyclic}, we use the $\imath$Hall algebras $\cs\cd\widetilde{\ch}_{\Z_2}(\bfk C_n)$ and $\iH(\bfk C_n,\Id)$ of cyclic quiver $C_n$ to realize the Drinifeld presentation of quantum groups $\tU_v(\widehat{\mathfrak{sl}}_n)$ and $\imath$quantum groups  $\tUi_v(\widehat{\mathfrak{sl}}_n)$. Some root vectors are also described in Hall algebras.
	The map $\Omega: \tUiD_\sqq\rightarrow \iH(\X_\bfk,\varrho)$ in Theorem~{\bf G} is formulated in Section~\ref{sec:hom}, and the proof of $\Omega$  being an algebra homomorphism is given in \S\ref{sec:Relationtube}-- \S\ref{sec:Relationsstari1 II}.

	\subsection{Acknowledgments}
	ML is partially supported by the National Natural Science Foundation of China (No. 12171333, 12261131498). 
	SR is partially supported by the National Natural Science Foundation of China (No. 11801473) and the Fundamental Research Funds for Central Universities of China (No. 20720180006).
	
	\section{Categories with involutions}
	\label{sec:cat-inv}
	
	Throughout this paper we always assume $\bfk=\F_q$ is a finite field with $q$ elements and its characteristic $\Char \bfk\neq2$. 
	
	In this section, we shall introduce $\imath$-categories from categories with involutions. In fact, the general framework of
	categories with automorphisms is introduced, as a generalization of categories of (bounded) complexes and $m$-periodic complexes.

	\subsection{$m$-periodic complexes}
	
	We assume that $\ca$  is an exact category. For the basics of exact categories, we refer to \cite{Buh}. For $m\geq1$, let $\Z_m=\Z/m\Z=\{0,1,\cdots,m-1\}$.
	Let $\cc_{\Z_m}(\ca)$ be the exact category of $m$-periodic complexes (also called $\Z_m$-graded complexes) over $\ca$. Namely, an object $M^\bullet=(M^i,d^i=d^i_{M^\bullet})_{i\in\Z_m}$ of this category consists of $M^i\in\ca$ and morphisms $d^i:M^i\rightarrow M^{i+1}$ with $d^{i+1}d^i=0$ for any $i\in\Z_m$. For convenience, we express $M^\bullet$ to be the following:
	$$ \xymatrix{ M^0\ar[r]^{d^0} & M^1\ar[r]^{d^1}&\cdots \ar[r]^{d^{m-2}} &M^{m-1}\ar[r]^{d^{m-1}} &M^0.}$$ 
	
	Let $M^\bullet=(M^i,d^i)$, $N^\bullet=(N^i,e^i)$. 
	A {morphism} $f^\bullet=(f^i):M^\bullet\rightarrow N^\bullet$ is a family of morphisms $f^i:M^i\rightarrow N^i$ such that $e^{i}f^i=f^{i+1}d^i$ for any $i\in\Z_m$. Two morphisms
	$f^\bullet,g^\bullet:M^\bullet\rightarrow N^\bullet$ are called homotopic equivalent if there exist $h^i:M^i\rightarrow N^{i-1}$ such that
	$$f^i-g^i= h^{i+1}d^i+e^{i-1}h^i$$ for all $i\in\Z_m$.   
	Then the $m$-periodic homotopy category $\ck_{\Z_m}(\ca)$ is defined.
	
	A morphism $f^\bullet=(f^i):M^\bullet\rightarrow N^\bullet$ is called a quasi-isomorphism if the induced morphisms 
	$H^i(f^\bullet):H^i(M^\bullet)\rightarrow H^i(N^\bullet)$ are isomorphisms for all $i\in\Z_m$. The $m$-periodic derived category $\cd_{\Z_m}(\ca)$ is the localization of $\ck_{\Z_m}(\ca)$ with respect to all quasi-isomorphisms.

	The shift functor $\Sigma$ on complexes is an (exact) functor 
	$\Sigma:\cc_{\Z_m}(\ca)\rightarrow \cc_{\Z_m}(\ca)$	defined by
	$(\Sigma M^\bullet)^i=M^{i+1}$, and $d^i_{\Sigma M^\bullet}=-d^{i+1}$ for any $i\in\Z_m$.
	Then $\ck_{\Z_m}(\ca)$ and $\cd_{\Z_m}(\ca)$ are triangulated categories with suspension functors $\Sigma$; cf. \cite{PX97,Ke05}.

	
	We denote by $\cc^b(\ca)$ the category of bounded complexes over $\ca$, and by $\ck^b(\ca)$, $\cd^b(\ca)$ the bounded homotopy category and the bounded derived category respectively.

	%
	%
	\subsection{Categories with automorphisms}
	\label{subsec:ca-auto}
	Let $\ca$ be an exact category, and $\varrho: \ca\rightarrow \ca$ an automorphism.
	We call the pair $(\ca,\varrho)$ a category with automorphism. A $\varrho$-complex $M^\bullet$ is a pair $(M,d)$ such that $M\in\ca$, and $d:M\rightarrow \varrho(M)$ is a morphism satisfying $\varrho(d)d=0$ ($d$ is called a differential). A morphism of $\varrho$-complexes $f: (M,d)\rightarrow (N,e)$ is  a morphism $f:M\rightarrow N$ in $\ca$ such that $\varrho(f) d=e f$, i.e., the following diagram commutes:
	\[\xymatrix{ M\ar[r]^d \ar[d]^f&\varrho M \ar[d]^{\varrho (f)}
		\\
		N\ar[r]^e& \varrho N.}\]
	
	Denote by  $\cc_\varrho(\ca)$ the category of $\varrho$-complexes. 

	
	\begin{example}
		If $\varrho=\Id$, then $\varrho$-complexes are just $1$-periodic complexes, and $\cc_\varrho(\ca)$ is the category $\cc_{\Z_1}(\ca)$ of $1$-periodic complexes.
	\end{example}
	
	A $\varrho$-complex $M^\bullet=(M,0)$ is called a stalk complex, which is just denoted by $M$. Then $\ca$ is a full subcategory of $\cc_{\varrho}(\ca)$. The homology group of $M^\bullet=(M,d)$ is $H(M^\bullet):=\Ker (d)/\Im (\varrho^{-1} (d))\in\ca$. We call $(M,d)$ acyclic if its homology group is zero. We denote by $\cc_{\varrho,ac}(\ca)$ the subcategory formed by acyclic $\varrho$-complexes.
	
	A morphism $f:(M,d)\rightarrow (N,e)$ is called null-homotopic if there  is some $h:\varrho M\rightarrow N$ such that $f=hd+\varrho^{-1}(e)\circ \varrho^{-1}(h)$, and called quasi-isomorphic if $f$ induces an isomorphism on homology groups. Similar to $m$-periodic complxes, we can define the relative homotopy category $\ck_\varrho(\ca)$ and the relative derived category $\cd_\varrho(\ca)$, and one can prove that $\ck_\varrho(\ca)$ and $\cd_\varrho(\ca)$ are triangulated categories with the suspension functor given by 
	\begin{align}
		\Sigma: (M,d)\mapsto (\varrho(M),-\varrho(d));
	\end{align}
	cf. \cite{Ha88,PX97,St17}.
	
	Let
	\begin{align}
		\res:\cc_{\varrho}(\ca)\longrightarrow \ca
		, \qquad (M,d)\mapsto M	
	\end{align}
	be	the forgetful functor.

	\begin{example}
		\label{ex:2-periodic}
		Let $\ca$ be a category, and $m\geq1$. Let $\cb=\ca^{\times m}$ be the product of $m$ copies of $\ca$. We see that  $\cb$ admits a natural antomorphism $\varrho$ defined by $$
		\varrho((M^0,M^1,\cdots,M^{m-1}))=(M^1,M^2,\cdots, M^{m-1},M^0).$$ 
		Then $\cc_{\varrho}(\cb)$ is isomorphic to the category of $m$-periodic complexes over $\ca$ with the automorphic functor given by
		\begin{align}
			F:	\cc_{\varrho}(\cb)&\longrightarrow \cc_{\Z_m}(\ca),\\
			\big((M^0,\dots,M^{m-1}),(d^0,\cdots,d^{m-1})\big)&\mapsto (M^i,d^i)_{i\in\Z_m}.
		\end{align}
	\end{example}
	
	Obviously, if $\ca$ is an abelian category with an automorphism $\varrho$, then $\cc_\varrho(\ca)$ is also an abelian category.

	In the following, we assume that $\ca$ is a $\bfk$-linear hereditary abelian category, with Hom-spaces finite-dimensional, and an automorphism $\varrho$. The $\varrho$-complexes are also called $\imath$-complexes, and the category $\cc_\varrho(\ca)$ is called $\imath$-category if $\varrho$ is an involution, i.e., $\varrho^2=\Id$.

	\subsection{Another description of $\imath$-categories}
	Unless otherwise specified, we always assume that $\varrho$ is an involution in the following.
	We describe $\imath$-complexes by using the language of invariant subcategories.
	
	Recall that $\cc_{\Z_2}(\ca)$ is the category of $2$-periodic complexes $(\xymatrix{ M^0 \ar@<0.5ex>[r]^{d^0}& M^1 \ar@<0.5ex>[l]^{d^1}  })$. 
	Let $(\ca,\varrho)$ be a category with involution. Then $\varrho$ induces an involution $\varrho^\sharp$ on $\cc_{\Z_2}(\ca)$ by mapping
	$$ (\xymatrix{ M^0 \ar@<0.5ex>[r]^{d^0}& M^1 \ar@<0.5ex>[l]^{d^1}  })\longleftrightarrow (\xymatrix{ \varrho(M^1) \ar@<0.5ex>[r]^{\varrho(d^1)}& \varrho(M^0) \ar@<0.5ex>[l]^{\varrho(d^0)}  }).$$
	Let $\cc_{\Z_2}(\ca)^{\varrho^\sharp}$ be the $\varrho^\sharp$-invariant subcategory of $\cc_{\Z_2}(\ca)^{\varrho^\sharp}$. In fact, $\cc_{\Z_2}(\ca)^{\varrho^\sharp}$ is formed by $(\xymatrix{ M \ar@<0.5ex>[r]^{d}& \varrho(M) \ar@<0.5ex>[l]^{\varrho(d)}  })$ for $M\in\ca$.
	
	\begin{lemma}
		\label{lem:invariant}
		For any $(\ca,\varrho)$ with $\varrho$ an involution, we have 
		that $\cc_{\varrho}(\ca)$ is equivalent to $\cc_{\Z_2}(\ca)^{\varrho^\sharp}$.
	\end{lemma}
	
	\begin{proof}
		By definitions, there exists a natural functor
		\begin{align}
			\label{functor:F}
			F:\cc_{\varrho}(\ca)\longrightarrow\cc_{\Z_2}(\ca)
		\end{align} sending $(M,d)\rightarrow (\xymatrix{ M \ar@<0.5ex>[r]^{d}& \varrho(M) \ar@<0.5ex>[l]^{\varrho(d)}  })$. One can show that $F$ induces an equivalence $\cc_{\varrho}(\ca)\simeq \cc_{\Z_2}(\ca)^{\varrho^\sharp}$.
	\end{proof}
	
	The functor $F$ in \eqref{functor:F} admits a left and also right adjoint functor
	\begin{align}
		\label{functor:pi2}
		\pi^2:\cc_{\Z_2}(\ca)\longrightarrow \cc_{\varrho}(\ca)
	\end{align}
	by sending
	$$ \big(\xymatrix{ M^0 \ar@<0.5ex>[r]^{d^0}& M^1 \ar@<0.5ex>[l]^{d^1}  }\big)\mapsto \big(M^0\oplus \varrho M^1, \begin{bmatrix} & \varrho(d^1)\\
		d^0&\end{bmatrix}
	\big).$$	
	There exists an exact functor 
	\begin{align}
		\label{eq:pii}
		\pi^b: \cc^b(\ca)\longrightarrow \cc_{\varrho}(\ca),
	\end{align}
	which sends a complex $\xymatrix{\cdots\ar[r]& M^{i-1} \ar[r]^{d^{i-1}} & M^i\ar[r]^{d^i} &M^{i+1}\ar[r] & \cdots}$ to the $\varrho$-complex
	$$\Big(\bigoplus_{i\in\Z} M^{2i}\oplus \varrho(M^{2i+1}), \begin{bmatrix}&\varrho (d^{2i+1}) \\ d^{2i}  & \end{bmatrix}_{i\in\Z}\Big).$$
	
	
	Note that the functors $F,\pi^2,\pi^b$ preserve acyclic complexes.	
	
	In this paper, we also need the natural covering functor
	$\pi^b_2: \cc^b(\ca)\rightarrow \cc_{\Z_2}(\ca)$; see \cite[\S2.1]{LP16}. Obviously,  $\pi^2\pi^b_2=\pi^b$.

	\subsection{$\imath$Quiver algebras}
	\label{subsec:i-quiver}
	
	Let $Q=(Q_0,Q_1)$ be a quiver (not necessarily acyclic), where $Q_0$ is the set of vertices and $Q_1$ is  the set of arrows. A representation of $Q$ over $\bfk$ is $(M_i,x_\alpha)_{i\in\I,\alpha\in\Omega}$ by assigning a finite-dimensional $\bfk$-linear space $M_i$ to each vertex $i\in Q_0$, and a  $\bfk$-linear maps $x_\alpha: M_i\rightarrow M_j$ to each arrow $\alpha:i\rightarrow j\in Q_1$.
	A representation $(M_i,x_\alpha)$ of $Q$ is called nilpotent if $x_{\alpha_r}\cdots x_{\alpha_{2}}x_{\alpha_1}$ is nilpotent for any cyclic paths $\alpha_r\cdots \alpha_2\alpha_1$.
	Let $\mod_\bfk^{\rm nil}(Q)$ be the category of finite dimensional nilpotent representations of $Q$ over $\bfk$. Then $\rep_\bfk^{\rm nil}(Q)$ is a hereditary abelian category.

	Let $(Q,\varrho)$ be an $\imath$quiver \cite{LW19a}. That is, $\varrho$ is an involution of $Q$. The $\imath$quiver algebra $\Lambda^\imath$ \cite{LW19a} is defined to be $\bfk \ov{Q}/\ov{I}$, where the bound quiver $(\ov{Q},\ov{I})$ is defined as follows:
	\begin{itemize}
		\item[(a)] $\ov{Q}$ is constructed from $Q$ by adding a loop $\varepsilon_i$ at the vertex $i\in Q_0$ if $\varrho( i)=i$, and adding an arrow $\varepsilon_i: i\rightarrow \varrho(i)$ for each $i\in Q_0$ if $\varrho( i)\neq i$;
		\item[(b)] $\ov{I}$ is generated by
		\begin{itemize}
			\item[(1)] (Nilpotent relations) $\varepsilon_{i}\varepsilon_{\varrho(i)}$ for any $i\in\I$;
			\item[(2)] (Commutative relations) $\varepsilon_i\alpha-\varrho(\alpha)\varepsilon_j$ for any arrow $\alpha:j\rightarrow i$ in $Q_1$.
		\end{itemize}
	\end{itemize}
	
	Let $\mod^{\rm nil}(\Lambda^\imath)$ be the category of finite-dimensional nilpotent representations of $(\ov{Q},\ov{I})$. Obviously, 
	$\varrho$ induces an involution of $\rep_\bfk^{\rm nil}(Q)$, which is also denoted by $\varrho: \rep_\bfk^{\rm nil}(Q)\rightarrow \rep_\bfk^{\rm nil}(Q)$. Let $\cc_{\varrho}(\rep_\bfk^{\rm nil}(Q))$ be the $\imath$-category of $(\rep_\bfk^{\rm nil}(Q),\varrho)$. 
	
	\begin{lemma}
		\label{lem:iQ-iCat}
		Retain the notations as above. Then $\mod^{\rm nil}(\Lambda^\imath)\simeq \cc_{\varrho}(\rep_\bfk^{\rm nil}(Q))$.
	\end{lemma}

	\begin{proof}
		Let $R_2:=\bfk(\xymatrix{1 \ar@<0.5ex>[r]^{\varepsilon} & 1' \ar@<0.5ex>[l]^{\varepsilon'}}) /(\varepsilon' \varepsilon,\varepsilon\varepsilon ')$. Let $\Lambda=\bfk Q\otimes_\bfk R_2$. Then we have an equivalence of categories $\mod^{\rm nil}(\Lambda)\simeq \cc_{\Z_2}(\rep_\bfk^{\rm nil}(Q))$. Moreover, $\varrho$ induces an involution $\varrho^\sharp:\mod^{\rm nil}(\Lambda)\rightarrow \mod^{\rm nil}(\Lambda)$; see \cite[Lemma 2.4]{LW19a}. By \cite[Remark 2.11]{LW19a}, we have that $\mod^{\rm nil}(\Lambda^{\imath})$ is equivalent to the $\varrho^\sharp$-invariant subcategory $\mod^{\rm nil}(\Lambda)^{\varrho^\sharp}$ of $\mod^{\rm nil}(\Lambda)$. 
		The equivalence $\mod^{\rm nil}(\Lambda)\simeq \cc_{\Z_2}(\rep_\bfk^{\rm nil}(Q))$ induces that $\mod^{\rm nil}(\Lambda)^{\varrho^\sharp}\simeq \cc_{\Z_2}( \rep_\bfk^{\rm nil}(Q))^{\varrho^\sharp} \simeq\cc_\varrho(\rep_\bfk^{\rm nil}(Q))$ by using Lemma \ref{lem:invariant}.
	\end{proof}

	\subsection{Homological dimensions of acyclic complexes}

	For any $f:M\rightarrow N$, we denote $C_f:=\Big(M\oplus \varrho(N),\begin{pmatrix} 0 &0\\
		f&0 \end{pmatrix}\Big)$. Then there exists an exact sequence in $\cc_\varrho(\ca)$:
	\begin{align}
		0\longrightarrow \varrho(N)\longrightarrow C_f\longrightarrow M\longrightarrow0.
	\end{align}
	Note that $H(C_f)\cong \Ker (f)\oplus \coker (\varrho (f))$. 
	For any $X\in\ca$, define
	\begin{align}
		\label{eq:KX}
		K_X:=(X\oplus \varrho( X), d ), \text{ where }d=\begin{bmatrix} 0 &0 \\ \Id &0 \end{bmatrix}: X\oplus \varrho (X)\rightarrow \varrho (X)\oplus  X.
	\end{align}
	Then $K_X$ is acyclic.
	
	For any $X^\bullet=(X,d)$, define 
	$$\varrho(X^\bullet):=(\varrho(X),-\varrho(d)).$$ 
	Then there exists a short exact sequence
	\begin{align}
		\label{eq:KXses}
		\xymatrix{0\ar[r]& \varrho( X^\bullet) \ar[rr]^{\begin{bmatrix} -\varrho( d)\\ \Id \end{bmatrix}}&& K_X\ar[rr]^{\begin{bmatrix}\Id&\varrho(d))\end{bmatrix}}& &X^\bullet\ar[r]&0}.
	\end{align}

	\begin{lemma}
		\label{lemma extension 2 zero}
		Let $\ca$ be a hereditary abelian category.
		For any $M^\bullet,K_X\in\cc_{\varrho}(\ca)$ and $p\geq2$, we have
		\begin{align*}
			\Ext^{p}_{\cc_{\varrho}(\ca)}(K_X,M^\bullet)=0,\qquad 
			\Ext^{p}_{\cc_{\varrho}(\ca)}(M^\bullet,K_X)=0. 
		\end{align*}
	\end{lemma}
	
	\begin{proof}
		We only prove the first equality since the second one is dual. 
		
		Using Yoneda product, it is enough to prove 
		$\Ext^2_{\cc_{\varrho}(\ca)}(K_X,M^\bullet)=0$.
		Given a short exact sequence 
		$$\xi:\quad 0\longrightarrow M^\bullet\stackrel{f_3}{\longrightarrow} M_2^\bullet\stackrel{f_2}{\longrightarrow} M_1^\bullet\stackrel{f_1}{\longrightarrow} K_X\longrightarrow0,$$
		assume $M^\bullet=(M,d)$, $M_i^\bullet=(M_i,d_i)$, for $i=1,2$. 
		By \eqref{eq:KX}, we have $f_1=\begin{bmatrix} \varrho(g)d_1\\ g \end{bmatrix}$ for some $g:M_1\rightarrow \varrho(X)$. 
		Using pullback, we have the following commutative diagram in $\ca$ with rows short exact:
		\[\xymatrix{ M_2\ar[r]^{g_2} \ar@{=}[d] & W_1\ar[rr]^{g_1} \ar[d]^{h_1} && X\ar[d]^{\tiny\begin{bmatrix}
					1\\0
			\end{bmatrix}} \\
			M_2\ar[r]^{f_2} & M_1\ar[rr]_{\tiny \begin{bmatrix} \varrho(g)d_1\\ g \end{bmatrix} } && X\oplus \varrho(X) ,}\]
		which yields the following commutative diagram with rows exact:
		\[\xymatrix{ 0\ar[r] & M\ar[r]^{f_3} \ar@{=}[d]  & M_2\ar[r]^{g_2} \ar@{=}[d]  & W_1\ar[r]^{g_1} \ar[d]^{h_1} & X\ar[r] \ar[d]^{\tiny\begin{bmatrix}
					1\\0
			\end{bmatrix} }&0
			\\
			0\ar[r]& M\ar[r]^{f_3} & M_2\ar[r]^{f_2} &M_1\ar[r]^{f_1} & X\oplus \varrho(X)\ar[r] &0.}\]
		From the exact sequence of the first row, we obtain an exact sequence in $\cc_\varrho(\ca)$:
		\[\eta: \quad 0\longrightarrow K_M \stackrel{\alpha_3}{\longrightarrow} K_{M_2}\stackrel{\alpha_2}{\longrightarrow} K_{W_1}\stackrel{\alpha_1}{\longrightarrow} K_X\longrightarrow0,\]
		where $\alpha_3=\begin{bmatrix} f_3&0 \\ 0&\varrho(f_3)
		\end{bmatrix}$, $\alpha_2=\begin{bmatrix} g_2&0 \\ 0&\varrho(g_2)
		\end{bmatrix}$, $\alpha_1=\begin{bmatrix} g_1&0 \\ 0&\varrho(g_1)
		\end{bmatrix}$. Note that $[\eta]=0$ in $\Ext^2_{\cc_\varrho(\ca)}(K_X,K_M)$ since $\ca$ is hereditary. 
		In addition we get a morphism from $\eta$ to $\xi$ as follows:
		\[\xymatrix{ 0\ar[r] & K_M\ar[r]^{\alpha_3} \ar[d]^{h_3}  & K_{M_2}\ar[r]^{\alpha_2} \ar[d]^{h_2}  & K_{W_1}\ar[r]^{\alpha_1} \ar[d]^{h} & K_X\ar@{=}[d] \ar[r]&0
			\\
			0\ar[r]& M^\bullet\ar[r]^{f_3} & M_2^\bullet\ar[r]^{f_2} &M_1^\bullet\ar[r]^{f_1} & K_X\ar[r] &0,}\]
		where $h_3=\begin{bmatrix}
			1&\varrho(d)     
		\end{bmatrix}$, $h_2=\begin{bmatrix}
			1&\varrho(d_2)     
		\end{bmatrix}$, $h=\begin{bmatrix}
			h_1&\varrho(d_1)\varrho(h_1)     
		\end{bmatrix}$.
		Therefore, we have $[\xi]=0$ in $\Ext^2_{\cc_\varrho(\ca)}(K_X,M^\bullet)$ by using $[\eta]=0$ and $\Ext^2_{\cc_\varrho(\ca)}(K_X,h_3)([\eta])=[\xi]$.
	\end{proof}

	\begin{proposition}
		\label{prop:Extp=0}
		Let $\ca$ be a hereditary $\bfk$-linear abelian   category. For any $M^\bullet,K^\bullet\in\cc_{\varrho}(\ca)$ with $K^\bullet$ acyclic and $p\geq2$, we have
		\begin{align}
			\label{eq:Extp=0}
			\Ext^{p}_{\cc_{\varrho}(\ca)}(K^\bullet,M^\bullet)=0,\qquad 
			\Ext^{p}_{\cc_{\varrho}(\ca)}(M^\bullet,K^\bullet)=0. 
		\end{align}  
	\end{proposition}

	\begin{proof}
		We denote by $\cc^{can}_{\varrho,ac}(\ca)$ the finite extension closure of the acyclic $\varrho$-complexes $K_X$  for all $X\in\ca$. Then $\cc^{can}_{\varrho,ac}(\ca)$ is closed under taking isomorphisms. By Lemma \ref{lemma extension 2 zero}, we can obtain that \eqref{eq:Extp=0} holds for any $K^\bullet\in \cc^{can}_{\varrho,ac}(\ca)$. 
		
		By induction on the width of bounded complexes, it is routine to prove that $\pi^b(X^\bullet)\in\cc^{can}_{\varrho,ac}(\ca)$ for any acyclic complex $X^\bullet$ in $\cc^b(\ca)$; see \eqref{eq:pii} for the definition of $\pi^b$. So \eqref{eq:Extp=0} holds for any $\pi^b(X^\bullet)$.
		
		For any acyclic complex $Y^\bullet=\xymatrix{(Y^0\ar@<0.5ex>[r]^{d^0} &  Y^1)\ar@<0.5ex>[l]^{d^1}}\in\cc_{\Z_2}(\ca)$, by using \cite[Lemma 2.1]{LP16}, we have the following
		short exact sequences
		\begin{equation*}
			0\longrightarrow \pi^b_2(Z^\bullet)\longrightarrow \pi^b_2(U^\bullet)\longrightarrow Y^\bullet\longrightarrow0 \text{ and }0\longrightarrow Y^\bullet\longrightarrow \pi^b_2(V^\bullet) \longrightarrow\pi^b_2(W^\bullet)\longrightarrow0
		\end{equation*}
		with $U^\bullet,V^\bullet,Z^\bullet,W^\bullet\in\cc^{b}(\ca)$ acyclic. 
		It follows that  \eqref{eq:Extp=0} holds for any $\pi^2(Y^\bullet)$.
		
		By \eqref{functor:pi2}, we have that $\pi^2:\cc_{\Z_2}(\ca)\rightarrow \cc_{\varrho}(\ca)$ is a covering functor with order $2$. And $\pi^2$ admits a left and also right adjoint functor $F: \cc_\varrho(\ca)\rightarrow \cc_{\Z_2}(\ca)$. Then $K^\bullet$ is a direct summand of $\pi^2 F(K^\bullet)$ by the assumption $\Char \bfk\neq 2$. So \eqref{eq:Extp=0} holds for any $K^\bullet$.
	\end{proof}
	
	\begin{remark}
		We expect that Proposition \ref{prop:Extp=0} also holds for  arbitrary field $\bfk$ with $\Char\bfk =2$. In fact, it holds for any $\imath$quiver algebras by \cite{LW20a}, and for any hereditary $\bfk$-linear abelian   category $\ca$ if $\varrho=\Id$ by \cite{LinP}; see also \cite[Lemma 2.2]{LRW20a}. 
	\end{remark}

	\subsection{Triangulated orbit categories}
	
	\begin{lemma}
		\label{lem:dense}
		If $\ca$ is a hereditary abelian category, then  $M^\bullet\cong H(M^\bullet)$ in $\cd_\varrho(\ca)$
		for any $M^\bullet\in\cc_{\varrho}(\ca)$.
	\end{lemma}

	\begin{proof}
		The proof is inspired by \cite[Lemma 5.1]{St17}.
		
		Assume $M^\bullet=(M,d)$. Denote by $p:M\rightarrow \Im (d)$ the epimorphism induced by $d$.
		Since $\ca$ is a hereditary abelian category, the epimorphism $\varrho (p):\varrho(M)\rightarrow \Im (\varrho (d))$ induces an epimorphism
		$$\xymatrix{\Ext^1_\ca\big(H(M^\bullet), \varrho(M)\big)\ar@{->>}[r] &\Ext^1_\ca\big(H(M^\bullet),\Im (\varrho (d))\big).}$$
		In particular, we obtain the commutative diagram with rows short exact:
		\begin{align}
			\label{eq:comdiag1}
			\xymatrix{0\ar[r]&\varrho (M)\ar[r]^f \ar[d]^{\varrho (p)} & E\ar[r]^g\ar[d]^s &H(M^\bullet)\ar@{=}[d] \ar[r]&0\\
				0\ar[r] &	\Im (\varrho (d))\ar[r] & \Ker (d)\ar[r] & H(M^\bullet) \ar[r]&0 }
		\end{align}
		Then there is a short exact sequence in $\cc_{\varrho}(\ca)$:
		\begin{align*}
			\xymatrix{0\ar[r] &K_{M}\ar[rr]^{\begin{bmatrix} \Id&0\\0&f \end{bmatrix}}&& C_{\varrho (f)}\ar[rr]^{\begin{bmatrix}0&g\end{bmatrix}}&& H(M^\bullet)\ar[r]&0.}
		\end{align*}
		So $H(M^\bullet)\cong C_{\varrho (f)}$ in $\cd_\varrho(\ca)$.
		
		On the other hand, consider the morphism $h:=\begin{bmatrix}\Id_M&\tilde{s}\end{bmatrix}: C_{\varrho (f)}\rightarrow M^\bullet$, where $\tilde{s}$ is the morphism induced by $s:E\rightarrow \Ker (d)$.
		The map induced in homology by $h$ therefore appears as the cokernel
		\begin{align*}
			\xymatrix{0\ar[r]&\Im (f)\ar[r] \ar[d]^{s'} & E\ar[r]\ar[d]^s &H(C_{\varrho (f)})\ar[d]^{H(h)}\ar[r]&0 \\
				0\ar[r]&\Im (\varrho (d))\ar[r] & \Ker (d)\ar[r] & H(M^\bullet) \ar[r]&0 }
		\end{align*}
		which is nothing but the diagram in \eqref{eq:comdiag1}, and then $h$ is a quasi-isomorpism.
	\end{proof}

	\begin{lemma}
		\label{lem:Extcrho}
		For any $X,Y\in\ca$, we have 
		$$\Ext^1_{\cc_\varrho(\ca)}(X,Y)\cong \Ext^1_\ca(X,Y)\oplus \Hom_\ca(X,\varrho (Y)).$$
	\end{lemma}
	
	\begin{proof}
		The forgetful functor induces an epimorphism of groups $\Ext^1_{\cc_\varrho(\ca)}(X,Y)\rightarrow \Ext^1_\ca(X,Y)$, one can prove that its kernel is isomorphic to $\Hom_\ca(X,\varrho (Y))$. So the desired formula follows.	
	\end{proof}
	
	Let $M^\bullet=(M,d)$ and $N^\bullet=(N,e)$.
	For any morphism $a:M^\bullet\rightarrow N^\bullet$ in $\cc_\varrho(\ca)$, the cone of $a$ is 
	\begin{align}
		\cone(a):=\big(N\oplus \varrho M, \begin{bmatrix} e &\varrho a \\ &-\varrho (d) \end{bmatrix}\big)\in\cc_\varrho(\ca).
	\end{align}
	Then 
	\begin{align}
		\xymatrix{M^\bullet\ar[r]^a&N^\bullet\ar[rr]^{\begin{bmatrix} 1\\0 \end{bmatrix}}&& \cone(a)\ar[rr]^{\begin{bmatrix}0&1\end{bmatrix}}&&\Sigma M^\bullet }
	\end{align}
	is a distinguished triangle in $\cd_\varrho(\ca)$.
	
	\begin{proposition}
		\label{prop:Extconincide}
		Let $(\ca,\varrho)$ be a category with involution, where $\ca$ is a hereditary $\bfk$-linear abelian  category.
		For any $X,Y\in\ca$, we have
		$\Ext^1_{\cc_\varrho(\ca)}(X,Y)\cong \Hom_{\cd_\varrho(\ca)}(X,\varrho (Y))
		$. 
	\end{proposition}
	
	\begin{proof}
		
		Let $A^\bullet$, $B^\bullet$ be any two complexes. 
		For any short exact sequence $\xi: 0\rightarrow B^\bullet\xrightarrow{a} C^\bullet \xrightarrow{b} A^\bullet\rightarrow0$,  we have a morphism $w=\begin{bmatrix}0&1\end{bmatrix}:\cone(a)\rightarrow \Sigma B^\bullet$, and a quasi-isomorphism $s=\begin{bmatrix}b&0\end{bmatrix}:\cone(a)\rightarrow A^\bullet$, which yields a homomorphism 
		$\delta_\xi= ws^{-1}\in\Hom_{\cd_\varrho}(\ca)(A^\bullet,\Sigma B^\bullet)$.
		This induces a morphism of abelian groups
		\begin{align*}
			\phi: &\Ext^1_{\cc_\varrho(\ca)}(A^\bullet,B^\bullet)\rightarrow\Hom_{\cd_\varrho(\ca)}(A^\bullet,\Sigma B^\bullet)
			,\qquad [\xi]\mapsto \delta_\xi.
		\end{align*}
		In particular, we have a triangle $B^\bullet\xrightarrow{a} C^\bullet \xrightarrow{b} A^\bullet\stackrel{\delta_\xi}{\rightarrow} \Sigma B^\bullet$ in $\cd_\varrho(\ca)$. 
		
		The morphism $\delta_\xi$ is functorial on $\xi$. Explicitly, given another short exact sequence $\eta: 0\rightarrow M^\bullet\rightarrow L^\bullet\rightarrow N^\bullet\rightarrow0$, and a commutative diagram:
		\[\xymatrix{ 0\ar[r] & B^\bullet \ar[r] \ar[d]^f&C^\bullet \ar[r]\ar[d]^g & A^\bullet\ar[d]^h\ar[r]& 0
			\\
			0\ar[r]& M^\bullet\ar[r] &L^\bullet \ar[r]& N^\bullet \ar[r]&0}\]
		we have a commutative diagram in $\cd_\varrho(\ca)$:
		\[\xymatrix{ B^\bullet \ar[r] \ar[d]^f&C^\bullet \ar[r]\ar[d]^g & A^\bullet\ar[d]^h\ar[r]^{\delta_\xi}& \Sigma B^\bullet\ar[d]^{\Sigma f}
			\\
			M^\bullet\ar[r] &L^\bullet \ar[r]& N^\bullet \ar[r]^{\delta_\eta}&\Sigma M^\bullet}\]

		Given a short exact sequence $\xi:0\rightarrow Y\xrightarrow{a} Z^\bullet \xrightarrow{b} X\rightarrow0$, denote $Z^\bullet=(Z,d)$. 
		The differentials of $X,Y$ are zero, so there exists a unique morphism $c: X\rightarrow \varrho (Y)$ such that $d=\varrho(a) c b$. Note that $c$ is the connecting morphism of $\xi$.
		
		If $\de_\xi=0$, then the triangle $Y\xrightarrow{a} Z^\bullet \xrightarrow{b} X\rightarrow \varrho (Y)$ splits. It yields that $Z^\bullet \cong X\oplus Y$ in $\cd_\varrho(\ca)$, and then $H(Z^\bullet)\cong X\oplus Y$. Moreover, $H(a):Y\rightarrow H(Z^\bullet)$ is a section, and $H(b):H(Z^\bullet)\rightarrow X$ is a retraction. Together with the connecting morphism of $\xi$, we have the differential $d=0$ and then $Z^\bullet\cong(X\oplus Y,0)$ in $\cc_\varrho(\ca)$. So $[\xi]=0$ in $\Ext^1_{\cc_\varrho(\ca)}(X,Y)$, and then $\phi$ is injective. 
		
		It remains to prove that $\phi$ is surjective. For any morphism in $\Hom_{\cd_\varrho(\ca)}(X,\varrho (Y))$, we assume that it is of the form 
		$ws^{-1}$, where $s: U^\bullet\rightarrow X$ is a quasi-isomorphism, and $w: U^\bullet\rightarrow \varrho (Y)$ is a morhism of complexes. Let $C^\bullet$ be the cone of the identity morphism of $U^\bullet$. Then $C^\bullet$ is acyclic. We have a short exact sequence $0\rightarrow Y\rightarrow \Sigma \cone(w)\rightarrow U^\bullet \rightarrow 0$. 
		Using the natural embedding $\iota: U^\bullet\rightarrow C^\bullet$ and $s:U^\bullet\rightarrow X$, we obtain an injective homomorphism $\begin{bmatrix} s\\ \iota \end{bmatrix}: U^\bullet\rightarrow X\oplus C^\bullet$ with its cokernel acyclic.
		Applying Proposition  \ref{prop:Extp=0}, there exists a commutative diagram in $\cc_\varrho(\ca)$ with all rows short exact:
		\[\xymatrix{ Y\ar[r]\ar@{=}[d]&\Sigma\cone( w) \ar[r] \ar[d]& U^\bullet \ar[d]^{\tiny\begin{bmatrix} s\\ \iota \end{bmatrix}} 
			\\
			Y\ar[r] \ar@{=}[d]  & W^\bullet \ar[r] & X\oplus C^\bullet 
			\\
			Y\ar[r] & L^\bullet \ar[r] \ar[u] &X \ar[u]_{\tiny\begin{bmatrix} 1\\0\end{bmatrix}} }\]
		Denote by $\xi_i$ the short exact sequence of the $i$-th row in the above diagram for $i=1,2,3$. 
		Then we have a commutative diagram in $\cd_\varrho(\ca)$ with all rows triangles:
		\[\xymatrix{ Y\ar[r]\ar@{=}[d]&\Sigma\cone( w) \ar[r] \ar[d]& U^\bullet \ar[d]^{\tiny\begin{bmatrix} s\\ \iota \end{bmatrix}} \ar[r]^{\delta_{\xi_1}} &\varrho (Y)\ar@{=}[d]
			\\
			Y\ar[r] \ar@{=}[d]  & W^\bullet \ar[r] & X\oplus C^\bullet \ar[r]^{\delta_{\xi_2}} &\varrho (Y)\ar@{=}[d]
			\\
			Y\ar[r] & L^\bullet \ar[r] \ar[u] &X \ar[u]_{\tiny\begin{bmatrix} 1\\0\end{bmatrix}} \ar[r]^{\delta_{\xi_3}} &\varrho (Y) }\]
		One can check that $\delta_{\xi_1}=w$ by using the definition.
		Then $\phi([\xi_3])=ws^{-1}$. 
		The proof is completed.	
	\end{proof}

	The involution $\varrho$ induces an auto-equivalence $\widehat{\varrho}: \cd^b(\ca)\rightarrow \cd^b(\ca)$.

	\begin{lemma}
		Let $(\ca,\varrho)$ be a category with involution, where $\ca$ is a hereditary abelian category. Then $\cd^b(\ca)/\Sigma \circ \widehat{\varrho}$ is a triangulated orbit category \`a la Keller \cite{Ke05}.
	\end{lemma}
	
	\begin{proof}
		The proof is the same as for \cite[Lemma 3.17]{LW19a}, and hence omitted here.
	\end{proof}

	\begin{theorem}
		\label{prop: orbit}
		Let $\ca$ be a hereditary abelian category with an involution $\varrho$. Then $\cd_\varrho(\ca)\simeq \cd^b(\ca)/\Sigma \circ \widehat{\varrho}$. 
	\end{theorem}
	
	\begin{proof}
		Consider $\cc_{\varrho}(\ca)$ and $\cc^b(\ca)$ as DG categories. Then the exact functor $\pi^b$ is a DG functor; see \eqref{eq:pii}.  $\pi^b$ induces a triangulated functor $\widetilde{\pi}^b: \cd^b(\ca)\rightarrow \cd_\varrho(\ca)$. By definition, $\widetilde{\pi}^b\circ(\Sigma \circ\widehat{\varrho})\cong \widetilde{\pi}^b$ in $\cd_\varrho(\ca)$. So \cite[\S9.4]{Ke05} shows that there exists a triangulated functor
		$\bar{\pi}^b: \cd^b(\ca)/\Sigma\circ \widehat{\varrho}\rightarrow \cd_\varrho(\ca)$ such that the following diagram commutes:
		\[
		\xymatrix{ \cd^b(\ca) \ar[r]^{\tilde{\pi}^b} \ar[d]_{\text{nat.}}& \cd_\varrho(\ca)
			\\
			\cd^b(\ca)/\Sigma\circ \widehat{\varrho} \ar[ur]_{\bar{\pi}^b} & }
		\]
		Lemma \ref{lem:dense} shows that $\widetilde{\pi}^b$ is dense and then so is $\bar{\pi}^b$.
		
		It remains to prove $\bar{\pi}^b$ to be fully faithful. Since $\ca$ is a hereditary abelian category, we only need to consider stalk complexes. For any $X,Y\in\ca$, we have 
		$\Hom_{\cd^b(\ca)/\Sigma\circ \widehat{\varrho}}(X,Y)=\Hom_\ca(X,Y)\oplus \Ext^1_\ca(X,\varrho (Y))\cong \Hom_{\cd_\varrho(\ca)}(X,Y)$, where the last identity follows from Proposition \ref{prop:Extconincide} and 	Lemma \ref{lem:Extcrho}. The proof is completed.
	\end{proof}

	\subsection{Euler forms}

	
	For any $M^\bullet$, its projective dimension $\pd_{\cc_{\varrho}(\ca)} M^\bullet$ is defined  to be the smallest integer $p \geq 0$ with the
	property that $\Ext^{p+1}_{\cc_\varrho(\ca)}(M^\bullet,-)=0$. Its injective dimension 	$\ind_{\cc_{\varrho}(\ca)} M^\bullet$ is defined dually.

	\begin{proposition}
		\label{cor:projfinite}
		For any $M^\bullet\in\cc_\varrho(\ca)$, the following statements are equivalent.
		\begin{enumerate}
			\item $\pd_{\cc_{\varrho}(\ca)} M^\bullet<\infty$,
			\item $\pd_{\cc_{\varrho}(\ca)} M^\bullet\leq1$,
			\item $\ind_{\cc_{\varrho}(\ca)} M^\bullet<\infty$,
			\item $\ind_{\cc_{\varrho}(\ca)} M^\bullet\leq1$,
			\item $M^\bullet$ is acyclic.
		\end{enumerate}
	\end{proposition}
	
	\begin{proof}
		It is enough to prove ``(1)$\Rightarrow$(5)'' and ``(3)$\Rightarrow$(5)''  by Proposition \ref{prop:Extp=0}.
		
		We only prove  ``(1)$\Rightarrow$(5)'' since the other one is similar. Assume $M^\bullet=(M,d)$. 
		By \eqref{eq:KXses}, we have short exact sequences
		\begin{align}
			\label{eq:seseq1}
			0\longrightarrow M^\bullet\longrightarrow K_{\varrho(M)}\longrightarrow \varrho( M^\bullet)\longrightarrow0,
			\\
			\label{eq:seseq2}
			0\longrightarrow\varrho( M^\bullet)\longrightarrow K_{M}\longrightarrow M^\bullet\longrightarrow0,
		\end{align}
		By Lemma \ref{lemma extension 2 zero}, we have $\Ext^p_{\cc_{\varrho}(\ca)}(M^\bullet,-)\cong \Ext^{p+1}_{\cc_{\varrho}(\ca)}( \varrho( M^\bullet),-)\cong \Ext^{p+2}_{\cc_{\varrho}(\ca)}( M^\bullet,-)$ for any $p\geq2$. Since $\pd_{\cc_{\varrho}(\ca)} M^\bullet<\infty$, we have $\Ext^p_{\cc_{\varrho}(\ca)}(M^\bullet,-)=0$ for any $p\geq2$. So $\pd_{\cc_{\varrho}(\ca)} M^\bullet\leq1$. It follows from \eqref{eq:seseq2} that  $\pd_{\cc_{\varrho}(\ca)} \varrho( M^\bullet)\leq1$ by noting that $\pd_{\cc_{\varrho}(\ca)} K_M\leq1$. 
		
		By applying $\Hom_{\cc_{\varrho}(\ca)}(-,\varrho(M^\bullet))$ to \eqref{eq:seseq1}, we have an epimorphism $$\Ext^1_{\cc_{\varrho}(\ca)}(K_{\varrho(M)}, \varrho(M^\bullet))\longrightarrow \Ext^1_{\cc_{\varrho}(\ca)}(M^\bullet,\varrho(M^\bullet)).$$
		considering \eqref{eq:seseq2}, we have the following commutative diagram:
		\[\xymatrix{\varrho(M^\bullet) \ar@{=}[r] \ar[d] & \varrho(M^\bullet)\ar[d] \\
			K_M\ar[r]\ar[d] & U^\bullet \ar[r] \ar[d] & \varrho(M^\bullet) \ar@{=}[d]\\
			M^\bullet\ar[r]& K_{\varrho(M)}\ar[r] & \varrho(M^\bullet)}\]
		which is  a pullback and pushout diagram. Therefore, we have a short exact sequence
		$$0\longrightarrow K_M\longrightarrow U^\bullet\oplus M^\bullet\longrightarrow K_{\varrho(M)}\longrightarrow0.$$
		It follows that $ U^\bullet\oplus M^\bullet$ is acyclic since $K_M,K_{\varrho(M)}$ are. In particular $M^\bullet$ is acyclic.
	\end{proof}

	For any $K^\bullet$, $M^\bullet$ with $K^\bullet$ acyclic, by Proposition \ref{cor:projfinite}, we define
	the Euler forms
	\begin{align}\label{left Euler form}
		\langle K^\bullet,M^\bullet\rangle=&\sum_{i=0}^{+\infty}(-1)^i \dim_\bfk\Ext_{\cc_{\varrho}(\ca)}^i(K^\bullet,M^\bullet)
		\\\notag
		=&\dim_\bfk\Hom_{\cc_\varrho(\ca)}(K^\bullet,M^\bullet)-\dim_\bfk\Ext_{\cc_{\varrho}(\ca)}^1(K^\bullet,M^\bullet),
		\\\label{right Euler form}
		\langle M^\bullet,K^\bullet\rangle=&\sum_{i=0}^{+\infty}(-1)^i \dim_\bfk\Ext^i_{\cc_{\varrho}(\ca)}(M^\bullet,K^\bullet)
		\\\notag
		=&\dim_\bfk\Hom_{\cc_{\varrho}(\ca)}(M^\bullet,K^\bullet)-\dim_\bfk\Ext^1_{\cc_{\varrho}(\ca)}(M^\bullet,K^\bullet).
	\end{align}
	These forms descend to bilinear Euler forms on the Grothendieck groups $K_0(\cc_{\varrho,ac}(\ca))$ and $K_0(\cc_{\varrho}(\ca))$, denoted by the same symbol:
	\begin{align}
		\langle\cdot,\cdot\rangle: K_0(\cc_{\varrho,ac}(\ca))\times K_0(\cc_{\varrho}(\ca))\longrightarrow \Z,
		\\
		\langle\cdot,\cdot\rangle: K_0(\cc_{\varrho}(\ca))\times K_0(\cc_{\varrho,ac}(\ca))\longrightarrow \Z.
	\end{align}
	We have used the same symbol by noting that these two forms coincide when restricting to $K_0(\cc_{\varrho,ac}(\ca))\times K_0(\cc_{\varrho,ac}(\ca))$.
	
	By abusing notations, we also denote by 
	\begin{align}
		\langle\cdot,\cdot\rangle: K_0(\ca)\times K_0(\ca)\longrightarrow \Z
	\end{align}
	the Euler form of $\ca$, that is,
	\begin{align}
		\langle A,B\rangle=\dim_\bfk\Hom_\ca(A,B)-\dim_\bfk \Ext^1_\ca(A,B).
	\end{align}
	Let $(\cdot,\cdot)$ be the symmetric Euler form, that is,
	\begin{align}
		( A,B)=\langle A,B\rangle+\langle B,A\rangle.
	\end{align}

	\begin{lemma}
		\label{lem:Euler form}
		For any $X,Y\in \ca$, we have
		\begin{align}
			\label{eq:KXYEuler}
			\langle K_X, Y \rangle&=\langle X,Y\rangle,\qquad \langle Y, K_X\rangle=\langle Y,\varrho (X)\rangle,
			\\
			\label{eq:KXKYEuler}
			&\langle K_X,K_Y\rangle=\langle X,Y\rangle +\langle X,\varrho(Y)\rangle.
		\end{align}
	\end{lemma}
	
	\begin{proof}
		
		For \eqref{eq:KXYEuler}, We only prove the first one and omit the similar proof of the other one. Obviously, $\Hom_{\cc_\varrho(\ca)}(K_X,Y)\cong \Hom_\ca(X,Y)\cong\Hom_{\cc_{\varrho}(\ca)}(X,Y)$.
		It is enough to prove 
		$$\dim_\bfk\Ext^1_{\cc_\varrho(\ca)}(K_X,Y)=\dim_\bfk\Ext^1_\ca(X,Y).$$
		
		We have a short exact sequence in $\cc_\varrho(\ca)$: 
		\begin{align}
			\label{eq:sses}
			\xymatrix{ 0\ar[r] & \varrho(X)\ar[rr]^{\begin{bmatrix} 0\\1\end{bmatrix}} &&K_X\ar[rr]^{ \begin{bmatrix} 1&0\end{bmatrix}} & &X\ar[r] &0. }
		\end{align}
		It induces a long exact sequence
		\begin{align*}
			0\longrightarrow\Hom_{\cc_\varrho(\ca)}(\varrho(X),Y)\longrightarrow\Ext^1_{\cc_\varrho(\ca)}(X,Y)\longrightarrow \Ext^1_{\cc_\varrho(\ca)}(K_X,Y)\stackrel{\psi}{\longrightarrow} \Ext^1_{\cc_\varrho(\ca)}(\varrho(X),Y)
		\end{align*}
		by noting that the map $\Hom_{\cc_\varrho(\ca)}(K_X,Y)\rightarrow \Hom_{\cc_\varrho(\ca)}(\varrho(X),Y)$ is zero. 
		
		We claim $\psi=0$. Then 
		\begin{align*}
			\dim_\bfk\Ext^1_{\cc_\varrho(\ca)}(K_X,Y)&=\dim_\bfk\Ext^1_{\cc_\varrho(\ca)}(X,Y)-\dim_\bfk \Hom_{\cc_\varrho(\ca)}(\varrho(X),Y)
			\\
			&=\dim_\bfk\Ext^1_{\ca}(X,Y)+\dim_\bfk \Hom_{\ca}(X,\varrho(Y))-\dim_\bfk \Hom_{\ca}(\varrho(X),Y)
			\\
			&=\dim_\bfk\Ext^1_{\ca}(X,Y),
		\end{align*}
		where the second equality follows from 
		Lemma \ref{lem:Extcrho}.

		Let us prove the claim. For any short exact sequence 
		\[\xymatrix{\xi:& 0\ar[r] &Y\ar[r]^g &Z^\bullet \ar[rr]^{\begin{bmatrix}f_1 \\ f_2 \end{bmatrix}} && K_X\ar[r] &0, }\]
		denote by $Z^\bullet=(Z,d)$. We obtain $dg=0$, and then there exists $\begin{bmatrix} h_1&h_2 \end{bmatrix}:X\oplus\varrho(X)\rightarrow \varrho(Z)$ such that
		the following  diagram in $\ca$ commutes:
		\[\xymatrix{ Y\ar[rr]^0\ar[d]^g && \varrho (Y)\ar[d]^{\varrho(g)}
			\\
			Z\ar[rr]^d \ar[d]^{\tiny \begin{bmatrix} f_1\\ f_2\end{bmatrix}}   && \varrho(Z) \ar[d]^{\tiny\begin{bmatrix} \varrho(f_1)\\\varrho(f_2)\end{bmatrix}} 
			\\
			X\oplus \varrho(X) \ar[rr]^{\tiny \begin{bmatrix} 0&0\\\Id&0\end{bmatrix}}&&\varrho(X)\oplus X}\]
		In particular, 
		\begin{align*}
			\begin{bmatrix} \varrho(f_1) h_1& \varrho(f_1)h_2\\\varrho(f_2)h_1& \varrho(f_2)h_2\end{bmatrix}=\begin{bmatrix} 0&0\\\Id&0\end{bmatrix},
		\end{align*}
		and then $\varrho(f_2)h_1=\Id$, $\varrho(f_1)h_1=0$. It follows that
		\begin{align}
			\label{eq:pb}
			\begin{bmatrix}f_1\\ f_2\end{bmatrix}\varrho(h_1)=\begin{bmatrix}0\\ 1\end{bmatrix}
		\end{align}

		Using pullback, we have the following  commutative diagram  in $\cc_\varrho(\ca)$:
		\[\xymatrix{ Y \ar[r] \ar@{=}[d] & L^\bullet \ar[r]^{t_1} \ar[d]_{t_2}  & \varrho(X) \ar[d]^{\tiny\begin{bmatrix} 0\\1\end{bmatrix}} 
			\\
			Y\ar[r]^g & Z^\bullet \ar[r]^{\tiny\begin{bmatrix} f_1\\f_2\end{bmatrix} }\ar[d]^{f_1} &K_X\ar[d]^{\tiny\begin{bmatrix} 1&0\end{bmatrix}}
			\\
			&X\ar@{=}[r]&X
		}\]
		Denote by $\eta$ the short exact sequence at the first row. Then $\psi([\xi])=[\eta]$. 
		Using  \eqref{eq:pb} and the universal property of pullback, there exists a unique morphism $s:\varrho(X)\rightarrow L^\bullet$ such that the following diagram commutes:
		\[\xymatrix{
			&\varrho(X)\ar@/_/[ddr]_{\varrho(h_1)}\ar[dr]^{s}\ar@/^/[drr]^{\Id}&&
			\\
			& & L^\bullet \ar[r]^{t_1} \ar[d]_{t_2}  & \varrho(X) \ar[d]^{\tiny\begin{bmatrix} 0\\1\end{bmatrix}} 
			\\
			& & Z^\bullet \ar[r]^{\tiny\begin{bmatrix} f_1\\f_2\end{bmatrix} }\ar[d]^{f_1} &K_X\ar[d]^{\tiny\begin{bmatrix} 1&0\end{bmatrix}}
			\\
			&&X\ar@{=}[r]&X
		}\]
		In particular, $t_1s=\Id$, and then $\eta$ splits, which gives $\psi([\xi])=0$. So $\psi=0$. The claim is proved. 
		
		For \eqref{eq:KXKYEuler}, using \eqref{eq:sses}, we have
		\begin{align*}
			\langle K_X,K_Y\rangle =&\langle X,K_Y\rangle +\langle\varrho(X), K_Y\rangle
			\\
			=& \langle X,\varrho(Y)\rangle +\langle X, Y\rangle.
		\end{align*}
	\end{proof}
	
	For $X\in\ca$, $Y^\bullet=(Y,e)\in\cc_\varrho(\ca)$, we have a short exact sequence in $\cc_\varrho(\ca)$:
	\begin{align*}
		0\longrightarrow \Ker(e)\longrightarrow Y^\bullet  \longrightarrow \Im(e)\longrightarrow0.
	\end{align*}
	By Lemma \ref{eq:KXYEuler},  we have
	\begin{align}
		\label{eq:KXYcompEuler}
		\langle K_X, Y^\bullet\rangle =&\langle K_X, \Ker(e)\rangle+\langle K_X,\Im(e)\rangle
		\\\notag
		=&\langle X, \Ker(e)\rangle+\langle X,\Im(e)\rangle
		\\\notag
		=&\langle X,Y\rangle=\langle X,\res(Y^\bullet)\rangle.
	\end{align}
	Similarly, we have
	\begin{align}
		\langle Y^\bullet,K_X\rangle=\langle \res(Y^\bullet), \varrho(X)\rangle.
	\end{align}
	
	We denote by $\widehat{A}$ the class of $A\in\ca$ in the Grothendieck group $K_0(\ca)$. 
	Similar to \cite[Lemma 3.13]{LP16}, we have the following lemma.
	\begin{lemma}
		\label{lem:equiSES}
		For objects $U^\bullet=(U,d),V^\bullet=(V,e),W^\bullet=(W,f)\in\cc_{\varrho}(\ca)$, 
		if there is a short exact sequence $0\rightarrow U^\bullet\xrightarrow{h_1} V^\bullet\xrightarrow{h_2} W^\bullet\rightarrow0$, then the following statements are equivalent.
		\begin{itemize}
			\item[(i)] $0\rightarrow \Im (d)\rightarrow \Im(e) \rightarrow \Im (f)\rightarrow 0$ is  exact;
			
			\item[(ii)] $0\rightarrow \Ker(d)\rightarrow\Ker(e)\rightarrow \Ker(f)\rightarrow 0$ is  exact;
			
			\item[(iii)] $0\rightarrow H(U^\bullet)\rightarrow H(V^\bullet)\rightarrow H(W^\bullet)\rightarrow 0$ is  exact.
		\end{itemize}
		Here the morphisms are induced by $h_1$ and $h_2$.
		In particular, if $U^\bullet$ or $W^\bullet$ is acyclic,
		then
		$\widehat{\Im (e)}=\widehat{\Im (d)}+ \widehat{\Im (f)}$ in $K_0(\ca)$.
	\end{lemma}
	
	\begin{proof}
		The proof is the same as \cite[Lemma 3.13]{LP16}, hence omitted here.
	\end{proof}
	
	Now we can prove the following proposition which is crucial to construct the $\imath$Hall algebra in the next section.
	
	\begin{theorem}
		\label{prop:Euler}
		For $K^\bullet=(K,d)$ and $M^\bullet \in\cc_{\varrho}(\ca)$ with $K^\bullet$ acyclic, we have 
		\begin{align}
			\label{eq:KMeuler1}
			\langle K^\bullet,M^\bullet\rangle=\langle K_{\Im(d)},M^\bullet\rangle=\langle \Im(d),\res(M^\bullet)\rangle ,
			\\
			\label{eq:MKeuler2}
			\langle M^\bullet,K^\bullet\rangle=\langle M^\bullet,K_{\Im(d)}\rangle =\langle \res(M^\bullet),\varrho(\Im(d))\rangle.
		\end{align}
		In particular, if $M^\bullet$ is also acyclic, then 
		\begin{align}
			\label{eq:MKeuler3}
			\langle K^\bullet,M^\bullet\rangle=\frac{1}{2}\langle \res (K^\bullet),\res (M^\bullet)\rangle.
		\end{align}
	\end{theorem}
	
	\begin{proof}
		Denote by $p:K\rightarrow \Im(d)$ and $l:\Im(d)\rightarrow \varrho(K)$ the canonical morphisms induced by $d:K\rightarrow\varrho(K)$. First, we prove that $\langle K^\bullet, A\rangle =\langle \Im(d),A\rangle$ for $A\in\ca\subseteq\cc_{\varrho}(\ca)$.

		For any $f:K^\bullet\rightarrow A$, we have $f\circ \varrho (d)=0$ by noting that the differential of $A$ is $0$. So there exists a unique $g:\Im(d)\rightarrow A$ such that $gp=f$. Conversely, 
		for any $g:\Im(d)\rightarrow A$, set $f=gp$. Obviously, $f:K^\bullet\rightarrow A$ is a morphism 
		of $\varrho$-complexes.
		So $\Hom_{\cc_\varrho(\ca)}(K^\bullet,A)\cong \Hom_\ca(\Im(d),A)$.

		Let $0\rightarrow A\xrightarrow{\alpha'}C\xrightarrow{\beta'}\Im(d)\rightarrow0$ be a short exact sequence in $\ca$.
		Then we have the commutative diagram by doing pullback:
		\[\xymatrix{ & \varrho(\Im(d))\ar@{=}[r] \ar[d]_{\varrho(\delta)} &\varrho(\Im(d)) \ar[d]^{\varrho (l)} \\
			A\ar[r]^\alpha \ar@{=}[d] & B\ar[r]^\beta \ar[d]^\mu &K \ar[d]^p \\
			A\ar[r]^{\alpha'} &C\ar[r]^{\beta'} &\Im(d)
		}\]
		Let $b:=\delta\beta'\mu:B\rightarrow \varrho(B)$. Then $\varrho(b)b=0$, and so $B^\bullet=(B,b)$ is a $\varrho$-complex. Furthermore, $b\alpha=\delta\beta'\mu\alpha=\delta\beta'\alpha'=0$, and $\varrho(\beta)b=\varrho(\beta)\delta\beta'\mu=l\beta'\mu=lp\beta=d\beta$. So $\alpha: A\rightarrow B^\bullet$, $\beta:B^\bullet\rightarrow K^\bullet$ are morphisms of $\varrho$-complexes. Then we have a short exact sequence $0\rightarrow A\xrightarrow{\alpha} B^\bullet\xrightarrow{\beta} K^\bullet\rightarrow0$. This induces a morphism $\varphi:\Ext^1_\ca(\Im(d),A)\rightarrow \Ext^1_{\cc_\varrho(\ca)}(K^\bullet,A)$.
		
		Let us prove that $\varphi$ is surjective. In fact, for any short exact sequence 
		$$\xi:\quad 0\longrightarrow A\stackrel{\alpha}{\longrightarrow} (B,b)\stackrel{\beta}{\longrightarrow}K^\bullet\longrightarrow0,$$ we have $\beta$ induces an isomorphism $\widetilde{\beta}:\Im(b)\rightarrow \Im(d)$ by Lemma \ref{lem:equiSES} since $K^\bullet$ is acyclic. Denote by $\delta':\Im(b)\rightarrow \varrho(B)$ the canonical embedding. Then $l\widetilde{\beta}=\varrho(\beta)\delta'$. It follows that we have the following commutative diagram:
		\begin{align}
			\label{eq:pullback}
			\xymatrix{ & \varrho(\Im(b))\ar[r]^{\varrho(\widetilde{\beta})} \ar[d]_{\varrho(\delta')} &\varrho(\Im(d)) \ar[d]^{\varrho (l)} \\
				A\ar[r]^\alpha \ar@{=}[d] & B\ar[r]^\beta \ar[d]^\mu &K \ar[d]^p \\
				A\ar[r]^{\alpha'} &C\ar[r]^{\beta'} &\Im(d)
			}
		\end{align}
		The diagram
		\eqref{eq:pullback} is certainly a pullback, and $b=\delta' \widetilde{\beta} ^{-1}p\beta$.
		Let 
		$$\xi':\quad 0\longrightarrow A\stackrel{\alpha'}{\longrightarrow}C\stackrel{\beta}{\longrightarrow}\Im(d)\longrightarrow0.$$ 
		Then we have $\varphi([\xi'])=[\xi]$ by definition. Here, we denote $[\xi]$ its class in $\Ext_{\cc_\varrho(\ca)}(K^\bullet,A)$. Therefore, $\varphi$ is surjective.
		
		In conclusion, we have $\langle K^\bullet,A\rangle \geq\langle \Im(d),A\rangle$.

		Using \eqref{eq:KXses}, we have $\langle K^\bullet,A\rangle +\langle \varrho(K^\bullet),A\rangle =\langle K_K,A\rangle =\langle K,A\rangle$, where the last equality follows from Lemma 
		\ref{lem:Euler form}. Note that $\langle K,A\rangle=\langle \Im(d),A\rangle +\langle \varrho(\Im(d)),A\rangle$. Then $\langle K^\bullet,A\rangle +\langle \varrho(K^\bullet),A\rangle =\langle \Im(d),A\rangle +\langle \varrho(\Im(d)),A\rangle$. 
		
		Using \eqref{eq:KXses}, we have $$\langle K^\bullet,A\rangle +\langle \varrho(K^\bullet),A\rangle =\langle K_K,A\rangle =\langle K,A\rangle =\langle \Im(d),A\rangle +\langle \varrho(\Im(d)),A\rangle,$$ where the second equality follows from Lemma 
		\ref{lem:Euler form}.
		Together with 
		$\langle K^\bullet,A\rangle \geq\langle \Im(d),A\rangle$, and $\langle \varrho(K^\bullet),A\rangle \geq\langle \varrho(\Im(d)),A\rangle$, we have 
		$\langle K^\bullet,A\rangle =\langle \Im(d),A\rangle$. 
		
		Using Lemma \ref{lem:Euler form}, we also have $\langle K_{\Im(d)},A\rangle=\langle \Im(d),A\rangle$.
		
		For general $M^\bullet$, assume $M^\bullet=(M,m)$. Then we have a short exact sequence of $\varrho$-complexes: $0\rightarrow \Ker (m)\rightarrow M^\bullet\rightarrow \Im(m)\rightarrow0$.
		So $\langle K^\bullet, M^\bullet\rangle=\langle K^\bullet, \Ker(m)\rangle+\langle K^\bullet,\Im(m)\rangle=\langle \Im(d),\Ker(m)\rangle+\langle \Im(d),\Im(m)\rangle=\langle \Im(d),M\rangle$.
		This proves \eqref{eq:KMeuler1}.
		
		It is similar to prove \eqref{eq:MKeuler2}.

		For \eqref{eq:MKeuler3}, it is enough to consider $K^\bullet=K_X$ and $M^\bullet=K_Y$. By \eqref{eq:KXKYEuler}, we have 
		\begin{align*}
			\langle K_X,K_Y\rangle=\langle X,Y\oplus \varrho(Y)\rangle=\frac{1}{2}\langle X\oplus \varrho(X),Y\oplus \varrho(Y) \rangle=\langle \res(K_X),\res(K_Y)\rangle.
		\end{align*}
	\end{proof}
	
	
	

	\section{$\imath$Hall algebras}
	\label{sec:Semi-derived}
	
	\subsection{Hall algebras}
	
	Let $\ce$ be an essentially small exact category, linear over the finite field $\bfk=\F_q$.
	Assume that $\ce$ has finite morphism and extension spaces:
	$$|\Hom_\ce(A,B)|<\infty,\quad |\Ext^1_\ce(A,B)|<\infty,\,\,\forall A,B\in\ce.$$

	Given objects $A,B,C\in\ce$, define $\Ext^1_\ce(A,C)_B\subseteq \Ext^1_\ce(A,C)$ to be the subset parameterising extensions with the middle term  isomorphic to $B$. We define the Hall algebra (also called Ringel-Hall algebra) $\ch(\ce)$ to be the $\Q$-vector space whose basis is formed by the isomorphism classes $[A]$ of objects $A$ of $\ce$, with the multiplication
	defined by
	\begin{align}
		\label{eq:mult}
		[A]\diamond [C]=\sum_{[B]\in \Iso(\ce)}\frac{|\Ext_\ce^1(A,C)_B|}{|\Hom_\ce(A,C)|}[B].
	\end{align}
	It is well known that
	the algebra $\ch(\ce)$ is associative and unital. The unit is given by $[0]$, where $0$ is the zero object of $\ce$; see \cite{Rin90,Br}. 

	For any three objects $A,B,C$, let
	\begin{align}
		\label{eq:Fxyz}
		F_{AC}^B:= \big |\{L\subseteq B \mid L \cong C,  B/L\cong A\} \big |.
	\end{align}
	The Riedtmann-Peng formula states that
	\[
	F_{AC}^B= \frac{|\Ext^1(A,C)_B|}{|\Hom(A,C)|} \cdot \frac{|\Aut(B)|}{|\Aut(A)| |\Aut(C)|}.
	\]
	
	For any object $A$, 
	let
	\begin{align*}
		[\![A]\!]:=\frac{[A]}{|\Aut(A)|}.
	\end{align*}
	Then the Hall multiplication \eqref{eq:mult} can be reformulated to be
	\begin{align}
		[\![A]\!]\diamond [\![C]\!]=\sum_{[\![B]\!]}F_{AC}^B[\![B]\!],
	\end{align}
	which is the version of Hall multiplication used in \cite{Rin90}.

	
	



	\subsection{$\imath$Hall algebras}
	In this section, we always assume $\ca$ is a hereditary $\bfk$-linear abelian category, and shall construct a kind of Hall algebras for the category $\cc_{\varrho}(\ca)$.
	
	Let
	\[
	\sqq=\sqrt{q}.
	\]
	Generalizing \cite{LP16}, Lu defined a (twisted) semi-derived Ringel-Hall algebra of a category satisfying some conditions; see \cite[\S A.2]{LW19a}. By \eqref{eq:KXses}, Proposition \ref{cor:projfinite} and Theorem \ref{prop:Euler}, $\cc_{\varrho}(\ca)$ satisfies the conditions there, and we can define a semi-derived Ringel-Hall algebra for $\cc_{\varrho}(\ca)$. The $\imath$Hall algebra $\iH(\ca,\varrho)$  is by definition the twisted semi-derived Ringel-Hall algebra of $\cc_{\varrho}(\ca)$.
	For convenience, we recall it here briefly.
	
	Let $\ch(\cc_{\varrho}(\ca))$ be the Ringel-Hall algebra of $\cc_{\varrho}(\ca)$, i.e.,
	$$\ch(\cc_{\varrho}(\ca))=\bigoplus_{[M^\bullet]\in\Iso(\cc_{\varrho}(\ca))}\Q(\sqq)[M^\bullet],$$
	with the multiplication defined by (see \cite{Br})
	\[
	[M^\bullet]\diamond [N^\bullet]=\sum_{[M^\bullet]\in\Iso(\cc_{\varrho}(\ca))}\frac{|\Ext^1(M^\bullet,N^\bullet)_{L^\bullet}|}{|\Hom(M^\bullet,N^\bullet)|}[L^\bullet].
	\]
	Following \cite{LP16,LW19a,LW20a,LinP}, we consider the ideal $\cI$ of $\ch(\cc_{\varrho}(\ca))$ generated by
	\begin{align}
		\label{eq:ideal}
		&\Big\{ [M^\bullet]-[N^\bullet]\mid H(M^\bullet)\cong H(N^\bullet), \quad \widehat{\Im d_{M^\bullet}}=\widehat{\Im d_{N^\bullet}} \Big\}.
	\end{align}
	Consider the following multiplicatively closed subset $\cs$ of $\ch(\cc_{\varrho}(\ca))/\cI$:
	\begin{equation}
		\label{eq:Sca}
		\cs = \{ a[K^\bullet] \in \ch(\cc_{\varrho}(\ca))/\cI \mid a\in \Q(\sqq)^\times, K^\bullet\in \cc_{\varrho,ac}(\ca)\}.
	\end{equation}
	The semi-derived Ringel-Hall algebra of $\cc_{\varrho}(\ca)$ \cite{LW19a} is defined to be the localization
	$$\cs\cd\ch(\cc_{\varrho}(\ca)):= (\ch(\cc_{\varrho}(\ca))/\cI)[\cs^{-1}].$$
	We define the $\imath$Hall algebra $\iH(\ca,\varrho)$ (or denoted by $\tMH$) \cite[\S4.4]{LW19a} to be the $\Q(\sqq)$-algebra on the same vector space as $\utMH$ but with twisted multiplication given by
	\begin{align}
		\label{eqn:twsited multiplication}
		[M^\bullet]* [N^\bullet] =\sqq^{\langle \res(M^\bullet),\res(N^\bullet)\rangle} [M^\bullet]\diamond[N^\bullet].
	\end{align}

	For any $\alpha\in K_0(\ca)$,  there exist $X,Y\in\ca$ such that $\alpha=\widehat{X}-\widehat{Y}$. Define $[K_\alpha]:=[K_X]*[K_Y]^{-1}$. This is well defined, see, e.g.,
	\cite[\S 3.2]{LP16}. 
	
	For any complex $M^\bullet$ and acyclic complex $K_X$, we have
	\begin{align}
		\label{eq:KXM}
		[K_X]*[M^\bullet]=&\sqq^{\langle \varrho(X),\res(M^\bullet)\rangle-\langle X,\res(M^\bullet)\rangle}[K_X\oplus M^\bullet]
		\\
		\notag
		=&\sqq^{( \varrho(X),\res(M^\bullet))-( X,\res(M^\bullet))}[M^\bullet]*[ K_X].
	\end{align}
	It follows that $[K_\alpha]\; (\varrho(\alpha)=\alpha\in K_0(\ca))$ and $[K_\alpha]*[K_{\varrho(\alpha)}]$ are central in the algebra $\iH(\ca,\varrho)$, for any $\alpha\in K_0(\ca)$.
	
	
	The {\em quantum torus} $\widetilde{\ct}(\ca,\varrho)$ is defined to be the subalgebra of $\iH(\ca,\varrho)$ generated by $[K_\alpha]$, for $\alpha\in K_0(\ca)$. 

	\begin{proposition}  [cf. \cite{LinP}]
		\label{prop:hallbasis}
		The folllowing hold in $\iH(\ca,\varrho)$:
		\begin{enumerate}
			\item
			The algebra $\widetilde{\ct}(\ca,\varrho)$ is isomorphic to the group algebra of the abelian group $K_0(\ca)$.
			\item
			$\iH(\ca,\varrho)$ has an ($\imath$Hall) basis given by
			\begin{align*}
				\{[M]*[K_\alpha]\mid [M]\in\Iso(\ca), \alpha\in K_0(\ca)\}.
			\end{align*}
		\end{enumerate}
	\end{proposition}
	
	\begin{proof}
		With the help of Lemma \ref{lem:equiSES}, Theorem \ref{prop:Euler} and Proposition \ref{prop:Extp=0}, the proof of (2) is the same as \cite[Theorem~ 3.12]{LinP} (see also \cite[Theorem 3.20]{LP16} and \cite[Theorem 3.6]{LW20a}), and hence omitted here. Part (1) follows from (2).
	\end{proof}
	
	For any $M^\bullet=(M,d)\in \cc_\varrho(\ca)$, by definition, we have
	\begin{align}
		[M^\bullet]=[H(M^\bullet)\oplus K_{\Im d}]=\sqq^{\langle H(M^\bullet),\Im (\varrho(d))\rangle-\langle H(M^\bullet),\Im (d)\rangle}[H(M^\bullet)]*[K_{\Im d}].
	\end{align}
	
	\begin{remark}
		The semi-derived Hall algebras are defined for $\cc_{\mathbb{Z}_m}(\ca)$ for $m\geq1$, where $\ca$ is any arbitrary hereditary abelian category, and $\mod^{\rm nil}(\Lambda^\imath)$ for any $\imath$quiver algebra $\Lambda^\imath$; see \cite{LP16,LinP,LW19a,LW20a}. 
		The results in this subsection also hold for these semi-derived Hall algebras.
	\end{remark}

	Let $\widetilde{\ch}(\ca)$ be the twisted Ringel-Hall algebra, that is, the same vector space as $\ch(\ca)$ equipped with the twisted multiplication
	$$[A]* [B]=\sqq^{\langle A,B\rangle}[A]\diamond [B]$$
	for $[A],[B]\in\Iso(\ca)$.
	
	For $\beta,\gamma\in K_0(\ca)$, define $\beta\leq \gamma$ if there is $M\in\ca$ such that $\gamma-\beta=\widehat{M}$. This  endows $K_0(\ca)$ a partial order $\leq$. 
	For any $\gamma\in K_0(\ca)$, denote by $\iH(\ca,\varrho)_{\leq\gamma}$ (respectively, $\iH(\ca,\varrho)_{\gamma}$) the subspace of $\iH(\ca)$ spanned by elements from Proposition \ref{prop:hallbasis}~(2) for which $\widehat{M}\leq \gamma$ (respectively, $\widehat{M}= \gamma$) in $K_0(\ca)$.
	Then $\iH(\ca,\varrho)$ is a $K_0(\ca)$-graded linear space, i.e.,
	\begin{align}
		\label{eqn: graded linear MH}
		\iH(\ca,\varrho)=\bigoplus_{\gamma\in K_0(\ca)} \iH(\ca,\varrho)_\gamma.
	\end{align}
	Then $\iH(\ca,\varrho)$ is a filtered algebra with its associated graded algebra $\iH(\ca,\varrho)^{\gr}$ isomorphic to $\widetilde{\ch}(\ca)\otimes \widetilde{\ct}(\ca,\varrho)$; cf. \cite[Lemma 5.4]{LW19a}.

	Following \cite{LP16}, we also use $\cs\cd\widetilde{\ch}_{\Z_2}(\ca)$ to denote the twisted semi-derived Hall algebra of $\cc_{\Z_2}(\ca)$. Note that  $\cs\cd\widetilde{\ch}_{\Z_2}(\ca)$ is isomorphic to $\iH(\ca\times \ca,\swa)$; cf. Example \ref{ex:2-periodic}. Its quantum torus is denoted by $\widetilde{\ct}_{\Z_2}(\ca)$. The complex $M^\bullet$ is denoted by 
	$$\xymatrix{ M^0 \ar@<0.5ex>[r]^{d^0}& M^1 \ar@<0.5ex>[l]^{d^1}  },\quad d^1d^0=d^0d^1=0.$$
	
	For any object $X\in\ca$, we define
	\begin{align}
		\label{stalks}
		\begin{split}
			K_X:=&(\xymatrix{ X \ar@<0.5ex>[r]^{1}& X \ar@<0.5ex>[l]^{0}  }),\qquad \,\, K_X^*:=(\xymatrix{ X \ar@<0.5ex>[r]^{0}& X \ar@<0.5ex>[l]^{1}  }),
			\\
			C_X:=&(\xymatrix{ 0 \ar@<0.5ex>[r]& X \ar@<0.5ex>[l]  }),\qquad \quad C_X^*:=(\xymatrix{ X\ar@<0.5ex>[r]& 0 \ar@<0.5ex>[l]  })
		\end{split}
	\end{align}
	in $\cc_{\Z_2}(\ca)$
	
	Then the maps
	\begin{align}
		\label{eq:Rpm}
		R^+:\widetilde{\ch}(\ca)&\longrightarrow \cs\cd\widetilde{\ch}_{\Z_2}(\ca),\qquad\qquad R^-:\widetilde{\ch}(\ca)\longrightarrow \cs\cd\widetilde{\ch}_{\Z_2}(\ca),
		\\\notag
		[A]&\mapsto [C_A],\qquad\qquad\qquad\qquad\qquad [A]\mapsto [C_A^*],
	\end{align}
	are algebra embeddings.
	Moreover, we have the following triangular decomposition
	\begin{align}
		\label{triandecomp}
		\cs\cd\widetilde{\ch}_{\Z_2}(\ca)\cong \widetilde{\ch}(\ca)\otimes \ct_{\Z_2}(\ca)\otimes \widetilde{\ch}(\ca).
	\end{align}

	\subsection{An $\imath$Hall multiplication formula}
	
	Below, we present a fairly general multiplication formula in the $\imath$Hall algebra 	$\iH(\ca,\varrho)$ for $\cc_\varrho(\ca)$, generalizing \cite[Proposition 3.10]{LW20a}. 
	In concrete situations, the items appearing in RHS \eqref{Hallmult1} below are computable, and this makes Proposition~\ref{prop:iHallmult} useful and applicable.
	
	\begin{proposition}
		\label{prop:iHallmult}
		Let $\ca$ be a hereditary abelian category with an involution $\varrho$.
		For any $A,B\in\ca\subset \cc_\varrho(\ca)$, we have
		\begin{align}
			\label{Hallmult1}
			[A]*[B]=&
			\sum_{[L],[M],[N]\in\Iso(\ca)} \sqq^{\langle M,B\rangle-\langle M,L\rangle-\langle M,A\rangle+\langle M,N\rangle-\langle A,B\rangle}  q^{\langle N,L\rangle}\frac{|\Ext^1(N, L)_{M}|}{|\Hom(N,L)| }
			\\
			\notag
			&\cdot |\{s\in\Hom(A,\varrho B)\mid \Ker s\cong N, \coker s\cong \varrho L\}|\cdot [M]*[K_{\widehat{A}-\widehat{N}}]
			\\\label{Hallmult2}
			=&\sum_{[L],[M],[N],[X]\in\Iso(\ca)} \sqq^{\langle M,B\rangle-\langle M,L\rangle-\langle M,A\rangle+\langle M,N\rangle-\langle A,B\rangle} q^{\langle N,L\rangle} F_{N,L}^M F_{X,N}^AF_{\varrho L,X}^{\varrho B} 
			\\\notag
			&\cdot \frac{|\aut(L)|\cdot |\aut(N)|\cdot |\aut(X)|}{|\aut(M)|} \cdot[M]*[K_{\widehat{X}}]
		\end{align}
		in $\iH(\ca,\varrho)$.
	\end{proposition}

	\begin{proof}
		The proof of \cite[Proposition 3.10]{LW20a} for $\cc_{\Z_1}(\ca)$ (i.e., $\varrho=\Id$) can be easily generalized to $\cc_\varrho(\ca)$, we do not repeat it here.
	\end{proof}
	
	\begin{remark}
		By Example \ref{ex:2-periodic}, we can view the category of $2$-periodic complexes to be the category of $\varrho$-complexes. Let $\ca$ be a heredtiary category. Then we can obtain the Hall multiplication formula of $\cs\cd\widetilde{\ch}_{\Z_2}(\ca)$ by using Proposition \ref{prop:iHallmult}. 
	\end{remark}
	
	\section{Weighted projective lines and Coherent sheaves}
	\label{sec:WPL}
	
	In this section, we recall some basic facts on the category of coherent sheaves on a weighted projective line given by Geigle and Lenzing \cite{GL87}.

	\subsection{Weighted projective lines}
	Recall $\bfk=\F_q$.
	Fix a positive integer $\bt$ such that $2\leq \bt\leq q$ in the following.
	Let $\bp=(p_1,p_2,\dots,p_\bt)\in\Z_{+}^{\bt}$, and $L(\bp)$
	denote the rank one abelian group on generators $\vec{x}_1$, $\vec{x}_2$, $\cdots$, $\vec{x}_{\bt}$ with relations
	$p_1\vec{x}_1=p_2\vec{x}_2=\cdots=p_{\bt}\vec{x}_{\bt}$. We call $\vec{c}=p_i\vec{x}_i$ the {canonical element} of $L(\bp)$.
	Obviously, the polynomial ring $\bfk[X_1,\dots,X_{\bt}]$ is a $L(\bp)$-graded algebra by setting $\deg X_i=\vec{x}_i$, which is denoted by $\bS(\bp)$.
	
	Let $\ul{\bla}=\{\bla_1,\dots,\bla_{\mathbf{t}}\}$ be a collection of distinguished closed points (of degree one) on the projective line $\PL$, normalized such that $\bla_1=\infty$, $\bla_2=0$, $\bla_3=1$.
	Let $I(\bp,\ul{\bla})$ be the $L(\bp)$-graded ideal of $\bS(\bp)$ generated by
	$$\{X_i^{p_i}-(X_2^{p_2}-\bla_iX_1^{p_1})\mid 3\leq i\leq \bt \}.$$
	Then $\bS(\bp,\ul{\bla}):=\bS(\bp)/I(\bp,\ul{\bla})$ is an $L(\bp)$-graded algebra.
	For any $1\leq i\leq \bt$, denote by $x_i$ the image of $X_i$ in $\bS(\bp,\ul{\bla})$.
	
	The {weighted projective line} $\X:=\X_{\bp,\ul{\bla}}$ is the set of all non-maximal prime homogeneous ideals of $\bS:=\bS(\bp,\ul{\bla})$, which is also denoted by $\X_\bfk$ to emphasis the base field $\bfk$.
	For any homogeneous element $f\in \bS$, let $V_f:=\{\fp\in\X\mid f\in\fp\}$ and $D_f:=\X\setminus V_f$. The sets $D_f$ form a basis of the Zariski topology on $\X$.
	The pair $(\bp,\underline{\bla})$ (or just $\bp$ if $\underline{\bla}$ is obvious) is called the weight type of $\X_{\bp,\ul{\bla}}$, and $p_i$ is the weight of the point $\bla_i$.  
	In particular, the weighted projective line $\X_\bfk$ with $\bp=(1,1)$ is $\P^1_\bfk$.

	The classification of the closed points in $\X$ is provided in  \cite[Proposition 1.3]{GL87}. First, each $\bla_i$ corresponds to the prime ideal generated by $x_i$, called the {exceptional point}. Second, any other ideal is of the form $(f(x_1^{p_1},x_2^{p_2}))$, where $f\in \bfk [y_1,y_2]$ is an irreducible homogeneous polynomial in $y_1,y_2$, which is different from $y_1$ and $y_2$, called the {ordinary point}.
	
	For any closed point $x$ in $\X$, define
	\begin{align}
		\label{def:dx}
		d_x:=\begin{cases} 1, & \text{ if }x=\bla_i \text{ (exceptional)},
			\\
			\deg(f), & \text{ if }x=(f(x_1^{p_1},x_2^{p_2})) \text{ (ordinary)}.
		\end{cases}
	\end{align}

	\subsection{Coherent sheaves}
	
	The structure sheaf $\co:=\co_\X$ is the sheaf of $L(\bp)$-graded algebras over $\X$ associated to the presheaf $D_f\mapsto \bS_f$, where $\bS_f:=\{gf^{-l}\mid g\in \bS,l>0\}$ is the localization of $\bS$. Let $\Mod(\co_\X)$ be the category of sheaves of $L(\bp)$-graded $\co_\X$-modules.
	For any $\cm\in\Mod(\co_\X)$, we denote by $\cm(\vec{x})$ the shift of $\cm$ by $\vec{x}\in L(\bp)$.
	
	A sheaf $\cm\in\Mod(\co_\X)$ is called coherent if for any ${\blx}\in\X$, there exists a neighbourhood $U$ of ${\blx}$ and an exact sequence
	\begin{align}
		\bigoplus_{j=1}^n\co_\X(\vec{l}_j)|_{U}\longrightarrow\bigoplus_{i=1}^m\co_\X(\vec{k}_i)|_{U}\longrightarrow \cm|_{U}\longrightarrow 0.
	\end{align}
	We denote by $\coh(\X)$ the full subcategory of $\Mod(\co_\X)$ consisting of all coherent sheaves on $\X$. Then $\coh(\X)$ is a $\bfk$-linear hereditary, Hom-finite and Ext-finite abelian category; see \cite[Subsection 2.2]{GL87}. Moreover, by \cite[Proposition 1.8]{GL87}, $\coh(\X)$ is equivalent to the Serre quotient $\mod^{L(\bp)}(\bS)/ \mod^{L(\bp)}_0(\bS)$ of the category of finitely generated $L(\bp)$-graded
	$\bS$-modules modulo the category of finite length graded $\bS$-modules.

	A coherent sheaf $\cm\in\coh(\X)$ is called a \emph{vector bundle} of rank $n$ (or \emph{locally free sheaf}) if there is an open covering $\{U_i\}$ of $\X$, and an isomorphism
	$\bigoplus_{j=1}^n \co(\vec{l}_j)|_{U_i}\xrightarrow{\simeq} \cm|_{U_i}$ for some suitable $\vec{l}_j$ for each $U_i$; and is called a \emph{torsion sheaf} if it is a finite-length object in $\coh(\X)$.
	Let $\scrf$ be the full subcategory of $\coh(\X)$ consisting of all locally free sheaves, and $\scrt$ be the full subcategory consisting of all torsion sheaves. Then $\scrf$ and $\scrt$ are extension-closed. Moreover, $\scrt$ is a hereditary length abelian category.
	
	\begin{lemma}[\cite{GL87}]
		(1) For any sheaf $\cm\in\coh(\X)$, it can be decomposed as $\cm_t\oplus \cm_f$ where $\cm_f\in\scrf$ and $\cm_t\in \scrt$;
		
		(2) $\Hom(\cm_t,\cm_f)=0=\Ext^1(\cm_f,\cm_t)$ for any $\cm_f\in\scrf$ and $\cm_t\in\scrt$.
	\end{lemma}

	\subsection{Torsion sheaves}
	
	In order to describe the category $\scrt$, we shall introduce the representation theory of cyclic quivers.
	We consider the oriented cyclic quiver $C_n$ for $n\geq2$ with its vertex set $\Z_n=\{0,1,2,\dots,n-2,n-1\}$. 
	Let $C_1$ be the Jordan quiver, i.e., the quiver with only one vertex and one loop arrow.
	
	Denote by $\rep^{\rm nil}_\bfk(C_n)$ the category of finite-dimensional nilpotent representations of $C_n$ over the field $\bfk$.

	Then the structure of $\scrt$ is described in the following.
	\begin{lemma}[\cite{GL87}]
		\label{lem:isoclasses Tor}
		(1) The category $\scrt$ decomposes as a coproduct $\scrt=\coprod_{x\in\X} \scrt_{x}$, where $\scrt_{x}$ is the subcategory of torsion sheaves with support at $x$.
		
		(2) For any ordinary point $x$ of degree $d$, let $\bfk_{x}$ denote the residue field at $x$, i.e., $[\bfk_x:\bfk]=d$. Then $\scrt_{x}$ is equivalent to the category $\rep^{\rm nil}_{\bfk_{x}}(C_1)$.
		
		(3) For any exceptional point $\bla_i$ ($1\leq i\leq \bt$), the category $\scrt_{\bla_i}$ is equivalent to $\rep^{\rm nil}_{\bfk}(C_{p_i})$.
	\end{lemma}

	For any ordinary point $\blx$ of degree $d$, let $\pi_{\blx}$ be the prime homogeneous polynomial corresponding to $\blx$. Then multiplication by $\pi_{\blx}$ gives the exact sequence
	$$0\longrightarrow \co\stackrel{\pi_{\blx}}{\longrightarrow} \co(d\vec{c})\longrightarrow S_{\blx}\rightarrow0,$$
	where $S_{\blx}$ is the unique (up to isomorphisms) simple sheaf in the category $\scrt_{\blx}$. Then $S_{\blx}(\vec{l})=S_{\blx}$ for any $\vec{l}\in L(\bp)$.
	
	For any exceptional point $\bla_i$, multiplication by $x_i$ yields the short exact sequence
	$$0\longrightarrow \co((j-1)\vec{x}_i)\stackrel{x_i}{\longrightarrow} \co(j\vec{x}_i)\longrightarrow S_{ij}\rightarrow0,\text{ for } 1\le j\le p_i;$$
	where $\{S_{ij}\mid  j\in\Z_{p_i}\}$ is a complete set of pairwise non-isomorphic simple sheaves in the category $\scrt_{\bla_i}$ for any $1\leq i\leq \bt$. Moreover, $S_{ij}(\vec{l})=S_{i,j+l_i}$ for any $\vec{l}=\sum_{1\leq i\leq \bt} l_i\vec{x}_i$.

	In order to describe the indecomposable objects in $\scrt$, we need the following well-known results on the representation theory of cyclic quivers.
	\begin{itemize}
		\item[(1)] For $\rep^{\rm nil}_\bfk(C_1)$, the set of isomorphism classes of indecomposable objects is $\{S^{(a)}\mid a\in\N\}$, where $S=S^{(1)}$ is the unique simple representation, and $S^{(a)}$ is the unique indecomposable representation of length $a$. Define $S^{(\mu)}=\bigoplus_{i} S^{(\mu_i)}$ for any partition $\mu=(\mu_1\geq \cdots\geq \mu_l)$. Then any object in $\rep^{\rm nil}_\bfk(C_1)$ is isomorphic to $S^{(\mu)}$ for some partition $\mu$.
		\item[(2)]  For $\rep^{\rm nil}_\bfk(C_n)$, we have $n$ simple representations $S_j=S_j^{(1)}$ ($1\leq j\leq n$). Denote by $S_j^{(a)}$ the unique indecomposable representation with top $S_j$ and length $a$. Then $\{S_j^{(a)}\mid 1\le j\le n, a\in\N\}$ is  the set of isomorphism classes of indecomposable objects of $\rep^{\rm nil}_\bfk(C_n)$. For any $1\le j\le n$, define $S_j^{(\mu)}=\bigoplus_{i} S_j^{(\mu_i)}$ for any partition $\mu=(\mu_1\geq \cdots\geq \mu_l)$. 
		
	\end{itemize}
	
	Combining with Lemma \ref{lem:isoclasses Tor}, we can give a classification of indecomposable objects in $\scrt$. We denote by $S_{\blx}^{(a)}$ the unique indecomposable object of length $a$ in $\scrt_{\blx}$ for any ordinary point ${\blx}$; and denote by $S_{ij}^{(a)}$ the unique indecomposable object with top $S_{ij}$ and length $a$ in $\scrt_{\bla_i}$.

	\subsection{Automorphism groups}
	
	Let us recall 
	the automorphism groups $\Aut(\Coh(\X))$ and $\Aut(\X)$ introduced in \cite{LM} . By definition, the automorphism group $\Aut(\Coh(\X))$ of the category $\Coh(\X)$ consists of the isomorphism classes of $\bfk$-linear self-equivalences of $\Coh(\X)$. The automorphism group $\Aut(\X)$ is the subgroup of $\Aut(\Coh(\X))$ of automorphism $F$ fixing the structure sheaf, i.e., satisfying $F(\co)\cong \co$.

	\begin{proposition}[\cite{LM}]
		\label{prop:autogroupP1}
		The mapping $\pi: \Aut(\X)\rightarrow \Aut(\P^1_\bfk)$, where $\pi(F)=\ov{F}$ with $F(\scrt_x)=\scrt_{\ov{F}(x)}$, induces an isomorphism between the automorphism group $\Aut(\X)$ and the subgroup $\Aut_g(\X)$ of $\rm{PGL}(2,\bfk)$ of automorphisms of $\P_\bfk^1$ preserving weights.  
	\end{proposition}

	Set $\kappa(p)$ to count the number of exceptional points of weight $p$ and  define $\mathfrak{S}(\bp)$ as the direct product $\prod_p\mathfrak{S}_{\kappa(p)}$ of the symmetric groups $\mathfrak{S}_{\kappa(p)}$. 

	\begin{lemma}[\cite{LM}]
		If $\X$ has $\bt\geq3$ exceptional points, then $\Aut(\X)$ is isomorphic to a subgroup of the group $\mathfrak{S}(\bp)$ and in particular is finite.
	\end{lemma}	
	
	If $\bt\geq3$, we identify $\Aut(\X)$ as the subgroup of $\mathfrak{S}(\bp)$ in the following.

	For $\bp=(p_1,\dots,p_\bt)\in\Z_{+}^\bt$ with $p_i\geq 2$, let us consider the following star-shaped graph $\Gamma=\T_{p_1,\dots,p_\bt}$:
	
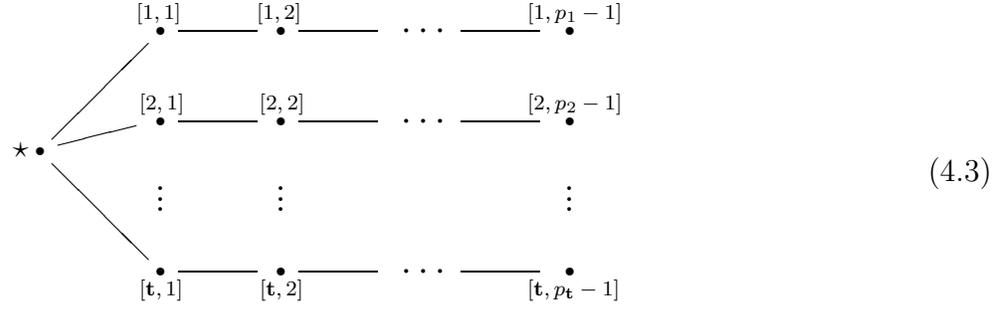
	\begin{center}\setlength{\unitlength}{0.8mm}
		\begin{equation}
			\label{star-shaped}
			\begin{picture}(110,30)(0,35)
				\put(0,40){\circle*{1.4}}
				\put(2,42){\line(1,1){16}}
				\put(20,60){\circle*{1.4}}
				\put(23,60){\line(1,0){13}}
				\put(40,60){\circle*{1.4}}
				\put(43,60){\line(1,0){13}}
				\put(60,58.5){\large$\cdots$}
				\put(70,60){\line(1,0){13}}
				\put(88,60){\circle*{1.4}}

				\put(3,41){\line(4,1){13}}
				\put(20,45){\circle*{1.4}}
				\put(23,45){\line(1,0){13}}
				\put(40,45){\circle*{1.4}}
				\put(43,45){\line(1,0){13}}
				\put(60,43.5){\large$\cdots$}
				\put(70,45){\line(1,0){13}}
				\put(88,45){\circle*{1.4}}
				
				\put(19,30){\Large$\vdots$}
				
				\put(39,30){\Large$\vdots$}
				
				\put(87,30){\Large$\vdots$}
				
				\put(2,38){\line(1,-1){16}}
				\put(20,20){\circle*{1.4}}
				\put(23,20){\line(1,0){13}}
				\put(40,20){\circle*{1.4}}
				\put(43,20){\line(1,0){13}}
				\put(60,18.5){\large$\cdots$}
				\put(70,20){\line(1,0){13}}
				\put(88,20){\circle*{1.4}}
				
				\put(-4.5,39){$\star$}
				
				\put(16,62){\tiny$[1,1]$}
				
				\put(36,62){\tiny$[1,2]$}
				\put(81,62){\tiny$[1,p_1-1]$}

				\put(16.5,47){\tiny$[2,1]$}
				\put(36.5,47){\tiny$[2,2]$}
				\put(81,47){\tiny$[2,p_2-1]$}

				\put(16.5,16){\tiny$[\bt,1]$}
				\put(36.5,16){\tiny$[\bt,2]$}
				\put(81,16){\tiny$[\bt,p_\bt-1]$}

				
			\end{picture}
		\end{equation}
		\vspace{1cm}
	\end{center}
	The set of vertices is denoted by $\II$. As marked in the graph, the central vertex is denoted by $\star$. Let $J_1,\dots,J_\bt$ be the branches, which are subdiagrams of type $A_{p_1-1},\dots,A_{p_\bt-1}$ respectively. Denote by $[i,j]$ the $j$-th vertex in the $i$-th branch. These examples include all finite-type Dynkin diagrams as well as the affine Dynkin diagrams of types $D_4^{(1)},E_6^{(1)},E_7^{(1)}$, and $E_8^{(1)}$.
	
	Let $Q_{\Gamma}$ be the quiver with underlying graph $\Gamma$ and all the orientations are from left to right. Note that the module category of $\bfk Q_{\Gamma}$ can be naturally embedded into the category of coherent sheaves $\coh(\X)$ of $\X=\X(\bp,\ul{\bla})$ for any parameter $\ul{\bla}$, namely, $\mod (\bfk Q_{\Gamma})$ can be identified with the perpendicular subcategory $\co(\vec{c})^{\perp}$ of $\coh(\X)$. More precisely, the simple modules at $[i,j]$ corresponds to simple sheaf $S_{ij}$, while the simple module 
	at $\star$ corresponds to the structure sheaf $\co$.

	For any involution $\varrho$ of the star-shaped graph $\Gamma$, $\varrho$ induces an involution on the module category of $\bfk Q_{\Gamma}$, still denoted by $\varrho$. We say that $\varrho$ can be lifted to $\coh(\X)$ if  there exists an involution $\tilde{\varrho}$ on $\coh(\X)$ for some parameter sequence $\ul{\bla}$, such that the restriction of $\tilde{\varrho}$ on $\co(\vec{c})^{\perp}$ coincides with $\varrho$.
	

	\begin{proposition}
		\label{prop:involution}
		Let $\bp=(p_1,\dots,p_\bt)\in\Z_+^\bt$ with $p_i\geq 2$, and $\varrho$ be a nontrivial involution of the star-shaped graph $\Gamma$. Then there exists a lifting of  $\varrho$ to $\coh(\X)$ for some weighted projective line $\X$ if and only if the number of $\varrho$-fixed branches $\leq 2$. In particular, any involution $\varrho$ of the star-shaped graph $\Gamma$ of finite and affine type can be lifted to $\coh(\X)$. 
	\end{proposition}
	
	\begin{proof}
		
		First, assume that there exists a lifting $\tilde{\varrho}$ of $\varrho$ to $\coh(\X)$. Note that $\tilde{\varrho}$ must fix the structure sheaf $\co$, hence belongs to $\Aut(\X)\subseteq\Aut(\P_\bfk^1)$. It follows from Proposition \ref{prop:autogroupP1} that the number of $\varrho$-fixed branches $\leq 2$ since $\varrho$ is a nontrivial involution.
		
		Conversely, let $G=\diag(1,-1)\in\rm{PGL}(2,\bfk)=\aut(P^1_\bfk)$. Then $G$ has two fixed points $0$ and $\infty$ of $\mathbf{P}^1_{\mathbf{k}}$, and exchanges the vertices $\lambda$ and $-\lambda$ for any $\lambda\in\bfk\setminus \{0\}\subset \P^1_\bfk$. Assume $\varrho$ exchanges the branches $J_{2i-1}$ with $J_{2i}$ for any $1\leq i\leq r$ and fix the other branches. Then  $p_{2i-1}=p_{2i}$ for any$1\leq i\leq r$. We can choose distinct points $\lambda_{i}\in\P^1_\bfk$ (and attach the weights $p_i$ to them) 
		such that $\lambda_{2i-1}=-\lambda_{2i}$ for $1\leq i\leq r$.
		Since the number of $\varrho$-fixed branches is less than or equal to 2, we can attach the associated weights to the points $0$ and/or $\infty$. 
		In this way we obtain a weighted projective line $\mathbb{X}(\bp,\ul{\bla})$ with an involution $\tilde{\varrho}=G$, which exchanges the weighted points $\lambda_{2i-1}$ and $\lambda_{2i}$ for $1\leq i\leq r$ and fixes the remaining weighted points. This yields an involution on the category of coherent sheaves over $\mathbb{X}(\bp,\ul{\bla})$, whose restriction on $\co(\vec{c})^{\perp}$ coincides with $\varrho$. Then we are done.
	\end{proof}

	\subsection{Grothendieck groups and Euler forms}
	The Grothendieck group $K_0(\coh(\X))$ of $\coh(\X)$ satisfies
	\begin{align}
		K_0(\coh(\X))\cong \Big(\Z\widehat{\co}\oplus \Z\widehat{\co(\vec{c})}\oplus \bigoplus_{i,j}\Z\widehat{S_{ij}}\Big)/I,
	\end{align}
	where $I$ is the subgroup generated by $\{\sum_{j=1}^{p_i} \widehat{S_{ij}} +\widehat{\co} -\widehat{\co(\vec{c})}\mid i=1,\dots,\bt\}$; see \cite{GL87}.
	Let $\de:=\widehat{\co(\vec{c})}-\widehat{\co}= \sum_{j=1}^{p_i} \widehat{S_{ij}}$.
	The Euler form $\langle-,-\rangle$ on $K_0(\coh(\X))$ is described in \cite{Sch04}, which is given by
	\begin{align}
		&\langle \widehat{\co},\widehat{\co}\rangle=1,\quad \langle \widehat{\co},\delta\rangle=1,\quad  \langle\delta,\widehat{\co}\rangle=-1,
		\\
		&\langle \delta,\delta\rangle=0,\quad \langle \widehat{S_{ij}},\delta\rangle=0,\quad  \langle \delta, \widehat{S_{ij}}\rangle=0,
		\\
		&\langle \widehat{\co}, \widehat{S_{ij}}\rangle =\begin{cases} 1& \text{ if } j=p_i, \\ 0& \text{ if }j\neq p_i , \end{cases}
		\qquad \langle \widehat{S_{ij}}, \widehat{\co}\rangle =\begin{cases} -1& \text{ if } j=1, \\ 0& \text{ if }j\neq 1 ,\end{cases}
		\\
		&\langle\widehat{ S_{ij}}, \widehat{S_{i',j'}}\rangle =\begin{cases} 1& \text{ if } i=i', j=j',
			\\ -1& \text{ if }i=i', j\equiv j'+1 (\mod p_i),
			\\
			0& \text{ otherwise. }\end{cases}
	\end{align}
	
	Recall the root system $\cR$ of $L\fg$. Then there is a natural isomorphism of $\Z$-modules $K_0(\coh(\X))\cong \cR$ given as below for $1\leq i\leq \bt,\;1\le j\le p_i-1, $ $r\in\N$ and $l\in\Z$:
	\begin{align*}
		&\widehat{S_{ij}} \mapsto \alpha_{ij}, 
		&\widehat{S_{i,0}}\mapsto \de-\sum_{j=1}^{p_i-1}\alpha_{ij},\quad
		&\widehat{S_{i,0}^{(p_i-1)}}\mapsto \de-\alpha_{i1}, 
		&\widehat{S_{i,0}^{(rp_i)}}\mapsto r\de,\quad
		&\widehat{\co(l\vec{c})}\mapsto \alpha_\star+l\de.
	\end{align*}
	Under this isomorphism, the symmetric Euler form on $K_0(\coh(\X))$ coincides with the Cartan form on $\cR$. So we always identify $K_0(\coh(\X))$ with $\cR$ via this isomorphism in the following.
	

	\section{Quantum groups and $\imath$quantum groups}
	\label{sec:QG-iQG}
	
	In this section, we review the basic materials on (affine) quantum groups and $\imath$quantum groups of quasi-split type, especially their Drinfeld type presentations.
	
	\subsection{Quantum groups}
	\label{subsec:QG}
	Denote, for $r,m \in \N$,
	\[
	[r]:=[r]_v=\frac{v^r-v^{-r}}{v-v^{-1}},
	\quad
	[r]!:=[r]^!_v=\prod_{i=1}^r [i], \quad \qbinom{m}{r} =\frac{[m][m-1]\ldots [m-r+1]}{[r]!}.
	\]
	For $A, B$ in a $\Q(v)$-algebra, we shall denote $[A,B]_{v^a} =AB -v^aBA$, and $[A,B] =AB - BA$. 
	
	Let $\Gamma$ be a graph with the vertex set $\I$.
	Let $n_{ij}$ be the number of edges connecting  vertices $i$ and $j$. Let $C=(c_{ij})_{\I\times \I}\in M_{\I\times \I}(\Z)$ be the symmetric generalised Cartan matrix of the graph $\Gamma$, defined by
	$$c_{ij}=2\delta_{ij}-n_{ij}.$$
	Let $\fg$ be the corresponding Kac-Moody Lie algebra, and $\tU_v( \mathfrak{g})$ its quantized enveloping algebra. 
	Then $\tU := \tU_v(\fg)$ is defined to be the $\Q(v)$-algebra generated by $E_i,F_i, \tK_i,\tK_i'$, $i\in \I$, where $\tK_i, \tK_i'$ are invertible, subject to the following relations:
	\begin{align}
		[E_i,F_j]= \delta_{ij} \frac{\tK_i-\tK_i'}{v-v^{-1}},  &\qquad [\tK_i,\tK_j]=[\tK_i,\tK_j']  =[\tK_i',\tK_j']=0,
		\label{eq:KK}
		\\
		\tK_i E_j=v^{c_{ij}} E_j \tK_i, & \qquad \tK_i F_j=v^{-c_{ij}} F_j \tK_i,
		\label{eq:EK}
		\\
		\tK_i' E_j=v^{-c_{ij}} E_j \tK_i', & \qquad \tK_i' F_j=v^{c_{ij}} F_j \tK_i',
		\label{eq:K2}
	\end{align}
	and the quantum Serre relations, for $i\neq j \in \I$,
	\begin{align}
		& \sum_{r=0}^{1-c_{ij}} (-1)^r \left[ \begin{array}{c} 1-c_{ij} \\r \end{array} \right]  E_i^r E_j  E_i^{1-c_{ij}-r}=0,
		\label{eq:serre1} \\
		& \sum_{r=0}^{1-c_{ij}} (-1)^r \left[ \begin{array}{c} 1-c_{ij} \\r \end{array} \right]  F_i^r F_j  F_i^{1-c_{ij}-r}=0.
		\label{eq:serre2}
	\end{align}
	Note that $\tK_i \tK_i'$ are central in $\tU$ for all $i$.

	
	
	Following Lusztig, for any $i\in\I$, there exists an algebra isomorphism $T_i:\tU\rightarrow\tU$ via the following rules (see \cite[Proposition 6.20]{LW22}):
	\begin{align*}
		T_i(E_i)=&-F_i(K'_i)^{-1},\qquad T_i(F_i)=-K_i^{-1}E_i,\\
		T_i(K_j)=&K_j(K_{i})^{-c_{ij}},\qquad T_i(K'_j)=K'_j(K_i')^{-c_{ij}},
		\\
		T_i(E_j)=&\sum_{r+s=-c_{ij}} (-1)^r v^{-r} E_i^{(s)} E_j E_i^{(r)}, \text{ for }i\neq j,
		\\
		T_i(F_j)=&\sum_{r+s=-c_{ij}} (-1)^r v^{-r} F_i^{(r)} F_j F_i^{(s)}, \text{ for }i\neq j,
	\end{align*}
	where
	$$E_i^{(r)}=\frac{E_i^r}{[r]!},\qquad F_i^{(r)}=\frac{F_i^r}{[r]!},\qquad \forall r\geq0.$$

	\subsection{Quantum loop algebras}

	\label{sec:QLA}
	
	Let $\fg$ be a Kac-Moody algebra, and $\Lg$ be its loop algebra. 
	The Drinfeld double quantum loop algebra $\tUD_v(\Lg)$
	is the $\Q(v)$-algebra generated by $x_{i k}^{\pm}$, $h_{i l}$ and the invertible elements $K_i$, $K_i'$, $C$, $\widetilde{C}$ for $i\in\II$, $k\in\Z$, $l\in\Z\backslash\{0\}$, subject to the following relations: 
	\begin{align}
		\label{DR1}
		C, \widetilde{C} \text{ are  central},\quad [K_i,K_j]  =& [K_i,K_j']=[K_i',K_j']=0,
		\\	\label{DR1a} [K_i,h_{j l}] =0=[K_i',h_{jl}],& 
		\\
		\label{DR2}
		[h_{il},h_{jm}]=0=[h_{i,-l},h_{j,-m}],\qquad [h_{il},h_{j,-m}] &= 
		\delta_{l, m} \frac{[l c_{ij}]}{l} \frac{C^l -\widetilde{C}^{l}}{v^{-1} -v},\quad \text{ for }l,m>0,
		\\
		\label{DR3}
		K_ix_{jk}^{\pm} =v^{\pm c_{ij}} x_{jk}^{\pm}K_i,\qquad K_i'x_{jk}^{\pm} & =v^{\mp c_{ij}} x_{jk}^{\pm}K_i',
		\\
		\label{DR4}
		[h_{i l},x_{j k}^{\pm}] =&\begin{cases}\pm\frac{[lc_{ij}]}{l} C^{\frac{l\mp|l|}{2}} x_{j,k+l}^{\pm},&\text{ if }l>0,
			\\
			\mp\frac{[lc_{ij}]}{l} \widetilde{C}^{-\frac{l\mp|l|}{2}} x_{j,k+l}^{\pm},&\text{ if }l<0,
		\end{cases}
		\\
		\label{DR5}
		[x_{i k}^+,x_{j t}^-] =\delta_{ij}& {(C^{-t} K_i\psi_{i,k+t} - \widetilde{C}^{k} K_i' \varphi_{i,k+t})}, 
		\\
		\label{DR6}
		x_{i,k+1}^{\pm} x_{j,t}^{\pm}-v^{\pm c_{ij}} x_{j,t}^{\pm} x_{i,k+1}^{\pm} &=v^{\pm c_{ij}} x_{i,k}^{\pm} x_{j,t+1}^{\pm}- x_{j,t+1}^{\pm} x_{i,k}^{\pm},
		\\
		\label{DR7}
		\Sym_{k_1,\dots,k_r}\sum_{s=0}^{r} (-1)^s \qbinom{r}{s}
		x_{i,k_1}^{\pm}\cdots
		x_{i,k_s}^{\pm} x_{j,t}^{\pm} & x_{i,k_s+1}^{\pm} \cdots x_{i,k_r}^{\pm} =0, \text{ for } r= 1-c_{ij}\; (i\neq j),
	\end{align}
	where
	$\Sym_{k_1,\dots,k_r}$ denotes the symmetrization with respect to the indices $k_1,\dots,k_r$; $\psi_{i,k}$ and $\varphi_{i,k}$ are defined by the following functional equations:
	\begin{align}
		\label{exp h+}
		1+ \sum_{m\geq 1} (v-v^{-1})\psi_{i,m}u^m &=  \exp\Big((v -v^{-1}) \sum_{m\ge 1}  h_{i,m}u^m\Big),
		\\
		\label{exp h-}
		1+ \sum_{m\geq1 }(v-v^{-1})\varphi_{i, -m}u^{-m} &= \exp \Big((v -v^{-1}) \sum_{m\ge 1} h_{i,-m}u^{-m}\Big).
	\end{align}
	
	The quantum loop algebra $\UD_v(\Lg)$ can be described as a quotient algebra of $\tUD_v(\Lg)$ modulo the ideal generated by
	$$C\widetilde{C}-1,\qquad K_iK_i'-1, \,\,\forall i\in\II.$$
	
	Assume $\fg$ is a simply-laced simple Lie algebra in the following.
	Let $\widehat{\fg}=\fg[t,t^{-1}]\oplus \C c$ be the affine Kac-Moody algebra.
	Then $\widehat{\fg}=\Lg$.
	It was stated by Drinfeld \cite{Dr88} and proved by Beck \cite{Be94} 
	that $\UD_v(\Lg)$ is isomorphic to $\U_v(\widehat{\fg})$ (by adding two central generators $C^{\pm\frac{1}{2}}$), which is  called the Drinfeld new presentation of $\U_v(\widehat{\fg})$.  (Comparing with \cite{Be94}, we slightly modify the imaginary root vectors $h_{i,l},\psi_{i,m},\varphi_{i,-m}$, and omit a degree operator $D$ in the version of ${}^{\text{Dr}}\U_v(\Lg)$ above.)

	Similarly, one can prove that $\tUD_v(\Lg)$ is isomorphic to $\tU_v(\widehat{\fg})$, which is  called the Drinfeld new presentation of $\tU_v(\widehat{\fg})$ (and also denoted by $\tUD_v(\widehat{\fg})$).
	This isomorphism is explicitly described by the root vectors of $\tU_v(\widehat{\fg})$  defined in the following. Define a sign function
	\[
	o(\cdot): \II \longrightarrow \{\pm 1\}
	\]
	such that $o(i) o(j)=-1$ whenever $c_{ij} <0$ (there are clearly exactly 2 such functions).
	For $i\in\II$ and $k\in\Z$, let
	\begin{align*}
		x_{i,k}^-=o(i)^{k} T_{\omega_i}^k(F_i),\qquad x_{i,k}^+=o(i)^{k}  T_{\omega_i}^{-k}(E_i).
	\end{align*}
	For $m>0$, define
	\begin{align*}
		\psi_{i,m}=o(i)^{m}\tK_i^{-1}[T_{\omega_i}^{-m}(E_i),F_i], \qquad \varphi_{i,-m}=o(i)^{m}(\tK_i')^{-1}[T_{\omega_i}^{-m}(F_i),E_i];
	\end{align*}
	and then define $h_{i,\pm m}$ via \eqref{exp h+}--\eqref{exp h-}.

	\subsection{Quantum symmetric pairs and $\imath$quantum groups}
	
	For a  (generalized) Cartan matrix $C=(c_{ij})$, let $\Aut(C)$ be the group of all permutations $\tau$ of the set $\I$ such that $c_{ij}=c_{\tau i,\tau j}$. An element $\tau\in\Aut(C)$ is called an \emph{involution} if $\tau^2=\Id$.
	
	Let $\tau$ be an involution in $\Aut(C)$. We define $\widetilde{\bU}^\imath$ to be the $\Q(v)$-subalgebra of $\tU$ generated by
	\[
	B_i= F_i +  E_{\tau i} \tK_i',
	\qquad \K_i = -v^{c_{i,\tau i}}\tK_i \tK_{\tau i}', \quad \forall i \in \I.
	\]

	We only consider symmetric Cartan matrices $C=(c_{ij})$ with $c_{ij}\in\{2,0,-1,-2\}$ in this paper, and assume the involution $\tau$ with $\tau i\neq i$ only if $c_{i,\tau i}=0$. In this situation,  $\tUi$ is generated by $B_i,\K_i$ ($i\in\I$) subject to
	\begin{align}
		\K_i\K_i^{-1} =\K_i^{-1}\K_i=1, & \quad \K_i\K_j=\K_j \K_i, \quad
		\K_j B_i=v^{c_{i,\tau j} -c_{ij}} B_i \K_j,
		\\
		B_iB_j -B_j B_i&=0, \qquad\qquad\qquad \text{ if } c_{ij}=0 \text{ and }\tau i\neq j,
		\label{eq:S1} \\
		\sum_{s=0}^{1-c_{ij}} (-1)^s \qbinom{1-c_{ij}}{s}& B_i^{s}B_jB_i^{1-c_{ij}-s}=0, \quad \text{ if } j \neq \tau i\neq i,
		\label{eq:S6}
		\\
		B_{\tau i}B_i -B_i B_{\tau i}& =   \frac{\K_i -\K_{\tau i}}{v-v^{-1}},
		\quad \text{ if } \tau i \neq i \text{ and }c_{\tau i,i}=0,
		\label{relation5}
		\\
		B_i^2 B_j -[2] B_i B_j B_i +B_j B_i^2 &= - v^{-1}  B_j \K_i,  \qquad \text{ if }c_{ij}=-1 \text{ and }\tau i=i,
		\label{eq:S2} \\
		\sum_{r=0}^3 (-1)^r \qbinom{3}{r} B_i^{3-r} B_j B_i^{r} &= -v^{-1} [2]^2 (B_iB_j-B_jB_i) \K_i,  \label{eq:S3}
		\\\notag
		&  \qquad\qquad\qquad\qquad\text{ if }c_{ij}=-2 \text{ and }\tau i=i.
	\end{align}
	
	
	For $\mu = \mu' +\mu''  \in \Z \I := \oplus_{i\in \I} \Z \alpha_i$,  define $\K_\mu$ such that
	\begin{align}
		\K_{\alpha_i} =\K_i, \quad
		\K_{-\alpha_i} =\K_i^{-1}, \quad
		\K_{\mu} =\K_{\mu'} \K_{\mu''}.
	\end{align}
	The algebra $\tUi$ is endowed with a filtered algebra structure
	\begin{align}  \label{eq:filt1}
		\widetilde{\U}^{\imath,0} \subset \widetilde{\U}^{\imath,1} \subset \cdots \subset \widetilde{\U}^{\imath,m} \subset \cdots
	\end{align}
	by setting 
	\begin{align}  \label{eq:filt}
		\widetilde{\U}^{\imath,m} =\Q(v)\text{-span} \{ B_{i_1} B_{i_2} \ldots B_{i_r} \K_\mu \mid \mu \in \N\I, i_1, \ldots, i_r \in \I, r\le m \}.
	\end{align}
	Note that
	\begin{align}  \label{eq:UiCartan}
		\widetilde{\U}^{\imath,0} =\bigoplus_{\mu \in \N\I} \Q(v) \K_\mu,
	\end{align}
	is the $\Q(v)$-subalgebra generated by $\K_i$ for $i\in \I$.
	Then, according to a basic result of Letzter and Kolb on quantum symmetric pairs adapted to our setting of $\tUi$ (cf. \cite{Let02, Ko14}), the associated graded $\gr \tUi$ with respect to \eqref{eq:filt1}--\eqref{eq:filt} can be identified with
	\begin{align}   \label{eq:filter}
		\gr \tUi \cong \U^- \otimes \Q(v)[\K_i^\pm | i\in \I],
		\qquad \overline{B_i}\mapsto F_i,  \quad
		\overline{\K}_i \mapsto \K_i \; (i\in \I).
	\end{align}
	
	\subsection{A Drinfeld type presentation of quasi-split affine $\imath$quantum groups}
	
	Let $C=(c_{ij})_{i,j\in \II}$ be a generalized Cartan matrix (GCM) of a {\em simply-laced} Kac-Moody algebra $\fg$. 
	Let $\tau$ be an involution of $C$. 
	
	The $\imath$quantum loop algebra $\tUiD$ of quasi-split type is the $\Q(v)$-algebra  generated by $\K_{i}^{\pm1}$, $C^{\pm1}$, $H_{i,m}$, $\Theta_{i,l}$ and $B_{i,l}$, where  $i\in \II$, $m \in \Z_{+}$, $l\in\Z$, subject to some relations.
	In order to give the explicit definition of $\tUiD$,  we introduce some shorthand notations below.
	
	Let $k_1, k_2, l\in \Z$ and $i,j \in \II$. Set
	\begin{align}
		\begin{split}
			S(k_1,k_2|l;i,j)
			&=  B_{i,k_1} B_{i,k_2} B_{j,l} -[2] B_{i,k_1} B_{j,l} B_{i,k_2} + B_{j,l} B_{i,k_1} B_{i,k_2},
			\\
			\SS(k_1,k_2|l;i,j)
			&= S(k_1,k_2|l;i,j)  + \{k_1 \leftrightarrow k_2 \}.
			\label{eq:Skk}
		\end{split}
	\end{align}
	Here and below, $\{k_1 \leftrightarrow k_2 \}$ stands for repeating the previous summand with $k_1, k_2$ switched if $k_1\neq k_2$, so the sums over $k_1, k_2$ are symmetric.
	We also denote
	\begin{align}
		\begin{split}
			R(k_1,k_2|l; i,j)
			&=   \K_i  C^{k_1}
			\Big(-\sum_{p\geq0} v^{2p}  [2] [\Theta _{i,k_2-k_1-2p-1},B_{j,l-1}]_{v^{-2}}C^{p+1}
			\label{eq:Rkk} \\
			&\qquad\qquad -\sum_{p\geq 1} v^{2p-1}  [2] [B_{j,l},\Theta _{i,k_2-k_1-2p}]_{v^{-2}} C^{p}
			- [B_{j,l}, \Theta _{i,k_2-k_1}]_{v^{-2}} \Big),
			\\
			\RR(k_1,k_2|l; i,j) &= R(k_1,k_2|l;i,j) + \{k_1 \leftrightarrow k_2\}.
		\end{split}
	\end{align}
	Sometimes, it is convenient to rewrite part of the summands in \eqref{eq:Rkk} as
	\begin{align*}
		&-\sum_{p\geq 1} v^{2p-1}  [2] [B_{j,l},\Theta _{i,k_2-k_1-2p}]_{v^{-2}} C^{p}
		- [B_{j,l}, \Theta _{i,k_2-k_1}]_{v^{-2}}\\
		=&
		-\sum_{p\geq 0} v^{2p-1}  [2] [B_{j,l},\Theta _{i,k_2-k_1-2p}]_{v^{-2}} C^{p}
		+v^{-2}[B_{j,l}, \Theta _{i,k_2-k_1}]_{v^{-2}}.
	\end{align*}
	
	
	\begin{definition} [$\imath$Quantum loop algebras of quasi-split type]
		\label{def:iDRA1}
		Let $\tUiD=\tUiD(\Lg)$ be the
		$\mathbb{Q}(v)$-algebra generated by
		the elements $B_{i,l}$, $H_{i,k}$, $\K_i^{\pm1}$, $C^{\pm1}$, where $i\in\II$, $l\in\Z$ and $k,m>0$, subject
		to the following relations:
		\begin{align}
			&C^{\pm1} \text{ is central,} \qquad [\K_i,\K_j]=0, \quad [\K_i,H_{j,k}]=0,\quad \K_iB_{j,l}=v^{c_{\tau i,j}-c_{ij}} B_{j,l} \K_i,\label{qsiDR0}
			\\\label{qsiDR1}
			& [H_{i,k},H_{j,m}]=0,
			\\\label{qsiDR2}
			&[H_{i,k},B_{j,l}]=\frac{[kc_{ij}]}{k} B_{j,l+k}-\frac{[kc_{\tau i,j}]}{k} B_{j,l-k}\K_{k\delta},
			\\\label{qsiDR3}
			&[B_{i,k},B_{j,l+1}]_{v}-v[B_{i,k+1},B_{j,l}]_{v^{-1}} =0, \quad \text{ if } c_{ij}=-1, \quad j \neq  \tau i,
			\\\label{qsiDR4}
			&[B_{i,k},B_{\tau i,l}]= \K_{\tau i} \K_{l \delta} \TH_{i,k-l}- \K_{i} \K_{k\delta} \Theta_{\tau i,l-k},
			\qquad \text{ if }c_{i,\tau i}=0,
			\\\label{qsiDR5}
			&[B_{i,k},B_{i,l+1}]_{v^{-2}}-v^{-2}[B_{i,k+1},B_{i,l}]_{v^{2}}=v^{-2} \Theta_{i,l-k+1}\K_{k\delta+\alpha_i} -  v^{-4} \Theta_{i,l-k-1}\K_{(k+1)\delta+\alpha_i}
			\\
			& \quad \qquad\qquad\quad +v^{-2} \Theta_{i,k-l+1}\K_{l\delta+\alpha_i} -v^{-4} \Theta_{i,k-l-1}\K_{(l+1)\delta+\alpha_i}, \; \qquad\qquad\text{ if } i =\tau i,
			\notag  \\
			&  [B_{i,k}, B_{j,l}] =0,
			\qquad \text{ if }c_{ij}=0 \text{ and }\tau i\neq j, \label{qsiDR7}
			\\\label{qsiDR8}
			&\SS(k_1,k_2|l;i,j)=\RR(k_1,k_2|l;i,j), \qquad\qquad\quad\text{ if }c_{ij}=-1 \text{ and }\tau i=i, 
			\\
			&\SS(k_1,k_2|l;i,j)
			= 0,\qquad\qquad\qquad\qquad\qquad
			\text{ if }c_{i,j}=-1,  \tau i\neq i. \label{qsiDR9}
		\end{align}
	\end{definition}
	
	In the definition above, $H_{i,m}$ are related to $\Theta_{i,m}$ by the following equation for generating functions in $u$: 
	\begin{align}
		\label{exp}
		1+ \sum_{m\geq 1} (v-v^{-1})\Theta_{i,m} u^m  = \exp\Big( (v-v^{-1}) \sum_{m\geq 1} H_{i,m} u^m \Big).
	\end{align}
	
	\begin{lemma}[\text{\cite[Lemmas 4.7,4.8]{LWZ23}}] \label{lem:equiv}
		(1) The relation \eqref{qsiDR1}  (for $k,m>0$) is equivalent to 
		\begin{align}
			[\Theta_{i,k},\TH_{j,m}]=0.
		\end{align}
		(2) The relation \eqref{qsiDR2} (for $k>0,l\in\Z$) is equivalent to 
		\begin{align}
			[\TH_{i,k}, B_{j,l}] &+v^{c_{i,j}-c_{\tau i,j}}[\TH_{i,k-2},B_{j,l}]_{v^{2(c_{\tau i,j}-c_{i,j})}}C
			-v^{c_{i,j}}[\TH_{i,k-1}, B_{j,l+1}]_{v^{-2c_{ i,j}}} 
			\notag \\
			&- v^{-c_{\tau i, j}}[\TH_{i,k-1},B_{j,l-1}]_{v^{2c_{\tau i,j}}}C
			=0.
			\label{THBij} 
		\end{align}
		
	\end{lemma}
	
	By \cite{LW21b,LWZ23}, there is a $\Q(v)$-algebra isomorphism $\Phi:\tUiD(\Lg)\rightarrow\tUi(\widehat{\fg})$ if $\fg$ is of type ADE; cf. \cite{Z22}.

	
	\section{$\imath$Hall algebras for cyclic quivers}
	\label{sec:cyclic}

	Let $C_n$ be the cyclic quiver with $n$ vertices. 
	Recall that $\rep^{\rm nil}_\bfk(C_n)$ is the category of finite-dimensional nilpotent representations of $C_n$. For the Grothendieck group $K_0(\rep^{\rm nil}_\bfk(C_n))$, we denote $\alpha_j=\widehat{S_j}$ and $\de=\sum_{j=1}^n\widehat{S_j}$ by abusing notations. For any positive real root $\beta$ of $\widehat{\mathfrak{sl}}_{n}$, by Gabriel-Kac Theorem, we denote by $M(\beta)$ the unique
	(up to isomorphism) indecomposable $\bfk C_n$-module with its class $\beta$ in $K_0(\rep^{\rm nil}_\bfk(C_n))$.
	
	\subsection{Affine $\imath$quantum groups of $\mathfrak{sl}_n$}

	Let $\tUi:=\tUi_v(\widehat{\mathfrak{sl}}_n)$ be the universal affine $\imath$quantum group of type $A_{n-1}$, and $\tUiD$ be its Drinfeld presentation throughout this subsection.
	By \cite[Theorem 2.16]{LW21b}, there exists an algebra isomorphism $\Phi:\tUiD\rightarrow\tUi$. We choose the sign function $o(\cdot)$ such that $o(j)=(-1)^j$ for any $1\leq j\leq n-1$ throughout this paper.
	
	Let $\iH(\bfk C_n,\Id)$ be the $\imath$Hall algebra of $\cc_{\Z_1}(\rep_\bfk^{\rm nil}(C_n))$. 
	For $n\geq2$, \cite[Theorem 9.6]{LW20a} gives an algebra monomorphism:
	\begin{align}
		\label{eq:isoCnsln}
		\widetilde{\psi}_{C_n}: \tUi_{\sqq}\longrightarrow &\iH(\bfk C_n,\Id)
		\\
		B_j \mapsto& \frac{-1}{q-1}[S_{j}],
		\quad  \K_j \mapsto [K_{S_j}], \text{ for }0\leq j\leq n-1. 
	\end{align}
	
	Define
	\begin{align}\label{the map Psi A}
		\Omega_{C_n}:=\widetilde{\psi}_{C_n}\circ \Phi:\tUiD_{\sqq}\longrightarrow \iH(\bfk C_n,\Id).
	\end{align}
	Then $\Omega_{C_n}$ sends:
	\begin{align}
		\label{eq:HaDrA1}\K_j\mapsto [K_{S_j}],\qquad C\mapsto [K_\de],\qquad
		B_{j,0}\mapsto \frac{-1}{q-1}[S_j], \qquad \forall \,\,1\leq j\leq n-1.
	\end{align}
	For any $1\leq j\leq n-1$, $l\in\Z$ and $r\geq 1$, we define
	\begin{align}
		\label{def:haBa}
		\haB_{j,l}:=(1-q)\Omega_{C_n}(B_{j,l}),\qquad
		\widehat{\Theta}_{j,r}:= \Omega_{C_n}(\Theta_{j,r}),\qquad
		\widehat{H}_{j,r}:= \Omega_{C_n}(H_{j,r}).
	\end{align}
	In particular, $\haB_{j,0}=[S_j]$ for any $1\leq j\leq n-1 $. 
	
	The morphism $\Omega_{C_n}:\tUiD_{\sqq}\longrightarrow \iH(\bfk C_n,\Id)$ and root vectors $\haB_{j,l},
	\widehat{\Theta}_{j,r},
	\widehat{H}_{j,r}$ are studied indepth in \cite{LR21}.

	\subsection{Affine quantum groups of $\mathfrak{sl}_n$}
	For a quiver $Q$, let $\cs\cd\widetilde{\ch}_{\Z_2}(\bfk Q)$ be the twisted semi-derived Ringel-Hall algebra of $\rep^{\rm nil}_\bfk(Q)$. The following proposition is a generalization of Bridgeland's result \cite{Br} (without assuming $Q$ to be acyclic). Let $\tU$ be the corresponding quantum group. 
	
	\begin{proposition}[\text{\cite[Theorem H]{LW19a}}]
		\label{lem:UslnSDH}
		There exists an algebra embedding
		\begin{align}
			\label{eq:semiXi}
			\widetilde{\psi}_Q:\tU_\sqq&\longrightarrow \cs\cd\widetilde{\ch}_{\Z_2}( \bfk Q)
			\\
			E_j\mapsto&\frac{-1}{q-1}[C_{S_j}], \qquad F_j\mapsto  \frac{\sqq}{q-1}[C_{S_{j}}^*]
			\\
			\tK_j\mapsto&[K_{S_j}], \qquad \qquad \tK_j'\mapsto[K^*_{S_j}],
			\qquad\forall 0\leq j\leq n-1.
		\end{align}
	\end{proposition}

	Let $\cc\widetilde{\ch}_{\Z_2}(\bfk Q)$ be the composition subalgebra of $ \cs\cd\widetilde{\ch}_{\Z_2}( \bfk Q)$, which is generated by $[C_{S_i}]$, $[C^*_{S_i}]$, $[K_{S_i}]^{\pm1}$ and $[K_{S_i}^*]^{\pm1}$ for any $i\in Q_0$. Then $\widetilde{\psi}_Q$ induces an isomorphism
	\begin{align}
		\widetilde{\psi}_Q:\tU_\sqq\stackrel{\cong}{\longrightarrow} \cc\widetilde{\ch}_{\Z_2}( \bfk Q).
	\end{align}
	Moreover,  there exists an isomorphism of algebras $\texttt{FT}_{C_n,Q}:\cc\widetilde{\ch}_{\Z_2}(\bfk Q)\stackrel{\cong}{\longrightarrow} \cc\widetilde{\ch}_{\Z_2}(\bfk C_n)$
	such that $\texttt{FT}_{C_n,Q}([C_{S_i}])=[C_{S_i}]$, $\texttt{FT}_{C_n,Q}([C^*_{S_{i}}])=[C^*_{S_{i}}]$, $\texttt{FT}_{C_n,Q}([K_{S_i}])=[K_{S_i}]$, $\texttt{FT}_{C_n,Q}([K^*_{S_i}])=[K^*_{S_i}]$ for any $0\leq i\leq n-1$.

	Let $\tU:=\tU_v(\widehat{\mathfrak{sl}}_n)$ be the universal affine quantum group of type $A_{n-1}$, and $\tUD$ be its Drinfeld presentation throughout this subsection.

	Recall that $\Phi:\tUD\rightarrow\tU$ is the Drinfeld-Beck isomorphism; see \cite{Be94}. 
	Define
	\begin{align}\label{the map Psi A}
		\Omega_{C_n}:=\widetilde{\psi}_{C_n}\circ \Phi:\tUD_{\sqq}\longrightarrow \cs\cd\widetilde{\ch}_{\Z_2}(\bfk C_n).
	\end{align}
	Then $\Omega_{C_n}$ sends:
	\begin{align}
		\label{eq:HaDrA1}K_j\mapsto [K_{S_j}],\qquad K_j' \mapsto [K_{S_j}^*],\qquad
		C\mapsto [K_\de],\qquad \widetilde{C}\mapsto [K_\de^*],
		\\
		x^+_{j,0}\mapsto \frac{-1}{q-1}[C_{S_j}], \qquad x^-_{j,0}\mapsto \frac{\sqq}{q-1}[C_{S_j}^*],\qquad \forall \,\,1\leq j\leq n-1.
	\end{align}
	For any $1\leq j\leq n-1$, $l\in\Z$ and $r\geq 1$, we define
	\begin{align}
		\label{def:haB}
		\haX^+_{j,l}:=(1-q)\Omega_{C_n}(x^+_{j,l}),&\qquad \haX^-_{j,l}:=(\sqq-\sqq^{-1})\Omega_{C_n}(x^-_{j,l}), \qquad \widehat{\bh}_{j,l}:= \Omega_{C_n}(h_{j,l}),
		\\
		\label{def:haphi}
		\widehat{\psi}_{j,r}:=& \Omega_{C_n}(\psi_{j,r}),\qquad\widehat{\varphi}_{j,-r}:= \Omega_{C_n}(\varphi_{j,-r}).
	\end{align}
	In particular, $\haX^+_{j,0}=[C_{S_j}]$, $\haX^-_{j,0}=[C^*_{S_j}]$ for any $1\leq j\leq n-1$. 
	
	In next subsection, we shall describe the root vectors $\haX^{\pm}_{j,l},\widehat{\psi}_{j,r},\widehat{\varphi}_{j,-r}$.


	%
	\subsection{Description of root vectors}
	
	\subsubsection{}

	Let $Q=(Q_0,Q_1)$ be a general quiver. For any sink $\ell \in Q_0$, define the quiver $\bs_\ell^+ Q$ by reversing all the arrows of $Q$ ending at $\ell $. 
	Associated to a sink $\ell \in Q_0$, the BGP reflection functor induces a {\rm reflection functor} (see \cite[\S3.2]{LW19b}):
	\begin{align}  \label{eq:Fl}
		F_\ell ^+: \cc_{\Z_2}(\rep_\bfk(Q))  \longrightarrow \cc_{\Z_2}(\rep_\bfk(\bs^+_\ell Q)).
	\end{align}
	The functor $F_\ell^+$ induces an isomorphism $\Gamma_\ell: \cs\cd\widetilde{\ch}_{\Z_2}(\bfk Q)  \stackrel{\sim}{\rightarrow} \cs\cd\widetilde{\ch}_{\Z_2}(\bfk \bs_\ell^+ Q)$ by \cite[Theorem 4.3]{LW22},
	and then a commutative diagram
	\begin{equation}
		\label{eq:comm}
		\xymatrix{ \tU_{ \sqq} \ar[r]^{ T_\ell} \ar[d]^{\widetilde{\psi}_{Q}} & \tU_{ \sqq} \ar[d]^{\widetilde{\psi}_{\bs_\ell Q}}
			\\
			\cs\cd\widetilde{\ch}_{\Z_2}(\bfk Q) \ar[r]^{\Gamma_\ell}  &\cs\cd\widetilde{\ch}_{\Z_2}(\bfk \bs_\ell^{+} Q)}
	\end{equation}
	where $T_\ell$ is defined in \S\ref{subsec:QG}.
	
	Dually, we can define $\bs^{-}_\ell Q$ for any source $\ell$ by reversing all the arrows of $Q$ starting from $\ell $, and a reflection functor $F_\ell^-$ associated to a source $\ell\in Q_0$ is also defined.
	This functor $F_\ell^-$ also induces an isomorphism $\Gamma_\ell^-: \cs\cd\widetilde{\ch}_{\Z_2}(\bfk Q)  \stackrel{\sim}{\rightarrow} \cs\cd\widetilde{\ch}_{\Z_2}(\bfk \bs_\ell^- Q)$.

	For any $1\leq j\leq n-1$, set $Q(j)$ to be the following quiver:
	\begin{center}\setlength{\unitlength}{0.5mm}
		\begin{equation}
			\label{quiver:affineA}
			\begin{picture}(100,30)
				\put(2,0){\circle*{2}}
				\put(51,0){\circle*{2}}
				
				\put(102,0){\circle*{2}}
				\put(52,25){\circle*{2}}
				
				\put(3.5,0){\vector(1,0){16}}
				\put(21,-2){$\cdots$}
				
				\put(32,0){\vector(1,0){16}}
				
				\put(71.5,-2){$\cdots$}
				\put(70,0){\vector(-1,0){16}}
				\put(100,0){\vector(-1,0){16}}
				
				\put(54,24.5){\vector(2,-1){47}}
				\put(50,24.5){\vector(-2,-1){47}}
				
				\put(50.5,27){\tiny $0$}
				\put(1,-6){\tiny $1$}
				\put(50,-6){\tiny $j$}
				\put(95,-6){\tiny $n-1$}
			\end{picture}
		\end{equation}
		\vspace{-0.2cm}
	\end{center}
	Using the isomorphisms $\Gamma_\ell,\Gamma_\ell^-$ of semi-derived Ringel-Hall algebras,	similarly to the proof of \cite[Proposition 5.5]{LR21}, for any $l>0$, we have
	\begin{align}
		\label{eq:hax--}
		\widehat{x}_{j,-l}^-=& (-1)^{jl}\texttt{FT}_{C_n,Q(j)}([C_{M(l\de+\alpha_j)}^*]),
		\\
		\label{eq:hax-+}
		\widehat{x}_{j,l}^-= &(-1)^{jl}\sqq[K_{l\de-\alpha_j}]^{-1}*\texttt{FT}_{C_n,Q(j)}([C_{M(l\de-\alpha_j)}]),
		\\
		\label{eq:hax+-}
		\widehat{x}_{j,-l}^+=& (-1)^{jl}\sqq[K^*_{l\de-\alpha_j}]^{-1}*\texttt{FT}_{C_n,Q(j)}([C_{M(l\de-\alpha_j)}^*]),
		\\
		\label{eq:hax++}
		\widehat{x}_{j,l}^+= &(-1)^{jl}\texttt{FT}_{C_n,Q(j)}([C_{M(l\de+\alpha_j)}]),
		\\
		\widehat{\psi}_{j,l}=&\frac{{(-1)^{jl}}}{(q-1)^2\sqq^{l-1}}\texttt{FT}_{C_n,Q(j)}\Big(\sum_{0\neq g: S_j \rightarrow M(l\de+\alpha_j)} [C_{\coker(g)}]\Big).
		\\
		\widehat{\varphi}_{j,-l}=&\frac{{(-1)^{jl}}}{(q-1)^2\sqq^{l-1}}\texttt{FT}_{C_n,Q(j)}\Big(\sum_{0\neq g: S_j \rightarrow M(l\de+\alpha_j)} [C_{\coker(g)}^*]\Big).
	\end{align}

	It is interesting to give explicit formulas for all root vectors in $\cs\cd\widetilde{\ch}_{\Z_2}(\bfk C_n)$, which is difficult in general. We describe some of them below.

	For $1\leq j\leq {n}$ and any $\alpha\in K_0(\rep^{\rm nil}_\bfk (C_n))$, set
	\begin{align}
		\label{def:Mjalpha}
		\cm_{j,\alpha}:=\{[M]\mid \widehat{M}=\alpha,{\rm soc}(M)\subseteq S_1\oplus  \cdots \oplus S_j\}.
	\end{align}
	Define
	\begin{align}
		\pi_{j,1}^+=\frac{-\sqq^{-j}}{\sqq-\sqq^{-1}}\sum\limits_{\cm_{j,\delta}}(-1)^{\dim\End(M)}[C_M],
		\\
		\pi_{j,1}^-=\frac{-\sqq^{-j}}{\sqq-\sqq^{-1}}\sum\limits_{\cm_{j,\delta}}(-1)^{\dim\End(M)}[C^*_M],
	\end{align}
	Set $\pi_{0,1}=0$.
	
	\begin{proposition}[Hubery, see also \text{\cite[\S6.3]{DJX12}}]
		\label{prop:DrGenA}
		For any $1\leq j\leq n-1$, we have
		\begin{align}
			&\widehat{x}^+_{j,-1}=\sqq^{-j+2}\sum\limits_{\cm_{j+1,\delta-\alpha_{j}}}(-1)^{\dim\End(M)}[K^*_{\de-\alpha_j}]^{-1}*[C_M^*],
			\\
			&\widehat{x}^-_{j,1}=\sqq^{-j+2}\sum\limits_{\cm_{j+1,\delta-\alpha_{j}}}(-1)^{\dim\End(M)}[K_{\de-\alpha_j}]^{-1}*[C_M],
			\\
			&\widehat{h}_{j,\pm1}=\pm\big(\pi^\pm_{j+1,1}-(\sqq+\sqq^{-1})\pi^\pm_{j,1}+\pi^\pm_{j-1,1}\big).
		\end{align}
	\end{proposition}
	
	Similarly to the proof of \cite[Proposition 5.1, Corollary 5.2]{LR21}, we can give another proof of Proposition \ref{prop:DrGenA}.

	\subsubsection{Root vectors at vertex $1$}
	
	In the following, we  shall describe in $\cs\cd\widetilde{\ch}_{\Z_2}(\bfk C_n)$ the root vectors at the vertex $1$. In fact,  $\widehat{x}^+_{1,r}$ (and $\widehat{x}^-_{1,-r}$) are obtained in \cite[\S6.5]{Sch04} for $r>0$.

	For any nilpotent $\bfk C_n$-module $M$, let $\ell(M)$ be the number of its indecomposable direct summands.	For any $r>0$, denote by
	\begin{align*}
		&\cm_{r\delta+\alpha_{1}}=\{[S_{1,1}^{(bn+1)}\oplus S_{1,0}^{(\nu)}]\mid b\geq 0, \;|\nu|+b=r\};\\
		&\cm_{r\delta-\alpha_{1}}=\{[S_{1,0}^{(an-1)}\oplus S_{1,0}^{(\nu)}]\mid a\geq 1, \;|\nu|+a=r\};
	\end{align*}
	and set
	\begin{align}
		\label{eq:Mrde+alpha11}
		[M_{r\delta\pm\alpha_{1}}]:=\sum\limits_{[M]\in\cm_{r\delta\pm\alpha_{1}}}\bn(\ell(M)-1)[\![M]\!],
	\end{align}
	where
	\[\bn(l)=\prod_{i=1}^l(1-\sqq^{2i}),\quad \forall l\geq1, \text{ and }\bn(0)=1.\]

	\begin{proposition}
		\label{prop:realroot1}
		For any $r>0$, we have
		\begin{align}
			\label{realroot2}
			\widehat{x}^+_{1,r}= (q-1)[C_{M_{r\de+\alpha_{1}}}], \qquad \widehat{x}^+_{1,-r}= (1-q)\sqq[K^*_{r\de-\alpha_{1}}]^{-1}*[C^*_{M_{r\de-\alpha_{1}}}],
			\\
			\widehat{x}^-_{1,r}=(1-q)\sqq [K_{r\de-\alpha_{1}}]^{-1}*[C_{M_{r\de-\alpha_{1}}}], \qquad \widehat{x}^-_{1,-r}=(q-1)[C^*_{M_{r\de+\alpha_{1}}}].
		\end{align}
	\end{proposition}

	\begin{proof}
		We only need to prove \eqref{realroot2}, since the other one is similar.
		
		From  \eqref{eq:hax+-}--\eqref{eq:hax++}, it is enough to compute $\texttt{FT}_{C_n,Q(j)}([M(l\de-\alpha_j)])$ and $\texttt{FT}_{C_n,Q(j)}([M(l\de+\alpha_j)])$ in $\widetilde{\ch}(\bfk C_n)$ for any $l>0$; see \eqref{triandecomp}. However, from \cite[Proposition 9.1]{LR21}, we know
		\begin{align*}&\texttt{FT}_{C_n,Q(j)}([M(l\de-\alpha_j)])=(-1)^{jl}(1-q)[M_{l\de-\alpha_1}],\\
			&\texttt{FT}_{C_n,Q(j)}([M(l\de+\alpha_j)])=(-1)^{jl}(q-1)[M_{l\de+\alpha_1}]
		\end{align*}
		in $\iH(\bfk C_n,\Id)$.
		Since $\iH(\bfk C_n,\Id)$ is a filtered algebra with its associated graded algebra $\iH(\bfk C_n,\Id)^{\gr}$ isomorphic to $\widetilde{\ch}(\bfk C_n)\otimes \ct(\bfk C_n,\Id)$, we have
		\begin{align*}
			&\texttt{FT}_{C_n,Q(j)}([M(l\de-\alpha_j)])=(-1)^{jl}(1-q)[M_{l\de-\alpha_1}],
			\\
			&\texttt{FT}_{C_n,Q(j)}([M(l\de+\alpha_j)])=(-1)^{jl}(q-1)[M_{l\de+\alpha_1}]
		\end{align*}
		in $\widetilde{\ch}(\bfk C_n)$.
	\end{proof}

	For any $r>0$, set
	\begin{align}
		\label{def:Mr}
		\cm_{r\de}:=&\big\{[S_{1,1}^{(\nua  p_1)}\oplus S_{1,0}^{(\nu)}]\mid a\geq 1, a+|\nu|=r\big\}
		\\\notag
		& \bigcup\big\{[S_{1,0}^{(\nua  p_1-1)}\oplus S_{1,1}^{(b  p_1+1)}\oplus S_{1,0}^{(\nu)}]\mid  a+b+|\nu|=r\big\}.
	\end{align}
	
	For a partition $\lambda$, let $\ell(\lambda)$ be the number of nonzero parts in $\lambda$.	
	
	\begin{proposition}
		\label{prop:imageroot}
		For any $r>0$, we have
		\begin{align}
			\label{eq:TH11r}
			\widehat{\psi}_{1,r}=
			&\frac{\sqq}{q-1} \sum\limits_{|\lambda|=r}\bn(\ell(\lambda))  [\![C_{S_{0}^{(\lambda)}}]\!]
			+\sqq^{-1} \sum\limits_{[M]\in\cm_{r\de}}\bn(\ell(M)-1) [\![C_M]\!],
			\\
			\widehat{\varphi}_{1,-r}=
			&\frac{\sqq}{q-1} \sum\limits_{|\lambda|=r}\bn(\ell(\lambda))  [\![C^*_{S_{0}^{(\lambda)}}]\!]
			+\sqq^{-1} \sum\limits_{[M]\in\cm_{r\de}}\bn(\ell(M)-1) [\![C^*_M]\!].
		\end{align}
	\end{proposition}
	
	\begin{proof}
		The proof follows by using the same argument of Proposition \ref{prop:realroot1} with the help of \cite[Proposition 9.2]{LR21}.
	\end{proof}

	\section{$\imath$Hall algebras and $\imath$quantum loop algebras}
	\label{sec:hom}
	
	In this section, we shall formulate the main result of this paper.
	
	\subsection{Root vectors in tubes}
	\label{subsec:embeddingtube}
	Let $\varrho$ be an involution of the star-shaped graph $\Gamma$, which can be lifted to an involution $\varrho$ of the weighted projective line $\X$ of weight type $(\bp,\ul{\bla})$.  Consider the Kac-Moody algebra $\fg=\fg(\Gamma)$ associated to the graph $\Gamma$. Let $C$ be the generalised Cartan matrix of $\Gamma$. Then $\varrho$ can be viewed as an involution of $C$, 
	and we obtain the $\imath$quantum group $\tUiD=\tUiD(\Lg)$.
	
	Recall that $\scrt_{\bla_i}$ is the Serre subcategory of $\coh(\X)$ consisting of torsion sheaves supported at $\bla_i$, and we have an equivalence $\scrt_{\bla_i}\cong\rep^{\rm nil}_\bfk(C_{p_i})$.
	
	\subsubsection{Tubes invariant under $\varrho$}
	We assume $\scrt_{\bla_i}$ is invariant under $\varrho$ in the following. Then the equivalence $\scrt_{\bla_i}\cong\rep^{\rm nil}_\bfk(C_{p_i})$ induces a canonical embedding of Ringel-Hall algebras
	$\ch(\cc_{\Z_1}(\rep_\bfk^{\rm nil} (C_{p_i})))\rightarrow \ch(\cc_\varrho(\coh(\X_\bfk)))$, and then an embedding of $\imath$Hall algebras:
	\begin{equation}
		\label{eq:embeddingx}
		\iota_i: \iH(\bfk C_{p_i},\Id)\longrightarrow\iH(\X_\bfk,\varrho).
	\end{equation}

	Inspired by \eqref{def:haBa}, we define
	\begin{align}
		\label{def:haBThH}
		\haB_{[i,j],l}:= \iota_i (\haB_{j,l}),\quad \widehat{\Theta}_{[i,j],r}:=\iota_i(\widehat{\Theta}_{j,r}),\quad \widehat{H}_{[i,j],r}:=\iota_{i}(\widehat{H}_{j,r}),
	\end{align}
	for any $1\leq j\leq p_i-1$, $l\in\Z$ and $r>0$.
	By definition and \eqref{exp}, we have
	\begin{align}
		\label{exp haH}
		1+ \sum_{m\geq 1} (\sqq-\sqq^{-1})\haTh_{[i,j],m} u^m  = \exp\Big( (\sqq-\sqq^{-1}) \sum_{m\geq 1} \widehat{H}_{[i,j],m} u^m \Big).
	\end{align}
	
	For convenience, we set $\haTh_{[i,j],0}= \frac{1}{\sqq-\sqq^{-1}}$ for any $1\leq i\leq \bt,1\le j\le p_i-1$.

	\subsubsection{Tubes not invariant under $\varrho$}
	We assume $\scrt_{\bla_i}$  to be not invariant  under $\varrho$ in the following. 
	Then the restriction of $\varrho$ to  $\ct_{\bla_i,\bla_{\varrho i}}:=\scrt_{\bla_i}\times \ct_{\bla_{\varrho i}}$ is just the swap, which gives an embedding $\cc_{\swa}(\ct_{\bla_i,\bla_{\varrho i}})\rightarrow \cc_\varrho(\coh(\X))$. 
	We identify  $\cc_{\swa}(\ct_{\bla_i,\bla_{\varrho i}})$ with $\cc_{\Z_2}(\rep^{\rm nil}_\bfk(C_{p_i}))$ by using the equivalence
	\begin{align}
		F_i:\rep^{\rm nil}_\bfk(C_{p_i})\stackrel{\cong}{\longrightarrow} \scrt_{\bla_i}.
	\end{align}
	We obtain an embedding of Ringel-Hall algebras
	$$\widetilde{\ch} \big(\cc_{\Z_2}(\rep^{\rm nil}_\bfk(C_{p_i}))\big)\longrightarrow \widetilde{\ch}(\cc_\varrho(\coh(\X_\bfk))),$$
	and then a morphism of semi-derived Ringel-Hall algebras:
	\begin{align}
		\iota_i:\cs\cd\widetilde{\ch}_{\Z_2}(\bfk C_{p_i})\longrightarrow \iH(\X_\bfk,\varrho).
	\end{align}
	In particular, for any $X\in\rep^{\rm nil}_\bfk(C_{p_i})$, we have
	\begin{align*}
		\iota_i([C_X])=[F_{\varrho i}(X)],\qquad \iota_i([C_X^*])=[F_i(X)],\\
		\iota_i([K_X])=[K_{F_i(X)}],\qquad\iota_i([K^*_X])=[K_{F_{\varrho i}(X)}].
	\end{align*}
	However, $\iota_i$ is not injective by noting that $$\iota_{i}([K_\de])=[K_\de]=\iota_i([K_\de^*]).$$

	Inspired by \eqref{def:haB}--\eqref{def:haphi}, we define, for any $1\leq j\leq p_i-1$, $l\in\Z$ and $r\geq 1$,
	\begin{align}
		\label{def:haBXno}
		\haB_{[i,j],l}:=&\iota_i(\haX^-_{j,-l}),\qquad \haB_{[\varrho i,j],l}:=\iota_i(\haX^+_{j,l}), \\
		\label{def:haHXno}
		\qquad \widehat{H}_{[i,j],r}:= &\iota_i(\widehat{\bh}_{j,-r}),  \qquad\widehat{H}_{[\varrho i,j],r}:=\iota_i(\widehat{\bh}_{j,r}),
		\\
		\label{def:haphino}
		\widehat{\Theta}_{[i,j],r}:=&\iota_i(\widehat{\varphi}_{j,-r}),\quad\quad\widehat{\Theta}_{[\varrho i,j],r}:= \iota_i(\widehat{\psi}_{j,r}).
	\end{align}
	In particular, $\haB_{[i,j],0}=[S_{ij}]$, $\haB_{[\varrho i,j],0}=[S_{\varrho i,j}]$ for any $1\leq j\leq p_i-1$.

	\subsection{$\imath$Hall algebra of the projective line and $\tUi(\widehat{\mathfrak{sl}}_2)$}
	
	In this subsection, we assume $\tUi=\tUi(\widehat{\mathfrak{sl}}_2)$ is the $q$-Onsager algebra. Let $\tUiD=\tUiD(\widehat{\mathfrak{sl}}_2)$ be the $\Q(v)$-algebra  generated by $\K_1^{\pm1}$, $C^{\pm1}$, $H_{m}$ and $B_{1,r}$, where $m\geq1$, $r\in\Z$, subject to the following relations, for $r,s\in \Z$ and $m,n\ge 1$:
	\begin{align}
		\K_1\K_1^{-1}=1, C C^{-1}& =1, \quad \K_1, C \text{ are central, }
		\\
		[H_m,H_n] &=0,  \label{iDR1On}
		\\
		[H_m, B_{1,r}] &=\frac{[2m]}{m} B_{1,r+m}-\frac{[2m]}{m} B_{1,r-m}C^m,
		\label{iDR2On}
		\\
		\label{iDR3On}
		[B_{1,r}, B_{1,s+1}]_{v^{-2}}  -v^{-2} [B_{1,r+1}, B_{1,s}]_{v^{2}}
		&= v^{-2}\Theta_{s-r+1} C^r \K_1-v^{-4} \Theta_{s-r-1} C^{r+1} \K_1 \\\notag
		&\quad +v^{-2}\Theta_{r-s+1} C^s \K_1-v^{-4} \Theta_{r-s-1} C^{s+1} \K_1.\notag
	\end{align}
	Here
	\begin{align}
		\label{eq:exp1}
		1+ \sum_{m\geq 1} (v-v^{-1})\Theta_{m} z^m  = \exp\Big( (v-v^{-1}) \sum_{m\ge 1}  H_m z^m \Big).
	\end{align}

	By Proposition \ref{prop:autogroupP1}, we can view
	$\varrho$ to be an involution of $\coh(\P^1_\bfk)$ fixing $\co$, and $\cc_\varrho(\coh(\P^1_\bfk))$ is the $\imath$-category. Let $\iH(\P^1_\bfk,\varrho)$ be the $\imath$Hall algebra of $\cc_\varrho(\coh(\P^1_\bfk))$. Inspired by \cite[(4.2), Proposition 6.3]{LRW20a}, we define  in $\iH(\P^1_\bfk,\varrho)$
	\begin{align}
		\label{eq:HaThetam}
		\haT_{m}= \frac{1}{(q-1)^2\sqq^{m-1}}\sum_{0\neq f:\co\rightarrow \co(m) } [\coker f],\qquad \forall m\geq1.
	\end{align}
	We denote by $\mathbb N (\PL)$ the set of all functions ${\bf n}: \PL \rightarrow \N$ such that ${\bf n}_x \neq 0$ for only finitely many $x\in \PL$. We sometimes write $\bn \in \mathbb N (\PL)$ as ${\bf n} =({\bf n}_x)_{x\in \PL}$ or  ${\bf n} =({\bf n}_x)_x$.  We define a partial order $\leq$ on $\mathbb N (\PL)$:
	\begin{align}  \label{eq:mn}
		\bn \leq \bm \text{ if and only if } \bn_x \le \bm_x \text{ for all }x\in \PL.
	\end{align}
	For $\bn \in \mathbb N (\PL)$, we denote the torsion sheaf
	\begin{align}  \label{eq:Sbn}
		S_\bn =\bigoplus\limits_{x\in\PL}S_{x}^{(\bn_x)},
	\end{align}
	whose degree is given by
	\begin{align*}
		|| {\bf n} || :=\sum_{x\in \P^1_{\mathbf{k}}} d_x {\bf n}_x.
	\end{align*}
	By \cite[(4.1)]{LRW20a}, we have
	\begin{align}
		\label{def:Theta}
		\haT_{m}= \frac{1}{(q-1)\sqq^{m-1}} \sum\limits_{||\bn||=m}[S_\bn],
		\qquad \text{ for } m\ge 1.
	\end{align}

	For convenience, set $\haT_0:=\frac{1}{\sqq-\sqq^{-1}}$.
	
	\begin{proposition}[cf. \text{\cite[Theorem 4.2]{LRW20a}}]
		\label{prop:P1-Onsager}
		There exists a $\Q(\sqq)$-algebra homomorphism
		\begin{align}
			\label{eq:phi}
			\Omega: \tUiD_{|_{\sqq}} \longrightarrow \iH(\P^1_\bfk,\varrho)
		\end{align}
		which sends, for all $r\in \Z$ and $m \ge 1$,
		\begin{align*}
			\K_1\mapsto [K_\co],  \quad
			C\mapsto [K_\de],  \quad
			B_{1,r} \mapsto -\frac{1}{q-1}[\co(r)],  \quad
			\Theta_{m}\mapsto  \haT_m.
		\end{align*}
	\end{proposition}

	\begin{proof}
		The involution $\varrho$ fixes $\co(r)$ for any $r\in\Z$, so one can prove $\Omega$ preserves \eqref{iDR3On} by using the same proof of \cite[\S4.2]{LRW20a}. Let us check \eqref{iDR1On} first in the following.
		
		Recall that the subcategory of torsion sheaves $\scrt=\coprod_{x\in\P_\bfk^1}\scrt_x$, where $\scrt_x\simeq \rep^{\rm nil}_{\bfk_x}(C_1)$; see Lemma \ref{lem:isoclasses Tor}.
		Note that $\coker( f)\in \scrt$ for any $0\neq f:\co\rightarrow \co(m)$. It is enough to prove that $[M]*[N]=[N]*[M]$ for any $M,N\in\scrt$.

		For any $M\in\scrt_x,N\in\scrt_y$, we have $$\Hom_{\cc_\varrho(\coh(\P^1_\bfk))}(M,N)=0\text{ if }x\neq y;\qquad
		\Ext^1_{\cc_\varrho(\coh(\P^1_\bfk))}(M,N)=0\text{ if }\varrho (y)\neq x\neq y.$$ 
		So we only need to consider the following two cases.
		
		Case 1: $\underline{M,N\in\scrt_x}$. If $\varrho(x)\neq x$, then $\Ext^1_{\cc_\varrho(\coh(\P^1_\bfk))}(M,N)=\Ext^1_{\P^1_\bfk}(M,N)$. Then \eqref{iDR1On} preserved by $\Omega$ follows from the classic result of Steinitz and Hall. If $\varrho(x)=x$, then it follows from \cite[Lemma 4.5]{LRW20a}.
		
		Case 2: $\underline{M\in\scrt_x, N\in\ct_{\varrho x}\text{ with }\varrho x\neq x}$. We have $\Ext^1_{\cc_\varrho(\coh(\P^1_\bfk))}(M,N)=\Hom(M,\varrho(N))$. So any short exact sequence in $\Ext^1_{\cc_\varrho(\coh(\P^1_\bfk))}(M,N)$ is of the form
		$0\rightarrow N\rightarrow C_f\rightarrow M\rightarrow0$ for $f:M\rightarrow\varrho(N)$. Note that the Euler form is trivial. So
		$$[M]*[N]=\sum_{f:M\rightarrow \varrho(N)}[C_f]=\sum_{f:M\rightarrow \varrho(N)}[\Ker (f)\oplus \coker (\varrho(f))]*[K_{\Im(f)}].$$
		Similarly, 
		$$[N]*[M]=\sum_{g:N\rightarrow \varrho(M)}[\Ker(g)\oplus \coker (\varrho(g))]*[K_{\Im(g)}].$$
		Identify $\scrt_x$ with $\rep_{\bfk_x}^{\rm nil}(C_1)$ in the following. 
		Let $D=\Hom_\bfk(-,\bfk)$ be the duality functor. 
		Then $D(M)\cong M$ and $D(N)\cong N$. For any $g: N\rightarrow \varrho(M)$, we can assume $g=D(\varrho(f))$ for some (unique) $f: M\rightarrow\varrho(N)$, which gives a bijection 
		$\Hom(N,\varrho(M))\leftrightarrow \Hom(M,\varrho(N))$. Note that
		$\Ker(g)\cong \coker \varrho (f)$, $\coker\varrho(g)\cong \Ker (f)$, and $\Im(g)\cong \Im\varrho(f)$.  Furthermore, $\widehat{M}=\widehat{\varrho(M)}=(d_x\dim_{\bfk_x}M)\de$ in $K_0(\coh(\P_\bfk^1))$. So $[K_{\Im(f)}]=[K_{\Im (\varrho(f))}]$. Then we have
		\begin{align*}
			[N]*[M]=&\sum_{g:N\rightarrow \varrho(M)}[\Ker(g)\oplus \coker (\varrho(g))]*[K_{\Im(g)}]
			\\
			=&\sum_{f:M\rightarrow \varrho(N)}[\Ker (f)\oplus \coker (\varrho(f))]*[K_{\Im(f)}]
			\\
			=&[M]*[N].
		\end{align*}
		
		For \eqref{iDR2On}, by using \eqref{THBij}, we need to prove the following identity: 
		\begin{align}  \label{TOTO}
			&\big[\haT_{m},[\co(r)] \big] + \big[\haT_{m-2},[\co(r)] \big]*[K_\de]
			\\\notag
			&=\sqq^{2} \big[\haT_{m-1},[\co(r+1)] \big]_{\sqq^{-4}} +\sqq^{-2} \big[\haT_{m-1},[\co(r-1)] \big]_{\sqq^{4}}*[K_\de],
		\end{align}
		for any $m>0,r\in\Z$.
		Since it is almost the same as \cite[Proposition 4.7]{LRW20a}, we only sketch its proof here. 
		It is enough to compute $\haT_m*[\co(r)]$ and $[\co(r)]*\haT_m$.
		
		For $\haT_m*[\co(r)]$, its computation is the same as \cite[Lemma 4.10]{LRW20a}, hence omitted here.
		
		Let us compute $[\co(r)]*\haT_m$. For any morphism $f:\co(r)\rightarrow  \varrho(S_{\bn})$, we claim $[C_f]=[\co(r-a)\oplus  S_{\bm-\bl}]*[K_{a\de}]$ for some $\bl\leq \bm$ with $||\bl||=a$. Then we can get $[\co(r)]*\haT_m$ by using the same proof of
		\cite[Proposition 4.13]{LRW20a}. Then \eqref{TOTO} follows by using the proof of \cite[Proposition 4.7]{LRW20a}.

		Let us prove the claim. 
		For a closed point $x$ in $\PL$, the torsion sheaf $S_{x}^{(n)}$ is uniserial, and hence its subsheaf and quotient sheaf are of the form $S_{x}^{(a)}$, for $0\leq a\leq n$. Recall there is no nonzero homomorphism between torsion sheaves supported on distinct closed points. Therefore,  as a subsheaf of $ S_{\bn}$, $\Im (\varrho(f))=S_{\bl}$, for some $\bl\leq \bn$. It follows that $\coker (\varrho(f))\cong S_{\bn-\bl}$.
		
		Note that $\widehat{\Im\varrho(f)}=a\de=\widehat{\Im(f)}$ by noting that $K_0(\coh(\P^1))$ is generated by $\widehat{\co}$ and $\de$, which are fixed by $\varrho$. Then 
		$[K_{\Im(f)}]=[K_{a\de}]$. 
		Observe that any subsheaf of a line bundle is again a line bundle, which is determined by its degree. Hence $\Ker (f)\cong \co(r-a)$.
		Then the claim follows by
		\begin{align*}
			[C_f]=&[\Ker (f)\oplus \coker\varrho(f)]*[K_{\Im(f)}]
			\\
			=&[\co(r-a)\oplus  S_{\bm-\bl}]*[K_{a\de}].
		\end{align*}
	\end{proof}

	Define $\haH_m:=\Omega(H_m)$ for $m\geq1$.	Then we have
	\begin{align}
		\label{eq:Tr}
		1+ \sum_{m\geq 1} (\sqq-\sqq^{-1})\haT_{m} z^m  = \exp\big( (\sqq-\sqq^{-1}) \sum_{m\ge 1} \haH_m z^m \big).
	\end{align}
	
	We shall describe the root vector $\haH_m$.

	Recall that any indecomposable object in $\tor_x(\PL)$, for $x\in\PL$, has the form $S_x^{(n)}$ of length $n\geq 1$. For any partition $\lambda=(\lambda_1,\dots,\lambda_r)$, define
	\[S_x^{(\lambda)}:=S_x^{(\lambda_1)}\oplus\cdots \oplus S_x^{(\lambda_r)}. \]
	
	For any $x\in\PL$ and $m\geq 1$, we define (compare \cite{LRW21})
	\begin{align}  \label{eq:PTh}
		\haP_{m,x}:=\sum\limits_{\lambda\vdash m} n_x(\ell(\lambda)-1)\frac{[S_x^{(\lambda)}]}{|\Aut(S_x^{(\lambda)})|},
		\qquad
		\haT_{m,x}:=\frac{[S_x^{(m)}]}{\sqq_x-\sqq_x^{-1}},
	\end{align}
	where
	\[
	n_x(l)=\prod_{i=1}^l(1-\sqq_x^{2i})=\prod_{i=1}^l(1-q_x^{i}).
	\]
	Introduce the generating functions
	\begin{align}
		\haP_x(z):=&\sum_{m\geq1} \haP_{m,x} z^{m-1},
		\label{eq:Px} \\
		\haT_x(z):=&1+\sum_{m\geq1} (\sqq_x-\sqq_x^{-1})\haT_{m,x} z^{m}
		=1+\sum_{m\geq1} [S_x^{(m)}]z^{m}.
		\label{eq:Thx}
	\end{align}
	
	\begin{lemma}
		We have
		\begin{align}\label{exp pression between theta and Hm}
			\haT_x(z)  =\exp\big( (\sqq_x-\sqq_x^{-1}) \sum_{m\ge 1} \haH_{m,x} z^{m} \big),
		\end{align}
		where
		\begin{align}  \label{eq:Tx}
			\haH_{m,x} =\sqq_x^{m} \frac{[m]_{\sqq_x}}{m} \haP_{m,x} -\delta_{m, ev} \delta_{x,\varrho x} \sqq_x^{\frac{m}{2}} \frac{[m/2]_{\sqq_x}}{m} [K_{\frac{m}{2}d_x\de}].
		\end{align}
	\end{lemma}

	\begin{proof}
		Identify $\tor_x(\P^1_\bfk)$ with $\rep^{\rm nil}_{\bfk_x}(C_1)$ in the following. 
		
		If $\varrho(x)=x$, it is equivalent to consider $\iH(\bfk_x C_{1},\Id)$, and then the desired formula follows from \cite[Proposition 6.7]{LRW21}.

		If $\varrho(x)\neq x$, it is equivalent to consider $\cs\cd\widetilde{\ch}_{\Z_2}(\bfk_x C_1)$. Let $\widetilde{\ch}(\bfk_x C_1)$ be the twisted Hall algebra of $\rep^{\rm nil}_{\bfk_x}(C_1)$. Using the algebra embedding $R^+:\widetilde{\ch}(\bfk_x C_1)\rightarrow \cs\cd\widetilde{\ch}_{\Z_2}(\bfk_x C_1)$, the formula follows from \cite[Example 4.12]{Sch12}.
	\end{proof}

	Recall $\haT_m$ from \eqref{eq:HaThetam}.
	Define the generating function
	\begin{align}  \label{eq:haT}
		\haT(z) =1 + \sum_{m\ge 1} (\sqq -\sqq^{-1}) \haT_m z^m.
	\end{align}
	Since the categories $\tor_x(\PL)$ for $x \in \PL$ are orthogonal, by \eqref{def:Theta} and \eqref{eq:Thx} we have
	\begin{align*}
		\haT(\sqq z) =1+\sum_{m\geq 1}\sum_{||\bn||=m}[S_{\bn}]z^m=\prod_{x\in \PL} \haT_x(z^{d_x}).\\
	\end{align*}

	\begin{proposition}
		\label{prop:HaH}
		For $m\geq1$, we have
		\begin{align}\label{formula for Hm}
			\haH_m
			=&\sum_{x,d_x|m} \frac{[m]}{m} d_x \sum_{|\lambda|=\frac{m}{{d_x}}} n_x(\ell(\lambda)-1)\frac{[S_x^{(\lambda)}]}{\big|\aut(S_x^{(\lambda)})\big|}  -\delta_{m, ev} \frac{1-\sqq^{-m}}{m} \sum_{x, \varrho(x)=x,d_x|\frac{m}{2}}d_x [K_{\frac{m}{2}\de}].
		\end{align}
	\end{proposition}
	
	\begin{proof}
		The proof is similar to \cite[Proposition 6.3]{LRW21}.
		By \eqref{exp pression between theta and Hm},
		we have
		\begin{align*}
			\haT(\sqq z)=&\prod_{x\in \PL} \haT_x(z^{d_x})\\
			=&\prod_{x\in \PL} \exp\big( (\sqq_x-\sqq_x^{-1}) \sum_{m\ge 1} \haH_{m,x} z^{md_x} \big)\\
			=&\exp\big(\sum_{x\in \PL} (\sqq_x-\sqq_x^{-1}) \sum_{m\ge 1} \haH_{m,x} z^{md_x} \big)\\
			=&\exp\big(\sum_{m\geq 1}  \sum_{x, d_x|m} (\sqq_x-\sqq_x^{-1})\haH_{\frac{m}{d_x},x} z^{m} \big)\\
			=&\exp\big(  \sum_{m\ge 1}\sum_{x, d_x|m}(\sqq_x-\sqq_x^{-1}) \big(\sqq^{m} \frac{[m/d_x]_{\sqq_x}}{m/d_x} \haP_{\frac{m}{d_x},x}-\delta_{\frac{m}{d_x}, ev} \de_{x,\varrho(x)}\sqq^{\frac{m}{2}} \frac{[m/(2d_x)]_{\sqq_x}}{m/d_x} [K_{\frac{m}{2}\de}]\big) z^m\big).
		\end{align*}
		Observe from \eqref{eq:Tr} that
		\begin{align*}
			\exp\big( (\sqq-\sqq^{-1}) \sum_{m\ge 1} \haH_m z^m \big)
			=&1+ \sum_{m\geq 1} (\sqq-\sqq^{-1})\haT_{m} z^m = \haT(z).
		\end{align*}
		So it suffices to show that
		\begin{align}
			\label{eq:derived}
			&{\sqq^{-m}}\sum_{x, d_x|m} (\sqq_x-\sqq_x^{-1})\big(\sqq^{m} \frac{[m/d_x]_{\sqq_x}}{m/d_x} \haP_{\frac{m}{d_x},x} -\delta_{\frac{m}{d_x}, ev} \de_{x,\varrho(x)}\sqq^{\frac{m}{2}} \frac{[m/(2d_x)]_{\sqq_x}}{m/d_x} [K_{\frac{m}{2}\de}]\big)\\\notag
			&\qquad\qquad={(\sqq-\sqq^{-1})}\Big(\sum_{x,d_x|m} \frac{[m]}{m} d_x \sum_{|\lambda|=\frac{m}{{d_x}}} n_x(\ell(\lambda)-1)\frac{[S_x^{(\lambda)}]}{\big|\aut(S_x^{(\lambda)})\big|}  \\
			&\notag \qquad\qquad\qquad\qquad\qquad\qquad-\delta_{m, ev} \frac{1-\sqq^{-m}}{m} \sum_{x, \varrho(x)=x,d_x|\frac{m}{2}}d_x [K_{\frac{m}{2}\de}]\Big).
		\end{align}
		
		We have
		\begin{align*}
			&\sqq^{-m}\sum_{x, d_x|m}(\sqq_x-\sqq_x^{-1}) \sqq^{m} \frac{[m/d_x]_{\sqq_x}}{m/d_x} \haP_{\frac{m}{d_x},x} \\
			&=\sum_{x, d_x|m}  \frac{\sqq^m-\sqq^{-m}}{m} d_x \sum\limits_{\lambda\vdash \frac{m}{d_x}} n_x(\ell(\lambda)-1)\frac{[S_x^{(\lambda)}]}{|\Aut(S_x^{(\lambda)})|} \\
			&=(\sqq-\sqq^{-1})\sum_{x,d_x|m} \frac{[m]}{m} d_x \sum_{|\lambda|=\frac{m}{{d_x}}} n_x(\ell(\lambda)-1)\frac{[S_x^{(\lambda)}]}{\big|\aut(S_x^{(\lambda)})\big|}.
		\end{align*}
		Moreover, we compute
		\begin{align*}
			&\sqq^{-m}\sum_{x, d_x|m}  (\sqq_x-\sqq_x^{-1})\delta_{\frac{m}{d_x}, ev}\de_{x,\varrho(x)}\sqq^{\frac{m}{2}} \frac{[m/(2d_x)]_{\sqq_x}}{m/d_x} [K_{\frac{m}{2}\de}]\\
			&=\sum_{x, d_x|m} \delta_{\frac{m}{d_x}, ev} \de_{x,\varrho(x)}\sqq^{-\frac{m}{2}} \frac{\sqq^{\frac{m}{2}}-\sqq^{-\frac{m}{2}}}{m} d_x [K_{\frac{m}{2}\de}]
			= \delta_{m, ev} \frac{1-\sqq^{-m}}{m} \sum_{x, \varrho(x)=x,d_x|\frac{m}{2}}d_x [K_{\frac{m}{2}\de}].
		\end{align*}
		So \eqref{eq:derived} holds. We are done.
	\end{proof}


	%
	%
	\subsection{Root vectors at $\star$ point} 
	\label{Embedding from projective line to weighted projective line}
	
	Let $\mathcal{C}$ be the Serre subcategory of $\coh(\X)$ generated by those simple sheaves $S$ satisfying $\Hom(\co, S)=0$. Recall from \cite{BS13} that the Serre quotient $\coh(\X)/\mathcal{C}$ is equivalent to the category $\coh(\P^1)$ and the canonical functor $\coh(\X)\rightarrow\coh(\X)/\mathcal{C}$ has an exact fully faithful right adjoint functor
	\begin{equation}
		\label{the embedding functor F}
		\mathbb{F}_{\X,\P^1}: \coh(\P^1)\rightarrow\coh(\X),
	\end{equation}
	which sends $$\co_{\P^1}(l)\mapsto\co(l\vec{c}),\quad S_{\bla_i}^{(r)}\mapsto S_{i,0}^{(rp_i)},\quad S_{x}^{(r)}\mapsto S_{x}^{(r)}$$ 
	for any $l\in\Z, r\geq 1$, $1\leq i\leq \bt$ and $x\in\PL\setminus\{\bla_1,\cdots, \bla_\bt\}$. Then $\mathbb{F}_{\X,\P^1}$ induces an exact fully faithful functor $\cc_\varrho(\coh(\P^1))\rightarrow \cc_\varrho(\coh(\X))$, which is also denoted by $\mathbb{F}_{\X,\P^1}$. This functor $\mathbb{F}_{\X,\P^1}$ induces a canonical embedding of Hall algebras $\ch\big(\cc_\varrho(\coh(\P^1))\big)\rightarrow \ch\big(\cc_\varrho(\coh(\X))\big)$, and then an embedding
	\begin{equation}
		\label{the embedding functor F on algebra}
		F_{\X,\P^1}:\iH(\P^1_\bfk,\varrho) \longrightarrow \iH(\X_\bfk,\varrho).
	\end{equation}
	
	Hence, we define
	\begin{align}
		\label{def:Theta star}
		\widehat{\Theta}_{\star,m}:=F_{\X,\P^1}(\haT_m)=\frac{1}{(q-1)^2\sqq^{m-1}}\sum_{0\neq f:\co\rightarrow \co(m\vec{c}) } [\coker f],
	\end{align}
	and
	\begin{align}
		\label{formula for Hm}
		\haH_{\star,m}:=&F_{\X,\P^1}(\haH_{m})
		\\\notag
		=&\sum_{x,d_x|m} \frac{[m]_\sqq}{m} d_x \sum_{|\lambda|=\frac{m}{{d_x}}} \bn_x(\ell(\lambda)-1)\frac{F_{\X,\P^1}([S_x^{(\lambda)}])}{\big|\Aut(S_x^{(\lambda)})\big|}  -\delta_{m, ev} \frac{1-\sqq^{-m}}{m} \sum_{x, \varrho(x)=x,d_x|\frac{m}{2}}d_x [K_{\frac{m}{2}\de}].
	\end{align}
	Here, for any partition $\lambda=(\lambda_1, \lambda_2,\cdots)$, $F_{\X,\P^1}([S_x^{(\lambda)}])=[S_x^{(\lambda)}]$ for any $x\in\PL\setminus\{\bla_1,\cdots, \bla_\bt\}$, while $F_{\X,\P^1}([S_{\bla_i}^{(\lambda)}])=S_{i,0}^{(\lambda)}:=\bigoplus_{i}S_{i,0}^{(\lambda_ip_i)}$ for $1\leq i\leq \bt$.
	Note that $\haTh_{\star,m}=0$ if $m<0$ and $\haTh_{\star,0}=\frac{1}{\sqq-\sqq^{-1}}$.

	\subsection{The homomorphism $\Omega$}
	\label{subsec:homo}
	Let $\Gamma$ be the star-shaped graph with an involution $\varrho$. 
	Recall the star-shaped graph $\Gamma=T_{p_1,\dots,p_\bt}$ in \eqref{star-shaped}. Let $\tUiD$ be the $\imath$quantum loop algebra of $(\Gamma,\varrho)$. Let $\X$ be the weighted projective line with an involution induced by $\varrho$ (which is still denoted by $\varrho$).  Let 
	$$\I_\varrho:=\{\text{the chosen representatives of $\varrho$-orbits of nodes in $\Gamma$} \}.$$
	Then we have the following theorem.

	\begin{theorem}
		\label{thm:morphi}
		There exists a $\Q(\sqq)$-algebra homomorphism
		\begin{align}
			\Omega: \tUiD_\sqq\longrightarrow \iH(\X_\bfk,\varrho),
		\end{align}
		which sends
		\begin{align}
			\label{eq:mor1}
			&\K_{\star}\mapsto [K_{\co}], \qquad \K_{[i,j]}\mapsto [K_{S_{ij}}], \qquad C\mapsto [K_\de];&
			\\
			\label{eq:mor2}
			&{B_{\star,l}\mapsto \frac{-1}{q-1}[\co(l\vec{c})]},\qquad\Theta_{\star,r} \mapsto {\widehat{\Theta}_{\star,r}}, \qquad H_{\star,r} \mapsto {\widehat{H}_{\star,r}};\\
			\label{eq:mor3}
			& \Theta_{[i,j],r}\mapsto\widehat{\Theta}_{[i,j],r}, \quad H_{[i,j],r}\mapsto \widehat{H}_{[i,j],r}
			,\qquad B_{[i,j],l}\mapsto \begin{cases}{\frac{-1}{q-1}}\haB_{[i,j],l}, \text{ if }[i,j]\in\I_\varrho,\\
				\frac{\sqq}{q-1}\haB_{[i,j],l},\text{ if }[i,j]\notin\I_\varrho;\end{cases}
		\end{align}
		for any $[i,j]\in\II-\{\star\}$, $l\in\Z$, $r>0$.
	\end{theorem}
	
	The proof of this theorem consists of the next 3 sections, where we shall verify that Relations~\eqref{qsiDR0}--\eqref{qsiDR9} are preserved by $\Omega$.
	
	\begin{remark}
		We expect $\Omega$ to be injective; cf. \cite{Sch04, LR21, LR23}. 
		In fact, one can follow the proof of \cite[Theorem 3.2]{LR23} to strengthen Theorem ~\ref{thm:morphi} by showing that $\Omega$ is injective if $\fg$ is of finite or affine type.
	\end{remark}
	
	At the end of this section, let us verify Relations~\eqref{qsiDR0} in the $\imath$Hall algebra.
	
	\begin{lemma}
		For any vertices $\mu,\nu \in \Gamma$, $l\in\Z$, we  have
		\begin{enumerate}
			\item $[K_\de]$ and $[K_\co]$ are central elements;
			\item $\big[[K_\mu],[K_\nu]\big]=0$;
			\item $\big[[K_\mu],\widehat{H}_{\nu,l}\big]=0$;
			\item $[K_\mu] *[\co(l\vec{c})] =[\co(l\vec{c})]*[K_\mu]$,\quad $[K_\mu] *\widehat{B}_{\nu,l}=\sqq^{c_{\mu,\varrho(\nu)}-c_{\mu,\nu}}\widehat{B}_{\nu,l}* [K_\mu]$.
		\end{enumerate}
	\end{lemma}
	
	\begin{proof}
		(2) follows from Proposition \ref{prop:hallbasis}~(1).
		
		(1) follows from \eqref{eq:KXM} by noting that
		$(\de,-)=0$ and $\varrho(\co)=\co$. Similarly, we can prove $[K_\mu] *[\co(l\vec{c})] =[\co(l\vec{c})]*[K_\mu]$ in (4).
		
		For the remaining formulas in (3)--(4),  it is enough to prove for the case $\mu,\nu\in\II$. We can assume $\mu=\nu$ or $\mu=\varrho(\nu)$ by noting that $\Ext^1_\X(M,N)=0=\Hom_\X(M,N)$ if $M\in \scrt_x$, $N\in\scrt_y$ with $x\neq y$.
		
		\underline{Case $\mu=\nu$}. If $\varrho(\mu)=\mu$, then $[K_\mu]$ is central. Otherwise, the formulas follow from the definition of $\widehat{B}_{\nu,l}$ and $\widehat{H}_{\nu,l}$ in \eqref{def:haBXno}--\eqref{def:haphino} and Relations \eqref{DR1a}, \eqref{DR3} in $\tUD(\widehat{\mathfrak{sl}}_n)$.  
		
		\underline{Case $\mu=\varrho(\nu)$}. Only need to consider $\mu\neq \nu$, whose proof is similar to the above.
	\end{proof}

	\section{Relations in tubes and commutative relations}
	\label{sec:Relationtube}

	In this section, we shall verify the relations ~\eqref{qsiDR1}--\eqref{qsiDR9} for the indexes $\mu,\nu\in\II$, where $\{\mu,\nu\}\neq\{\star, [i,1]\}$ for any $1\leq i\leq \bt$.
	More precisely, we will consider the following three cases: $\mu,\nu\in\II-\{\star\}$, or $\mu=\star=\nu$, or $\{\mu,\nu\}=\{\star, [i,j]\}$ for $2\leq j\leq p_i-1$ in \S\ref{subsec:Relations tor}, \S\ref{subsec:Relationstar} and \S\ref{subsec:relationsstarjneq1} respectively.
	

	Let $k_1, k_2, l\in \Z$ and $\mu,\nu \in \II$. Inspiring by \eqref{eq:Skk}--\eqref{eq:Rkk}, we introduce shorthand notations:
	\begin{align}
		\begin{split}
			\widehat{S}(k_1,k_2|l;\mu,\nu)
			&=  \haB_{\mu,k_1} * \haB_{\mu,k_2} *\haB_{\nu,l} -[2]_\sqq \haB_{\mu,k_1} *\haB_{\nu,l} * \haB_{\mu,k_2} + \haB_{\nu,l} *\haB_{\mu,k_1} * \haB_{\mu,k_2},
			\\
			\widehat{\SS}(k_1,k_2|l;\mu,\nu)
			&= \widehat{S}(k_1,k_2|l;\mu,\nu)  + \{k_1 \leftrightarrow k_2 \}.
			\label{eq:haSkk}
		\end{split}
	\end{align}
	We also denote
	\begin{align}
		\notag
		\widehat{R}(k_1,k_2|l; \mu,\nu)
		&=  [K_{k_1\de+\alpha_\mu}]* 
		\Big(-\sum_{p\geq0} \sqq^{2p}  [2]_\sqq \cdot [\haTh _{\mu,k_2-k_1-2p-1},\haB_{\nu,l-1}]_{\sqq^{-2}}*[K_\de]^{p+1}
		\\\notag
		&\quad -\sum_{p\geq 1} \sqq^{2p-1}  [2]_\sqq \cdot [\haB_{\nu,l},\haTh _{\mu,k_2-k_1-2p}]_{\sqq^{-2}} * [K_\de]^{p}
		- [\widehat{B}_{\nu,l}, \haTh _{\mu,k_2-k_1}]_{\sqq^{-2}} \Big),
		\\
		\widehat{ \R}(k_1,k_2|l; \mu,\nu) &=\widehat{ R}(k_1,k_2|l;\mu,\nu) + \{k_1 \leftrightarrow k_2\}.\label{eq:haRkk}
	\end{align}
	
	By direct computations, we have 
	\begin{align}
		\label{eqnSSkk1}
		&\widehat{ \R}(k,k|l; \mu,\nu)=- [\widehat{B}_{\nu,l}, \haTh _{\mu,0}]_{\sqq^{-2}}*[K_{k\de+\alpha_\mu}]=-\sqq^{-1}\widehat{B}_{\nu,l} *[K_{k\de+\alpha_\mu}];\\
		\label{eqnSSkk+1}
		&\widehat{ \R}(k,k+1|l; \mu,\nu)=\Big(-[2]_\sqq \cdot [\haTh _{\mu,0},\haB_{\nu,l-1}]_{\sqq^{-2}}*[K_\de]- [\widehat{B}_{\nu,l}, \haTh _{\mu,1}]_{\sqq^{-2}}\Big)*[K_{k\de+\alpha_\mu}];
	\end{align}
	and for $k_2\geq k_1+2$,
	\begin{align}
		\label{Theta and B k_2>k_1+1}
		\widehat{ \R}(k_1,k_2|l; \mu,\nu)&-q\widehat{ \R}(k_1,k_2-2|l; \mu,\nu)*[K_\de]=\Big(-[2]_\sqq \cdot [\haTh _{\mu,k_2-k_1-1},\haB_{\nu,l-1}]_{\sqq^{-2}}*[K_\de] \\\notag
		&-  [\haB_{\nu,l},\haTh _{\mu,k_2-k_1-2}]_{\sqq^{-2}}* [K_\de]
		- [\widehat{B}_{\nu,l}, \haTh _{\mu,k_2-k_1}]_{\sqq^{-2}} \Big)*[K_{k_1\de+\alpha_\mu}].
	\end{align}

	\subsection{Relations for vertices in $\II-\{\star\}$}
	\label{subsec:Relations tor}
	For any vertex $\mu=[i,j]\in\Gamma$, let $\varrho(\mu)=[\varrho(i),\varrho(j)]$.
	
	\begin{proposition}
		For any $\mu=[i,j],\nu=[i',j']\in\II-\{\star\}$, the following relations hold in $\iH(\X_\bfk,\varrho)$, where $k,l\in\Z$, $m,n\geq 1$.
		\begin{align}
			\label{HHmn}
			&[\widehat{H}_{\mu,m},\widehat{H}_{\nu,n}]=0,\\
			&[\widehat{H}_{\mu,m},\widehat{B}_{\nu,l}]=\frac{[mc_{\mu\nu}]_\sqq}{m} \widehat{B}_{\nu,l+m}-\frac{[mc_{\mu,\varrho(\nu)}]_\sqq}{m} \widehat{B}_{\nu,l-m}*[K_{m\de}],\label{humbvl}
			\\
			&[\widehat{B}_{\mu,k} ,\widehat{B}_{\nu,l}]=0,   \text{ if }c_{\mu\nu}=0, \text{ and }\mu\neq \varrho(\nu),  
			\\
			\label{haBBitaui}
			&[\widehat{B}_{\mu,k} ,\widehat{B}_{\varrho(\mu),l}]= [K_{S_{\varrho(\mu)}}]* [K_{l \delta}]* \widehat{\TH}_{\mu,k-l}- [K_{S_{\mu}}]* [K_{k\delta}]* \widehat{\Theta}_{\varrho(\mu),l-k},   \text{ if }\mu\neq \varrho(\mu), 
			\\
			\label{BBmuneqrhonu}
			&[\widehat{B}_{\mu,k}, \widehat{B}_{\nu,l+1}]_{\sqq}  -\sqq [\widehat{B}_{\mu,k+1}, \widehat{B}_{\nu,l}]_{\sqq}=0, \text{ if } c_{\mu\nu}=-1, \text{ and }\mu\neq \varrho(\nu),
			\\ 
			\label{Bmuk Bmul}
			&[\widehat{B}_{\mu,k}, \widehat{B}_{\mu,l+1}]_{\sqq^{-2}}  -\sqq^{-2} [\widehat{B}_{\mu,k+1}, \widehat{B}_{\mu,l}]_{\sqq^{2}}
			=(1-q)^2
			\Big(\sqq^{-2}\widehat{\Theta}_{\mu,l-k+1} *[K_{k\de+\alpha_{ij}}]
			\\
			&\qquad-\sqq^{-4}\widehat{\Theta}_{\mu,l-k-1} * [K_{(k+1)\de+\alpha_{ij}}]
			+\sqq^{-2}\widehat{\Theta}_{\mu,k-l+1}* [K_{l\de+\alpha_{ij}}]-\sqq^{-4}\widehat{\Theta}_{\mu,k-l-1}* [K_{(l+1)\de+\alpha_{ij}}]\Big), \notag
			\\
			\notag
			&\hspace{12cm}{\text{ if }\varrho(\mu)=\mu},
			\\
			\label{SerrehaBB}
			&    \widehat{\SS}(k_1,k_2|l; \mu,\nu) = (1-q)^2 \widehat{\R}(k_1,k_2|l; \mu,\nu), \qquad\qquad\text{ if }c_{\mu\nu}=-1, \varrho(\mu)=\mu,
			\\\label{SerrehaBB1}
			&\widehat{\SS}(k_1,k_2|l;\mu,\nu)
			= 0,\qquad\qquad\qquad\qquad\qquad\qquad\qquad
			\text{ if }c_{\mu,\nu}=-1,   \varrho(\mu)\neq \mu.
		\end{align}
	\end{proposition}
	
	\begin{proof}
		
		If $i=i'$, and $\varrho(\mu)=\mu$ (also $\varrho(\nu)=\nu$), then the proof is completely same as \cite[Proposition 7.1]{LR21}.

		If $\mu\neq \varrho(\mu)$ and $i=i'$ or $i=\varrho(i')$, then we have 
		the following algebra embeddings
		\[
		\tUD_\sqq(\widehat{\mathfrak{sl}}_{p_i})\stackrel{\Omega_{C_{p_i}}\,\,}{\longrightarrow}  \cs\cd\widetilde{\ch}_{\Z_2}(\bfk C_{p_i}) \stackrel{\iota_i}{\longrightarrow} \iH(\X_\bfk,\varrho).
		\]
		Now it follows from \eqref{def:haBXno}--\eqref{def:haphino} that all the relations \eqref{HHmn}--\eqref{BBmuneqrhonu} and \eqref{SerrehaBB1} hold in this case. 
		
		It remains to consider the case $\varrho(i')\neq i\neq i'$.
		Moreover, for any torsion sheaves $X,Y$ supported at distinct points, we have $\Ext^1_\X(X,Y)=0=\Hom_\X(X,Y)$.  It follows from Lemma \ref{lem:Extcrho} that $\Ext^1_{\cc_\varrho(\coh(\X))}(X,Y)\cong \Ext^1_\X(X,Y)\oplus \Hom_\X(X,\varrho(Y))$, and $\Hom_{\cc_\varrho(\coh(\X))}(X,Y)=\Hom_\X(X,Y)$. Hence $[X]*[Y]=[X\oplus Y]=[Y]*[X]$ in $\widetilde{\ch}(\cc_\varrho(\coh(\X)))$ and in $\iH(\X_\bfk,\varrho)$. Therefore, all the relations \eqref{HHmn}--\eqref{SerrehaBB1} hold in this case.
	\end{proof}
	
	\subsection{Relations for $\star$}
	\label{subsec:Relationstar}

	The following proposition shows that \eqref{qsiDR1}--\eqref{qsiDR9} hold  in $\iH(\X_\bfk,\varrho)$ for $\mu=\nu=\star$.
	\begin{proposition}
		\label{prop:OO}
		For any $m,n\geq1$, $k,l\in\Z$, the following relations hold in $\iH(\X_\bfk,\varrho)$:
		\begin{align}
			\label{eq:HaDr1}
			[\haH_{\star, m},\haH_{\star,n}]&=0,\\
			\label{eq:HaDr2}
			\big[\haH_{\star,m}, [\co(l\vec{c})]\big] &=\frac{[2m]_\sqq}{m} [\co((l+m)\vec{c})]-\frac{[2m]_\sqq}{m} [\co((l-m)\vec{c})]*[K_{m\de}],
		\end{align}
		\begin{align}
			\label{eq:HaDr3}
			\big[[\co(k\vec{c})], [\co&((l+1)\vec{c})]\big]_{\sqq^{-2}}  -\sqq^{-2} \big[[\co((k+1)\vec{c})], [\co(l\vec{c})]\big]_{\sqq^{2}}
			\\\notag
			&= (1-q)^2\Big(\sqq^{-2}\haTh_{\star,l-k+1} *[K_\de]^k* [K_{\co}]-\sqq^{-4} \haTh_{\star,l-k-1}* [K_\de]^{k+1} *[K_{\co}]
			\\\notag
			&+\sqq^{-2}\haTh_{\star,k-l+1}*[K_\de]^l *[K_{\co}]-\sqq^{-4} \haTh_{\star,k-l-1}*[ K_\de]^{l+1} *[K_\co]\Big).
		\end{align}
	\end{proposition}
	
	\begin{proof}
		By Proposition \ref{prop:autogroupP1}, we can view $\varrho$ to be an involution of $\coh(\P^1_\bfk)$. 
		Recall from \S \ref{Embedding from projective line to weighted projective line}
		that there is an embedding \begin{equation*}
			F_{\X,\P^1}:\iH(\P^1_\bfk,\varrho) \longrightarrow \iH(\X_\bfk,\varrho),
		\end{equation*} which sends
		$[\co_{\P^1}(l)]\mapsto[\co(l\vec{c})]$, $[S_{\bla_i}^{(r)}]\mapsto [S_{i,0}^{(rp_i)}]$, $\haT_{r}\mapsto \haT_{\star,r}$ and $\haH_r\mapsto \haH_{\star,r}$
		for any $l\in\Z, r\geq 1$.
		Then the formulas follow from Proposition \ref{prop:P1-Onsager} and the definition of $\tUiD_v(\widehat{\mathfrak{sl}}_2)$.
	\end{proof}
	
	\subsection{Relations between $\star$ and $[i,j]$ with $j\neq 1$}
	\label{subsec:relationsstarjneq1}

	In this subsection, we verify the relations \eqref{qsiDR1}--\eqref{qsiDR9} between $\star$ and $[i,j]$ with $j\neq 1$. In fact, we only need to consider \eqref{qsiDR1}, \eqref{qsiDR2} and \eqref{qsiDR7}. If $\varrho(i)=i$, then the proof is completely the same to \cite[\S7]{LR21}, hence omitted here. So we always assume $\varrho(i)\neq i$ in the following.
	
	First we give some preliminary results.
	
	\begin{lemma}
		\label{middle and ending terms same}
		For any $1\leq i\leq \bt$ and $1\leq j< k\leq p_i$ and $l\in\Z$, the following formulas hold in $\iH(\X_\bfk,\varrho)$:
		\begin{align}
			\label{eq:OSij}
			\big[[\co(l\vec{c})], [S_{i,j}^{(k)}]\big]=&\big[\co(l\vec{c})], [S_{i,j}^{(j)}\oplus S_{i,0}^{(k-j)}]\big],
			\\
			\label{eq:Sipiij}
			\big[[S_{i,0}^{(rp_i)}], [S_{i',j}^{(k)}]\big]=&\big[[S_{i,0}^{(rp_i)}], [S_{i',j}^{(j)}\oplus S_{i',0}^{(k-j)}]\big],\quad\text{ for }i'\in\{i,\varrho(i)\},r\geq 1.
		\end{align}
	\end{lemma}
	
	\begin{proof}
		Observe that we have the following exact sequence
		\begin{equation}
			\label{exact sequence for sij}
			0\rightarrow S_{i,0}^{(k-j)}\longrightarrow S_{ij}^{(k)} \longrightarrow S_{ij}^{(j)}\rightarrow0.
		\end{equation}
		
		By applying $\Hom_{\X}(\co(l\vec{c}),-)$ and $\Hom_{\X}(-,\co(l\vec{c}))$ to (\ref{exact sequence for sij}), we know that
		\begin{align*}
			\Hom_\X(\co(l\vec{c}),S_{ij}^{(k)})\cong \Hom_\X(\co(l\vec{c}),S_{i,0}^{(k-j)})\cong \bfk,\\
			\Ext^1_\X(S_{ij}^{(k)}, \co(l\vec{c}))\cong \Ext_\X^1(S_{ij}^{(j)}, \co(l\vec{c}))\cong \bfk.
		\end{align*}
		Let $X\in \{S_{ij}^{(k)}, S_{ij}^{(j)}\oplus S_{i,0}^{(k-j)}\}$.
		For any non-zero morphism $f:\co(l\vec{c})\rightarrow \varrho(X)$, we have
		$$\Im(f)=S_{\varrho(i),0}^{(k-j)}, \quad\Ker(f)=\co(l\vec{c}-(k-j)\vec{x}_{\varrho(i)}),\quad \coker(\varrho(f))=S_{ij}^{(j)},$$ while the non-trivial extension in $\Ext^1_\X(X, \co(l\vec{c}))$ is given by
		$$0\rightarrow \co(l\vec{c})\longrightarrow \co(l\vec{c}+j\vec{x}_i)\oplus S_{i,0}^{(k-j)} \longrightarrow X\rightarrow0.$$
		Therefore, by \eqref{Hallmult1} we have
		\begin{align*}
			&\big[[\co(l\vec{c})], [S_{ij}^{(k)}]\big]\\
			=&[\co(l\vec{c})]\ast [S_{ij}^{(k)}]-[S_{ij}^{(k)}]\ast[\co(l\vec{c})]\\
			=&\sqq^{\langle \co(l\vec{c}), S_{ij}^{(k)}\rangle}\cdot\Big(\frac{1}{q}[\co(l\vec{c})\oplus S_{ij}^{(k)}]+\frac{q-1}{q}[\co(l\vec{c}-(k-j)\vec{x}_{\varrho(i)})\oplus S_{ij}^{(j)}]\ast [K_{S_{\varrho(i),0}^{(k-j)}}]\Big)\\
			&\quad -\sqq^{\langle S_{ij}^{(k)},\co(l\vec{c})\rangle}\cdot\Big([\co(l\vec{c})\oplus S_{ij}^{(k)}]+(q-1)[\co(l\vec{c}+j\vec{x}_i)\oplus S_{i,0}^{(k-j)}]\Big)\\
			=&(\sqq-\sqq^{-1})\Big([\co(l\vec{c}-(k-j)\vec{x}_{\varrho(i)})\oplus S_{ij}^{(j)}]* [K_{S_{\varrho(i),0}^{(k-j)}}]-[\co(l\vec{c}+j\vec{x}_i)\oplus S_{i,0}^{(k-j)}]\Big).
		\end{align*}
		Similarly,
		\begin{align*}
			&\big[[\co(l\vec{c})], [S_{ij}^{(j)}\oplus S_{i,0}^{(k-j)}]\big]\\
			=&(\sqq-\sqq^{-1})\Big([\co(l\vec{c}-(k-j)\vec{x}_{\varrho(i)}))\oplus S_{ij}^{(j)}]\ast [K_{S_{\varrho(i),0}^{(k-j)}}]-[\co(l\vec{c}+j\vec{x}_i)\oplus S_{i,0}^{(k-j)}]\Big),
		\end{align*}
		and then \eqref{eq:OSij} follows.

		For \eqref{eq:Sipiij}, observe that the composition factors of $S_{ij}^{(j)}$ is $S_{i,1}, S_{i,2},\cdots, S_{ij}$ (from socle to top), while the composition factors of $S_{i,0}^{(k-j)}$ is $S_{i,p_i+1-k+j},\cdots, S_{i,p_i-1}, S_{i,0}$. Hence
		\begin{align*}
			\Hom_\X(S_{i,0}^{(rp_i)},S_{ij}^{(j)})=0= \Ext^1_\X(S_{i,0}^{(rp_i)},S_{ij}^{(j)}),
			\quad
			\Hom_\X(S_{ij}^{(k-j)}, S_{i,0}^{(rp_i)})=0=\Ext^1_\X(S_{ij}^{(k-j)}, S_{i,0}^{(rp_i)}).
		\end{align*}
		By applying $\Hom_{\X}(S_{i,0}^{(rp_i)},-)$ and $\Hom_{\X}(-, S_{i,0}^{(rp_i)})$ to (\ref{exact sequence for sij}), we obtain
		\begin{align*}
			&\Hom_\X(S_{i,0}^{(rp_i)},S_{ij}^{(k)})\cong\Hom_\X(S_{i,0}^{(rp_i)},S_{i,0}^{(k-j)}),
			&&\Ext^1_\X(S_{i,0}^{(rp_i)},S_{ij}^{(k)})\cong\Ext^1_\X(S_{i,0}^{(rp_i)},S_{i,0}^{(k-j)})),
			\\
			&\Hom_\X(S_{ij}^{(k)}, S_{i,0}^{(rp_i)})\cong \Hom_\X(S_{ij}^{(j)}, S_{i,0}^{(rp_i)}),
			&&\Ext^1_\X(S_{ij}^{(k)}, S_{i,0}^{(rp_i)})\cong \Ext^1_\X(S_{ij}^{(j)}, S_{i,0}^{(rp_i)}),
		\end{align*}
		which are of dimension $1$.

		For $i'=i$, we have 
		\begin{align*}
			\big[[S_{i,0}^{(rp_i)}], [S_{ij}^{(k)}]\big]=&\sqq^{-1}\Big([S_{i,0}^{(rp_i)}\oplus S_{ij}^{(k)}]+(q-1)[S_{i,0}^{(rp_i+k-j)}\oplus S_{ij}^{(j)}] \Big)
			\\
			&\quad-\sqq^{-1}\Big([S_{i,0}^{(rp_i)}\oplus S_{ij}^{(k)}]+(q-1)[S_{ij}^{(rp_i+j)}\oplus S_{i,0}^{(k-j)}] \Big)
			\\
			=&(\sqq-\sqq^{-1})\Big([S_{i,0}^{(rp_i+k-j)}\oplus S_{ij}^{(j)}] - [S_{ij}^{(rp_i+j)}\oplus S_{i,0}^{(k-j)}]\Big).
		\end{align*}
		Similarly, one can get that 
		\begin{align*}
			\big[[S_{i,0}^{(rp_i)}], [S_{ij}^{(j)}\oplus S_{i,0}^{(k-j)}]\big]=(\sqq-\sqq^{-1})\Big([S_{i,0}^{(rp_i+k-j)}\oplus S_{ij}^{(j)}] - [S_{ij}^{(rp_i+j)}\oplus S_{i,0}^{(k-j)}]\Big).
		\end{align*}

		For $i'=\varrho(i)$ ($i'\neq i$ by our assumption),  we have
		\begin{align*}
			\big[[S_{i,0}^{(rp_i)}], [S_{\varrho(i),j}^{(k)}]\big]=&[S_{i,0}^{(rp_i)}\oplus S_{\varrho(i),j}^{(k)}]+(q-1)[S_{i,p_i-k+j}^{(rp_i-k+j)}\oplus S_{\varrho(i),j}^{(j)}\oplus K_{S_{i,0}^{(k-j)}}] \Big)
			\\
			&\quad-\Big([S_{i,0}^{(rp_i)}\oplus S_{\varrho(i),j}^{(k)}]+(q-1)[S_{\varrho(i),0}^{(k-j)}\oplus S_{i,0}^{(rp_i-j)}\oplus K_{S_{\varrho(i),j}^{(j)}}] \Big)
			\\
			=&(q-1)\Big([S_{i,p_i-k+j}^{(rp_i-k+j)}\oplus S_{\varrho(i),j}^{(j)}\oplus K_{S_{i,0}^{(k-j)}}] -[S_{\varrho(i),0}^{(k-j)}\oplus S_{i,0}^{(rp_i-j)}\oplus K_{S_{\varrho(i),j}^{(j)}}]\Big).
		\end{align*}
		
		Similarly, we can get that
		\begin{align*}
			\big[[S_{i,0}^{(rp_i)}], &[S_{\varrho(i),j}^{(j)}\oplus S_{\varrho(i),0}^{(k-j)}]\big]
			\\
			&=(q-1)\Big([S_{i,p_i-k+j}^{(rp_i-k+j)}\oplus S_{\varrho(i),j}^{(j)}\oplus K_{S_{i,0}^{(k-j)}}] -[S_{\varrho(i),0}^{(k-j)}\oplus S_{i,0}^{(rp_i-j)}\oplus K_{S_{\varrho(i),j}^{(j)}}]\Big).
		\end{align*}
	\end{proof}
	
	The set $\cm_{j,\alpha}$ defined in \eqref{def:Mjalpha} can be viewed as a set of isoclasses of $\coh(\X)$ via $\iota_i$, which is denoted by $\cm_{[i,j],\alpha}$; see \S\ref{subsec:embeddingtube}.
	Recall from Proposition \ref{prop:DrGenA} that
	\begin{align}
		\label{haThij1}
		\widehat{\Theta}_{[i,j],1}=\widehat{H}_{[i,j],1}=\pi_{[i,j+1],1}-(\sqq+\sqq^{-1})\pi_{[i,j],1}+\pi_{[i,j-1],1},
	\end{align}
	where \begin{align*}
		\pi_{[i,j],1}=\frac{-\sqq^{-j}}{\sqq-\sqq^{-1}}\sum\limits_{\cm_{[i,j],\delta}}(-1)^{\dim\End(M)}[M].
	\end{align*}

	\begin{lemma} For any $l\in\mathbb{Z}$, $1\leq i\leq \bt$ and $1\leq j\leq p_i-1$, we have
		\begin{align}\label{pi_ij,1+ anc co}
			\big[\pi_{[i,j],1}, [\co(l\vec{c})]\big]=&\frac{\sqq^{-j}}{\sqq-\sqq^{-1}}\big[[S_{i,0}^{(p_i)}], [\co(l\vec{c})]\big],
			\\
			\label{pi and S(rp)}
			\big[\pi_{[i,j],1}, [S_{i',0}^{(rp_i)}]\big]=&0, \text{ for }i'\in\{i,\varrho(i)\}.
		\end{align}
	\end{lemma}
	
	\begin{proof}
		For any $1\leq j',k\leq p_i$, the torsion sheaf $S_{i,j'}^{(k)}$ is uniserial, and $S_{i,1}$ is a composition factor of $S_{i,j'}^{(k)}$ if and only if $k\geq j'$. Hence, for any $[M]\in\cm_{[i,j],\delta}$, there exists a unique direct summand $S_{i,j'}^{(k)}$ of $M$ with $k\geq j'$, namely, $M\cong S_{i,j'}^{(k)}\oplus M'$ for some $M'$.
		We emphasize that $$[S_{i,j'}^{(k)}\oplus M']\in\cm_{[i,j],\delta}\quad \text{iff} \quad[S_{i,j'}^{(j')}\oplus S_{i,0}^{(k-j')}\oplus M']\in\cm_{[i,j],\delta}.$$
		
		For any $k>j'$, we have 
		$$\dim\End(S_{i,j'}^{(j')}\oplus S_{i,0}^{(k-j')}\oplus M')=\dim\End(S_{i,j'}^{(k)}\oplus M')+1.$$ 
		Moreover, there are no non-trivial homomorphisms and extensions between $\co(l\vec{c})$ and $M'$ since both $S_{i,0}$ and $S_{i,1}$ are not composition factors of $M'$. Then by \eqref{eq:OSij} we get
		$$\big[[S_{i,j'}^{(k)}\oplus M'], [\co(l\vec{c})]\big]=\big[[S_{i,j'}^{(j')}\oplus S_{i,0}^{(k-j')}\oplus M'], [\co(l\vec{c})]\big].$$
		Therefore, only the case $k=j'=p_i$, i.e., the term $[S_{i,0}^{(p_i)}]$ in $\pi_{[i,j],1}$, has non-trivial contribution to $\big[\pi_{[i,j],1}, [\co(l\vec{c})]\big]$.
		Then (\ref{pi_ij,1+ anc co}) follows. 
		
		For \eqref{pi and S(rp)}, by similar arguments as above and using (\ref{eq:Sipiij}), we obtain that
		\begin{align*}
			\big[\pi_{[i,j],1}, [S_{i',0}^{(rp_i)}]\big]=\frac{\sqq^{-j}}{\sqq-\sqq^{-1}}\big[[S_{i,0}^{(p_i)}], [S_{i',0}^{(rp_i)}]\big]=0,
		\end{align*}
		where the last equality follows from the proof of Proposition \ref{prop:P1-Onsager} and the embedding (\ref{the embedding functor F on algebra}).
	\end{proof}
	
	\begin{lemma} For any $1\leq i\leq \bt$ and $2\leq j\leq p_i-1$, we have
		\begin{align}\label{x * and h_ij,1}
			\big[\widehat{H}_{[i,j], 1},[\co(k\vec{c})]\big] =&0, \text{ for }k\in\Z,
			\\\label{haThetaij1starm}
			\big[\widehat{\Theta}_{[i,j],1}, \widehat{\Theta}_{\star,m}\big]=&0, \text{ for } m\geq 1.
		\end{align}
	\end{lemma}
	
	\begin{proof}
		By (\ref{pi_ij,1+ anc co}), we have
		\begin{align*}
			\big[\widehat{H}_{[i,j],1}, [\co(k\vec{c})]\big]
			=&\big[\pi_{[i,j+1],1}-(\sqq+\sqq^{-1})\pi_{[i,j],1}+\pi_{[i,j-1],1},[\co(k\vec{c})]\big]\\
			=&\big(\frac{\sqq^{-(j+1)}}{\sqq-\sqq^{-1}}-(\sqq+\sqq^{-1})\cdot\frac{\sqq^{-j}}
			{\sqq-\sqq^{-1}}+\frac{\sqq^{-(j-1)}}{\sqq-\sqq^{-1}}\big)\big[[S_{i,0}^{(p_i)}], [\co(k\vec{c})]\big]\\
			=&0.
		\end{align*}
		
		For \eqref{haThetaij1starm}, observe that there are no non-zero homomorphisms and extensions between distinct tubes in $\coh(\X)$, then by the definition of $\widehat{\Theta}_{\star,m}$, it suffices to show that
		\begin{align*}\big[\widehat{\Theta}_{[i,j],1}, [S_{i',0}^{(rp_i)}]\big]=0, \quad\forall\; r\geq 1,
		\end{align*}
		for $i'\in\{i,\varrho(i)\}$, 
		which then follows from (\ref{pi and S(rp)}) immediately.
	\end{proof}
	
	\begin{proposition}
		For any $k,t\in\mathbb{Z}$, $l,m>0$, $1\leq i\leq \bt$ and $2\leq j\leq p_i-1$, we have
		
		\begin{align}\label{x* and x_ij}
			&\big[[\co(k\vec{c})],\widehat{B}_{[i,j], t}\big] = 0,
			\\\label{x* and h_ij}
			&\big[[\co(k\vec{c})], \widehat{H}_{[i,j],m}\big] =0,
			\\\label{h* and x_ij}
			&[\widehat{H}_{\star ,l},\widehat{B}_{[i,j], k}] =0,
			\\\label{h* and h_ij}
			&[\widehat{H}_{\star, l},\widehat{H}_{[i,j], m}] = 0.
		\end{align}
	\end{proposition}

	\begin{proof}
		For any $2\leq j\leq p_i-1$, there are no nonzero homomorphisms and extensions between $\co(k\vec{c})$ and $S_{ij}$. Hence
		$$\big[[\co(k\vec{c})],\widehat{B}_{[i,j], 0}\big]=0.$$
		Recall from \eqref{humbvl} that for any $t$,
		\begin{align}
			\label{from xij,1 to xij,m+1}
			&[\widehat{H}_{[i,j],1},\widehat{B}_{[i,j],t}] = [2]_\sqq  \widehat{B}_{[i,j],t+1},\qquad 
			[\widehat{H}_{[\varrho(i),j],1},\widehat{B}_{[i,j],t}] = -[2]_\sqq [K_\de]* \widehat{B}_{[i,j],t-1}.
		\end{align}
		Recall from \eqref{x * and h_ij,1} that $\big[\widehat{H}_{[i,j], 1},[\co(k\vec{c})]\big] =0=\big[\widehat{H}_{[\varrho(i),j], 1},[\co(k\vec{c})]\big]$. Applying $\big[[\co(k\vec{c})],-\big]$ to \eqref{from xij,1 to xij,m+1}, by induction on $|t|$ we have $\big[[\co(k\vec{c})],\widehat{B}_{[i,j], t}\big] = 0$.
		This proves \eqref{x* and x_ij}.
		
		Recall from \eqref{haBBitaui} that for $k>t$: 
		\begin{align}
			\label{haBBitaui>1}
			[\widehat{B}_{[i,j],k} ,\widehat{B}_{[\varrho(i),j],t}]= [K_{S_{[\varrho(i),j]}}]* [K_{t \delta}]* \widehat{\TH}_{[i,j],k-t}.
		\end{align}
		It follows from
		\eqref{x* and x_ij} that
		$\big[[\co(k\vec{c})], \widehat{\Theta}_{[i,j],m}\big] =0$ for any $m>0$, and then \eqref{x* and h_ij} follows.
		
		For \eqref{h* and x_ij}, it is equivalent to prove that $\big[\widehat{\Theta}_{\star, l},\haB_{[i,j],k}\big]=0$.
		By using similar arguments as for \eqref{x* and x_ij}, thanks to \eqref{haThetaij1starm}, we only need to prove the formula for $k=0$. For $k=0$, $\Big[\widehat{\Theta}_{\star, l},[S_{ij}]\Big]=0$  is obvious since $\big[ [S_{i,0}^{(rp_i)}], [S_{ij}] \big]=0=\big[[S_{\varrho(i),0}^{(rp_i)}], [S_{ij}] \big]$ for $2\leq j\leq p_i-1$, $r\geq1$. Hence \eqref{h* and x_ij} follows.
		
		Finally, \eqref{h* and h_ij} follows from \eqref{h* and x_ij} by using \eqref{haBBitaui>1}.
		
		The proof is completed. 
	\end{proof}

	\section{Relations between $\star$ and $[i,1]$, I}
	\label{sec:Relationsstari1 I}
	
	In this section and \S\ref{sec:Relationsstari1 II}, we shall verify the relations \eqref{qsiDR1}--\eqref{qsiDR9} between $\star$ and $[i,1]$. First, we shall describe the root vectors at $[i,1]$. 
	
	\subsection{Root vectors at $[i,1]$}
	\label{subsec:realroot}

	For any $r>0$, denote by
	\begin{align}
		\label{eq:realroots11set}
		\begin{split}
			&\cm_{r\delta+\alpha_{i1}}=\{[S_{i,1}^{(bp_i+1)}\oplus S_{i,0}^{(\nu)}]\mid b\geq 0, \;|\nu|+b=r\};\\
			&\cm_{r\delta-\alpha_{i1}}=\{[S_{i,0}^{(ap_i-1)}\oplus S_{i,0}^{(\nu)}]\mid a\geq 1, \;|\nu|+a=r\};
		\end{split}
	\end{align}
	and set
	\begin{align}
		\label{eq:Mrde+alpha11}
		[M_{r\delta\pm\alpha_{i1}}]:=\sum\limits_{[M]\in\cm_{r\delta\pm\alpha_{i1}}}\bn(\ell(M)-1)[\![M]\!],
	\end{align}
	where
	\[[\![M]\!]=\frac{[M]}{|\Aut(M)|},\qquad\qquad\bn(l)=\prod_{i=1}^l(1-\sqq^{2i}),\,\, \forall l\geq1.\]
	
	We also set
	\begin{align}
		\label{def:Mr}
		\cm_{i,r\de}:=&\big\{[S_{i,1}^{(\nua  p_1)}\oplus S_{i,0}^{(\nu)}]\mid a\geq 1, a+|\nu|=r\big\}
		\\\notag
		&\bigcup\big\{[S_{i,0}^{(\nua  p_i-1)}\oplus S_{i,1}^{(b  p_i+1)}\oplus S_{i,0}^{(\nu)}]\mid  a+b+|\nu|=r\big\}.
	\end{align}

	\begin{proposition}
		\label{prop:realroot}
		For any $r>0$, we have
		\begin{align}
			\label{realroot}
			\haB_{[i,1],r}= (q-1)[M_{r\de+\alpha_{i1}}], \qquad \haB_{[i,1],-r}=(1-q) \sqq^{1-\de_{i,\varrho(i)}}[K_{-r\de+\alpha_{i1}}]*[M_{r\de-\alpha_{\varrho(i),1}}].
		\end{align}
		\begin{align}
			\label{eq:THeta11r}
			\haTh_{[i,1],r}=
			&\frac{\sqq}{q-1} \sum\limits_{|\lambda|=r}\bn(\ell(\lambda))  [\![S_{i,0}^{(\lambda)}]\!]
			+\sqq^{-1} \sum\limits_{[M]\in\cm_{i,r\de}}\bn(\ell(M)-1) [\![M]\!].
		\end{align}
	\end{proposition}

	\begin{proof}
		If $\varrho(i)=i$, it is equivalent to proving the formulas in $\iH(\bfk C_{p_i},\Id)$, which is given by \cite[Propositions 9.1, 9.2]{LR21}.
		
		If $\varrho(i)\neq i$, it is enough to compute the formulas in $\cs\cd\widetilde{\ch}_{\Z_2}(\bfk C_{p_i})$, which is given by Propositions \ref{prop:realroot1}--\ref{prop:imageroot}.
	\end{proof}
	
	

	In this section, we shall prove the following formulas in $\iH(\X_\bfk,\varrho)$:
	\begin{align}
		\label{coBi1}
		&\big[[\co(k\vec{c})], \widehat{B}_{[i,1],l+1}\big]_{\sqq}  -\sqq \big[[\co((k+1)\vec{c})], \widehat{B}_{[i,1],l}\big]_{\sqq^{-1}} =0,
		\\
		\label{Serrestari1}
		&\widehat{\SS}(k_1,k_2\mid l;\star, [i,1])=(1-q)^2\widehat{\R}(k_1,k_2\mid l;\star, [i,1]),
		\\
		\label{humstarbvl}
		&[\widehat{H}_{\star,m},\widehat{B}_{[i,1],l}]=-\frac{[m]_\sqq}{m} \widehat{B}_{[i,1],l+m}+\frac{[m]_\sqq}{m} \widehat{B}_{[i,1],l-m}*[K_{m\de}],
		\\
		\label{eq:HHstari1}
		&[\widehat{H}_{\star,m},\widehat{H}_{[i,1],r}]=0,
	\end{align}
	for any $k, k_1,k_2,l\in\Z$, $m,r>0$ and $1\leq i\leq \bt$.

	In fact, if $\varrho(i)=i$, these formulas can be derived by the same computations as \cite[\S8]{LR21}.
	So we can assume $\varrho(i)\neq i$ in the remaining of this section. 
	Without loss of generality, we assume $i=1$ and $\varrho(i)=2$ throughout this section.

	\subsection{The relation \eqref{coBi1}}
	
	\begin{proposition}
		\label{prop:coBil}
		We have \eqref{coBi1} holds for any 
		$k,l\in\Z$.
	\end{proposition}

	\begin{proof}
		Obviously, \eqref{coBi1} is equivalent to
		\begin{align}
			\label{coBi1equiv}
			\big[[\co(k\vec{c})], \widehat{B}_{[1,1],l+1}\big]_{\sqq}  + \big[ \widehat{B}_{[1,1],l},[\co((k+1)\vec{c})]\big]_{\sqq} =0,\qquad\forall\,\, k,l\in\Z.
		\end{align}
		
		First, we prove \eqref{coBi1equiv} for $l=0$. 
		Recall from Proposition \ref{prop:realroot} that
		$$\haB_{[1,1],0}=[S_{1,1}],\text{ and } \haB_{[1,1],1}=\frac{1}{q}[S_{1,1}^{(p_1+1)}] -\frac{1}{q}[S_{1,1}\oplus S_{1,0}^{(p_1)}].$$

		We have
		\begin{align*}
			\big[[S_{1,1}],[\co((l+1)\vec{c})]\big]_\sqq=(\sqq-\sqq^{-1})\Big([\co((l+1)\vec{c}+\vec{x}_1)]-[\co((l+1)\vec{c})\oplus S_{1,1}]\Big).
		\end{align*}
		On the other hand, by definition we have
		\begin{align*}
			&\big[[S_{1,1}^{(p_1+1)}],[\co(l\vec{c})]\big]_{\sqq^{-1}}\\
			=&[S_{1,1}^{(p_1+1)}]*[\co(l\vec{c})]-\sqq^{-1} [\co(l\vec{c})]*[S_{1,1}^{(p_1+1)}]\\
			=&\sqq^{-2}\Big([S_{1,1}^{(p_i+1)}\oplus \co(l\vec{c})]+(q-1) [\co(l\vec{c}+\vec{x}_1)\oplus S_{1,0}^{(p_1)}]+(q^2-q)[\co((l+1)\vec{c}+\vec{x}_1)]\Big)\\
			&-\sqq^{-1}\cdot \frac{\sqq}{q} \Big([S_{1,1}^{(p_1+1)}\oplus \co(l\vec{c})]+(q-1)[S_{1,1}\oplus \co((l-1)\vec{c})]*[K_\de]\Big)\\
			=&(q-1)[\co((l+1)\vec{c}+\vec{x}_1)]+\frac{q-1}{q}[\co(l\vec{c}+\vec{x}_1)\oplus S_{1,0}^{(p_1)}]-\frac{q-1}{q}[S_{1,1}\oplus \co((l-1)\vec{c})]*[K_\de],
		\end{align*}
		and
		\begin{align*}
			&\big[[S_{1,1}\oplus S_{1,0}^{(p_1)}],[\co(l\vec{c})]\big]_{\sqq^{-1}}\\
			=&[S_{1,1}\oplus S_{1,0}^{(p_i)}]*[\co(l\vec{c})]-\sqq^{-1} [\co(l\vec{c})]*[S_{1,1}\oplus S_{1,0}^{(p_1)}]\\
			=&\sqq^{-2}\Big([S_{1,1}\oplus S_{1,0}^{(p_1)}\oplus \co(l\vec{c})]+(q-1) [\co(l\vec{c}+\vec{x}_1)\oplus S_{1,0}^{(p_1)}]+(q^2-q)[\co((l+1)\vec{c})\oplus S_{1,1}]\Big)\\
			&-\sqq^{-1}\cdot \frac{\sqq}{q} \Big([S_{1,1}\oplus S_{1,0}^{(p_1)}\oplus \co(l\vec{c})]+(q-1)[S_{1,1}\oplus \co((l-1)\vec{c})]*[K_\de]\Big)\\
			=&(q-1)[\co((l+1)\vec{c})\oplus S_{1,1}]+\frac{q-1}{q}[\co(l\vec{c}+\vec{x}_1)\oplus S_{1,0}^{(p_1)}]-\frac{q-1}{q}[S_{1,1}\oplus \co((l-1)\vec{c})]*[K_\de].
		\end{align*}
		Hence, by combining the above two formulas we obtain
		\begin{align*}
			\big[[S_{1,1}^{(p_1+1)}]-[S_{1,1}\oplus S_{1,0}^{(p_1)}],[\co(l\vec{c})]\big]_{\sqq^{-1}}= &(q-1)\Big([\co((l+1)\vec{c}+\vec{x}_1)]-[\co((l+1)\vec{c})\oplus S_{1,1}]\Big)\\
			=&\sqq\big[[S_{1,1}],[\co((l+1)\vec{c})]\big]_\sqq,
		\end{align*} 
		which is \eqref{coBi1equiv} for $l=0$.

		For any $l\in\Z$, from \eqref{humbvl}, we have
		\begin{align*}
			&\big[\widehat{\Theta}_{[1,1],1}, \haB_{[1,1],l}\big]=[2]_\sqq\haB_{[1,1],l+1},
			&\big[\widehat{\Theta}_{[2,1],1}, \haB_{[1,1],l}\big]=-[2]_\sqq\haB_{[1,1],l-1}\ast[K_{\delta}].
		\end{align*}
		Moreover, by \eqref{haThij1} and \eqref{pi_ij,1+ anc co}, 
		\begin{align}
			\label{eq:Theta111Ol}
			&\big[\widehat{\Theta}_{[1,1],1}, [\co(k\vec{c})]\big]\\\notag
			=&\big[\pi_{[1,2],1}-(\sqq+\sqq^{-1})\pi_{[1,1],1},[\co(k\vec{c})]\big]\\\notag
			=&\big(\frac{\sqq^{-2}}{\sqq-\sqq^{-1}}-\frac{\sqq^{-1}(\sqq+\sqq^{-1})}{\sqq-\sqq^{-1}}\big)\big[[S_{1,0}^{(p_1)}], [\co(k\vec{c})]\big]\\\notag
			=&-\frac{1}{\sqq-\sqq^{-1}}\cdot \sqq^{-1}(q-1)\big([\co((k+1)\vec{c})]-[\co((k-1)\vec{c})]\ast[K_{\delta}]\big)\\\notag
			=&-[\co((k+1)\vec{c})]+[\co((k-1)\vec{c})]\ast[K_{\delta}].
		\end{align}
		Therefore, we have
		\begin{align*}
			&\Big[\haTh_{[1,1],1}, \big[[\co(k\vec{c})], \widehat{B}_{[1,1],l}\big]_{\sqq}  + \big[ \widehat{B}_{[1,1],l-1},[\co((k+1)\vec{c})]\big]_{\sqq}\Big]
			\\
			=&\Big[\haTh_{[1,1],1}, \big[[\co(k\vec{c})], \haB_{[1,1],l}\big]_{\sqq}\Big]  + \Big[\haTh_{[1,1],1}, \big[\haB_{[1,1],l-1},[\co((k+1)\vec{c})] \big]_{\sqq}\Big]
			\\
			=&\Big[\big[\haTh_{[1,1],1}, [\co(k\vec{c})]\big], \haB_{[1,1],l}\Big]_{\sqq}+\Big[[\co(k\vec{c})], \big[\haTh_{[1,1],1},  \haB_{[1,1],l}\big]\Big]_{\sqq}
			\\
			&+\Big[\big[\haTh_{[1,1],1}, \haB_{[1,1],l-1}\big], [\co((k+1)\vec{c})]\Big]_{\sqq}+\Big[\haB_{[1,1],l-1}, \big[\haTh_{[1,1],1},  [\co((k+1)\vec{c})]\big]\Big]_{\sqq}\\
			=&\Big[- [\co((k+1)\vec{c})]+ [\co((k-1)\vec{c})]*[K_\de], \haB_{[1,1],l}\Big]_\sqq+[2]_\sqq\Big[[\co(k\vec{c})], \haB_{[1,1],l+1}\Big]_{\sqq}\\
			&+[2]_\sqq\Big[ \haB_{[1,1],l}, [\co((k+1)\vec{c})]\Big]_{\sqq}+\Big[\haB_{[1,1],l-1}, - [\co((k+2)\vec{c})]+ [\co(k\vec{c})]*[K_\de]\Big]_\sqq\\
			=& -\Big(\big[[\co((k+1)\vec{c})],\haB_{[1,1],l}\big]_{\sqq}+\big[\haB_{[1,1],l-1},[\co((k+2)\vec{c})]\big]_{\sqq}\Big)
			\\
			&+\Big(\big[[\co((k-1)\vec{c})],\haB_{[1,1],l}\big]_{\sqq}+\big[\haB_{[1,1],l-1},[\co(k\vec{c})]\big]_{\sqq}\Big)*[K_\de] \\
			&+[2]_\sqq\Big(\big[[\co(k\vec{c})],\haB_{[1,1],l+1}\big]_{\sqq}+\big[\haB_{[1,1],l},[\co((k+1)\vec{c})]\big]_{\sqq}\Big).
		\end{align*}
		Observe that twisting with $\vec{c}$ induces an automorphism on $\iH(\X_\bfk,\varrho)$. Then \eqref{coBi1equiv} holds for any $l\geq0$ by induction. Similarly, one can also prove it for $l<0$  by using 
		$\big[\widehat{\Theta}_{[2,1],1}, [\co(l\vec{c})]\big]=-[2]_\sqq\haB_{[1,1],l-1}\ast[K_{\delta}]$.
	\end{proof}

	\subsection{The relation \eqref{Serrestari1}}
	\label{sec:intrelationsI}

	We shall prove \eqref{Serrestari1} by induction on $l$. First we check \eqref{Serrestari1} for $l=0,-1$.
	For this we need to introduce the following notations.
	
	For any $m\geq 0$,  define
	\begin{align}
		\label{eq:THstarm+}
		\haTh_{\star,m}^+:=&\frac{1}{(q-1)^2\sqq^{m-1}}\sum_{0\neq f:\co(s\vec{c})\rightarrow \co((m+s)\vec{c}+\vec{x}_1) } [\coker f],
		\\
		\label{eq:THstarm-}
		\haTh_{\star,m}^-:=&\frac{1}{(q-1)^2\sqq^{m-1}}\sum_{0\neq f:\co(s\vec{c}+\vec{x}_1)\rightarrow \co((m+s)\vec{c}) } [\coker f],
	\end{align}
	which are independent of $s\in\Z$.
	For convenience we set $\haTh_{{\star},m}^\pm=0$ for $m<0$.

	\begin{lemma}\label{Theta-p-m}
		For any $m\in\Z$,
		\begin{align}\label{Theta+}
			\haTh_{\star,m}^+=&
			q \haTh_{\star, m-2}^+*[K_{\delta}]+\frac{q}{q-1}\big[[S_{1,1}],\haTh_{\star, m}\big]_{\sqq^{-2}}-\frac{q}{q-1}\big[\haTh_{\star, m-1}, [S_{2,0}^{(p_1-1)}]\big]_{\sqq^{-2}}*[K_{\alpha_{11}}],
			\\
			\label{Theta-}
			\haTh_{\star,m}^-=&
			q\haTh_{\star, m-2}^-*[K_{\delta}]+\frac{\sqq}{q-1}\big[\haTh_{\star, m-1}, [S_{1,0}^{(p_1-1)}]\big]_{\sqq^{-2}}-\frac{q\sqq}{q-1}\big[[S_{2,1}],\haTh_{\star, m-2}\big]_{\sqq^{-2}}*[K_{\de-\alpha_{21}}].
		\end{align}
	\end{lemma}
	
	\begin{proof}
		For $m\leq0$, the above formulas hold obviously by definition. We assume $m>0$ in the following.
		
		For any $r,s\geq 0$, an easy computation shows that
		\begin{align}
			\label{S10S10pS20p}
			[S_{1,0}^{(p_1-1)}]*[S_{1,0}^{(rp_1)}\oplus S_{2,0}^{(sp_1)}]=& [S_{1,0}^{(p_1-1)}\oplus S_{1,0}^{(rp_1)}\oplus S_{2,0}^{(sp_1)}],\\
			\label{S10pS20pS10}
			[S_{1,0}^{(rp_1)}\oplus S_{2,0}^{(sp_1)}]*[S_{1,0}^{(p_1-1)}]=&\frac{1}{q}\Big([S_{1,0}^{(rp_1)}\oplus S_{2,0}^{(sp_1)}\oplus S_{1,0}^{(p_1-1)}]+(q-1)[S_{1,0}^{((r+1)p_1-1)}\oplus S_{2,0}^{(sp_1)}]\notag\\
			&+q(q-1)\sqq[S_{1,0}^{(rp_1)}\oplus S_{2,1}^{((s-1)p_1+1)}]*[K_{\de-\alpha_{21}}]\Big).
		\end{align}
		Here and in the following, it is understood that the terms involving $S_{2,1}^{((s-1)p_1+1)}$ do not exist if $s=0$. 
		Hence it follows by combining \eqref{S10S10pS20p} and \eqref{S10pS20pS10} that
		\begin{align}
			\label{la3}
			\big[[S_{1,0}^{(rp_1)}\oplus S_{2,0}^{(sp_1)}],[S_{1,0}^{(p_1-1)}]\big]_{\sqq^{-2}}
			=&\frac{q-1}{q}[S_{1,0}^{((r+1)p_1-1)}\oplus S_{2,0}^{(sp_1)}]\notag\\
			&+(q-1)\sqq[S_{1,0}^{(rp_1)}\oplus S_{2,1}^{((s-1)p_1+1)}]*[K_{\de-\alpha_{21}}].
		\end{align}

		For any non-zero map $f:\co\rightarrow\co(k\vec{c})$,
		assume $f=x_1^{rp_1}x_2^{sp_1}\cdot g$ for some $r,s\geq 0$, $x_1\nmid g$ and $x_2\nmid g$.
		So $\coker(f)\cong S_{1,0}^{(rp_1)}\oplus S_{2,0}^{(sp_1)}\oplus M$, where $M$ has no direct summand supported at $\bla_1$ or $\bla_2$.
		Note that any two torsion sheaves supported at distinct points have zero Hom and $\Ext^1$-spaces. Hence $[M]*[X]=[M\oplus X]=[X]*[M]$ for any $X\in\scrt_{\bla_1}\coprod \scrt_{\bla_2}$.
		Therefore,
		\begin{align*}
			&\big[[\coker f],[S_{1,0}^{(p_1-1)}]\big]_{\sqq^{-2}}=[M]*\big[[S_{1,0}^{(rp_1)}\oplus S_{2,0}^{(sp_1)}],[S_{1,0}^{(p_1-1)}]\big]_{\sqq^{-2}}\\
			\stackrel{(\ref{la3})}{=} &[M]*\Big(\frac{q-1}{q}[S_{1,0}^{((r+1)p_1-1)}\oplus S_{2,0}^{(sp_1)}]
			+(q-1)\sqq[S_{1,0}^{(rp_1)}\oplus S_{2,1}^{((s-1)p_1+1)}]*[K_{\de-\alpha_{21}}]\Big)\\
			=&\frac{q-1}{q}[M\oplus S_{1,0}^{((r+1)p_1-1)}\oplus S_{2,0}^{(sp_1)}]+(q-1)\sqq [M\oplus S_{1,0}^{(rp_1)}\oplus S_{2,1}^{((s-1)p_1+1)}]*[K_{\de-\alpha_{21}}]\\
			=&\frac{q-1}{q}[\coker f_1]+(q-1)\sqq [\coker f_2]*[K_{\de-\alpha_{21}}],
		\end{align*}
		where $f_1=x_1^{(r+1)p_1-1}x_2^{sp_1}\cdot g:\co(\vec{x}_1)\rightarrow \co((k+1)\vec{c})$ and $f_2=x_1^{rp_1}x_2^{(s-1)p_1+1}\cdot g:\co(\vec{c})\rightarrow \co(k\vec{c}+\vec{x}_2)$. Here it is understood that $f_2$ and then the term involving $[\coker f_2]$ does not exist if $s=0$. 
		
		Hence,  \begin{align*}
			&\sum_{0\neq f:\co\rightarrow \co(m\vec{c})} \big[[\coker f],[S_{1,0}^{(p_1-1)}]\big]_{\sqq^{-2}}\\
			=&\frac{q-1}{q}\sum_{0\neq f_1:\co(\vec{x}_1)\rightarrow \co((m+1)\vec{c})} [\coker f_1]+(q-1)\sqq\sum_{0\neq f_2:\co(\vec{c})\rightarrow \co(m\vec{c}+\vec{x}_2)} [\coker f_2]*[K_{\de-\alpha_{21}}].
		\end{align*}
		
		Then it follows from \eqref{def:Theta star}, \eqref{eq:THstarm+} and \eqref{eq:THstarm-} that
		\begin{align}\label{S_10andTheta}
			\big[\haTh_{\star, m},[S_{1,0}^{(p_1-1)}]\big]_{\sqq^{-2}}
			=(\sqq-\sqq^{-1})\haTh_{\star, m+1}^-+(\sqq-\sqq^{-1})\sqq\varrho(\haTh_{\star, m-1}^+)*[K_{\de-\alpha_{21}}].\end{align}
		
		Similarly, for any $r\geq 0$ we have
		\begin{align*}
			\big[[S_{1,1}], [S_{1,0}^{(rp_1)}\oplus S_{2,0}^{(sp_1)}]\big]_{\sqq^{-2}}
			=\frac{q-1}{q}[S_{1,1}^{(rp_1+1)}\oplus S_{2,0}^{(sp_1)}]+(\sqq-\sqq^{-1}) [S_{1,0}^{(rp_1)}\oplus S_{2,0}^{(sp_1-1)}]*[K_{\alpha_{11}}],
		\end{align*}
		which implies
		\begin{align*}
			&\big[[S_{1,1}],\sum_{0\neq f:\co\rightarrow \co(m\vec{c})} [\coker f]\big]_{\sqq^{-2}}\\
			&\qquad=\frac{q-1}{q}\sum_{0\neq f_1:\co\rightarrow \co(m\vec{c}+\vec{x}_1)} [\coker f_1]
			+(\sqq-\sqq^{-1}) \sum_{0\neq f_2:\co(\vec{x}_2)\rightarrow \co(m\vec{c})} [\coker f_2]*[K_{\alpha_{11}}].
		\end{align*}
		Then by definition we have
		\begin{align}\label{S_11andTheta}
			\big[[S_{1,1}],\haTh_{\star, m}\big]_{\sqq^{-2}}=\frac{q-1}{q}\haTh_{\star,m}^++(\sqq-\sqq^{-1}) \varrho(\haTh_{\star,m}^-)*[K_{\alpha_{11}}].
		\end{align}
		
		Combining (\ref{S_10andTheta}) and (\ref{S_11andTheta}), we obtain (\ref{Theta+}) and (\ref{Theta-}).
	\end{proof}
	
	\begin{lemma} For any $l\geq 1$,
		\begin{align}
			\label{eq:commcolcox1}
			\big[ [\co(l\vec{c})] ,&[\co(\vec{x}_1)] \big]_{\sqq^{-1}} +\big[ [\co],[\co(l\vec{c}+\vec{x}_1)]\big]_{\sqq^{-1}}
			\\\notag
			&=
			(\sqq-\sqq^{-1})^2 \big( \haTh_{\star,l}^+*[K_\co]-{\sqq^{-1}}{\varrho}(\haTh_{\star,l}^-)*[K_{\co(\vec{x}_1)}]\big).
		\end{align}
	\end{lemma}
	
	\begin{proof}
		By using the same proof of \cite[Corollary A.2]{LR21}, we obtain
		\begin{align}
			\label{co(l) and x1}
			&[\co(l\vec{c})]* [\co(\vec{x}_1)]-\sqq^{-1}[\co(l\vec{c}+\vec{x}_1)]*[\co]
			\\\notag&\qquad=
			\sqq^{-l-1}\Big([\co(l\vec{c})\oplus\co(\vec{x}_1)]-[\co(l\vec{c}+\vec{x}_1)\oplus\co]\Big).
		\end{align}
		
		Observe that for $l\geq1$,
		\begin{align*}
			&[\co(\vec{x}_1)]*[\co(l\vec{c})]\\
			=&\sqq^{-l}[\co(\vec{x}_1)\oplus \co(l\vec{c})]
			+\sqq^{-l-1}\sum_{0\neq f:\co(\vec{x}_1)\rightarrow \co(l\vec{c})} [\coker {\varrho (f)}]*[K_{\co(\vec{x}_1)}]\\
			=&\sqq^{-l}[\co(\vec{x}_1)\oplus \co(l\vec{c})]+\sqq^{-2}(q-1)^2{\varrho}(\haTh_{\star,l}^-)*[K_{\co(\vec{x}_1)}],
		\end{align*}
		and
		\begin{align*}
			[\co]*[\co(l\vec{c}+\vec{x}_1)]
			=&\sqq^{-l-1}[\co\oplus \co(l\vec{c}+\vec{x}_1)]+\sqq^{-l-1}\sum_{0\neq g:\co\rightarrow {\varrho}(\co(l\vec{c}+\vec{x}_1))} [\coker {\varrho} (g)]*[K_{\co}]\\
			=&\sqq^{-l-1}[\co\oplus \co(l\vec{c}+\vec{x}_1)]+\sqq^{-2}(q-1)^2\haTh_{\star,l}^+*[K_{\co}].
		\end{align*}
		Combining the above formulas, we obtain \eqref{eq:commcolcox1}.
	\end{proof}

	Define
	\begin{align}
		\label{def:S}
		\mathcal{S}:=\{S_{1,j},S_{2,j}\mid 2\leq j\leq p_1-1\}\cup \{S_{ij}\mid 3\leq i\leq \bt, 1\leq j\leq p_i-1\}.
	\end{align}
	Then the left perpendicular category ${}^\perp\mathcal{S}$ is equivalent to {$\coh(\mathbb{Y}_{\bfk})$}, where $\mathbb{Y}_{\bfk}$ is a weighted projective line of weight type $\bp=(2,2)$. By similar arguments as in \S\ref{Embedding from projective line to weighted projective line}, we have an embedding
	\begin{equation}
		\label{the embedding functor from Y to X}
		F_{\mathbb{X},\mathbb{Y}}:\iH(\mathbb{Y}_{\bfk},\varrho) \longrightarrow \iH(\X_\bfk,\varrho).
	\end{equation}
	
	In $\coh(\mathbb{Y}_{\bfk})$, we have  $2\vec{x}_1=2\vec{x}_2=\vec{c}$. 
	In this case, by definition we have
	\begin{align*}
		\haTh_{\star,l}^+(\vec{x}_1)=\sqq\haTh_{\star,l+1}^-,\qquad \haTh_{\star,l}^-(\vec{x}_1)
		=\sqq^{-1} \haTh_{\star,l-1}^+.
	\end{align*}
	Then by twisting with  $\vec{x}_1$ in the formulas \eqref{co(l) and x1} and \eqref{eq:commcolcox1}, we obtain
	
	\begin{align}
		\label{co(l+x1) and c}
		&[\co(l\vec{c}+\vec{x}_1)]*[\co(\vec{c})]-\sqq^{-1}[\co((l+1)\vec{c})]* [\co(\vec{x}_1)]
		\\\notag
		&\qquad=
		\sqq^{-l-1}\Big([\co(l\vec{c}+\vec{x}_1)\oplus\co(\vec{c})]-[\co((l+1)\vec{c})\oplus\co(\vec{x}_1)]\Big)
	\end{align}
	and
	\begin{align}
		\label{eq:commcolx1co}
		\big[ [\co(l\vec{c}+\vec{x}_1)], &[\co(\vec{c})] \big]_{\sqq^{-1}} +\big[ [\co(\vec{x}_1)],[\co((l+1)\vec{c})]\big]_{\sqq^{-1}}
		\\\notag
		&=
		(\sqq-\sqq^{-1})^2 \big(\sqq \haTh_{\star,l+1}^-*[K_{\co(\vec{x}_1)}]-{\sqq^{-2}}{\varrho}(\haTh_{\star,l-1}^+)*[K_{\co(\vec{c})}]\big).
	\end{align}
	Observe that each term in \eqref{co(l+x1) and c} and \eqref{eq:commcolx1co} belongs to the image of $F_{\mathbb{X},\mathbb{Y}}$. By using $F_{\mathbb{X},\mathbb{Y}}$ we know that these two formulas hold in $\iH(\X_\bfk,\varrho)$.

	
	\begin{lemma}
		\label{lem:SerreO11}
		For any $k_1,k_2\in\Z$, we have
		\begin{align}
			\label{SerreO11}
			\widehat{\SS}(k_1,k_2\mid 0;\star,[1,1])=&(1-q)^2\widehat{\R}(k_1,k_2\mid 0;\star,[1,1]).
		\end{align}
	\end{lemma}
	
	\begin{proof}
		Note that 
		\begin{align}
			\label{exchange relation of co and S_11}
			[\co(k\vec{c})]*[S_{1,1}]=\sqq[S_{1,1}]*[ \co(k\vec{c})] -(q-1)  [\co(k\vec{c}+\vec{x}_1)].
		\end{align}
		Hence,
		\begin{align*}
			&\widehat{S}(k_1,k_2\mid 0;\star,[1,1])\\
			=&[S_{1,1}]*[\co(k_1\vec{c})]*[\co(k_2\vec{c})] -[2]_\sqq [\co(k_1\vec{c})]*[S_{1,1}]*[\co(k_2\vec{c})]
			+[\co(k_1\vec{c})]*[\co(k_2\vec{c})]*[S_{1,1}]
			\\
			=&\big( \sqq^{-1} [\co(k_1\vec{c})]*[S_{1,1}]+(\sqq-\sqq^{-1})[\co(k_1\vec{c}+\vec{x}_1)]\big)*[\co(k_2\vec{c})]-[2]_\sqq [\co(k_1\vec{c})]*[S_{1,1}]*[\co(k_2\vec{c})]
			\\
			&+[\co(k_1\vec{c})]*\big(\sqq[S_{1,1}]*[\co(k_2\vec{c})]-(q-1)[\co(k_2\vec{c}+\vec{x}_1)]\big)
			\\
			=&(\sqq-\sqq^{-1})  [\co(k_1\vec{c}+\vec{x}_1)]*[\co(k_2\vec{c})]-(q-1) [\co(k_1\vec{c})]* [\co(k_2\vec{c}+\vec{x}_1)].
		\end{align*}
		For $k_1=k_2$, the desired formula \eqref{SerreO11} holds by noting that
		\begin{align*}
			&\widehat{\SS}(k_1,k_1\mid 0;\star,[1,1])
			=\widehat{S}(k_1,k_1\mid 0;\star,[1,1])
			\\
			=&-(q-1)  
			\big[[\co(k_1\vec{c})], [\co(k_1\vec{c}+\vec{x}_1)]\big]_{\sqq^{-1}}
			\\
			=&-\sqq^{-1}(q-1)^2 [S_{1,1}]*[K_{ \co(k_1\vec{c}) }]
			\\
			{\stackrel{\eqref{eqnSSkk1}}{=}}&(1-q)^2\widehat{\R}(k_1,k_1\mid 0;\star,[1,1]).
		\end{align*}
		Without loss of generality, we assume $k_1<k_2$ in the following. Then we have
		\begin{align*}
			&\widehat{\SS}(k_1,k_2\mid 0;\star,[1,1])\\
			=&(\sqq-\sqq^{-1})  [\co(k_1\vec{c}+\vec{x}_1)]*[\co(k_2\vec{c})]-(q-1) [\co(k_1\vec{c})]* [\co(k_2\vec{c}+\vec{x}_1)]
			\\
			&+ (\sqq-\sqq^{-1})   [\co(k_2\vec{c}+\vec{x}_1)]*[\co(k_1\vec{c})]-(q-1) [\co(k_2\vec{c})]* [\co(k_1\vec{c}+\vec{x}_1)]
			\\
			=&-(q-1) \Big( \big[[\co(k_2\vec{c})], [\co(k_1\vec{c}+\vec{x}_1)]\big]_{\sqq^{-1}}+\big[[\co(k_1\vec{c})], [\co(k_2\vec{c}+\vec{x}_1)]\big]_{\sqq^{-1}}\Big)
			\\
			\stackrel{\eqref{eq:commcolcox1}}{=}&\sqq^{-2}(q-1)^3\Big({\sqq^{-1}}{\varrho}(\haTh_{\star,k_2-k_1}^-)*[K_{\co(k_1\vec{c}+\vec{x}_1)}]-
			\haTh_{\star,k_2-k_1}^+*[K_{\co(k_1\vec{c})}]\Big).
		\end{align*}
		If $k_2=k_1+1$, then
		\eqref{SerreO11} holds since
		\begin{align*}
			&\widehat{\SS}(k_1,k_1+1\mid 0;\star,[1,1])\\
			=&\sqq^{-2}(q-1)^3{\sqq^{-1}}{\varrho}(\haTh_{\star,1}^-)*[K_{\co(k_1\vec{c}+\vec{x}_1)}]-
			\sqq^{-2}(q-1)^3\haTh_{\star,1}^+*[K_{\co(k_1\vec{c})}]\\
			=&{\sqq^{-2}}[2]_\sqq(q-1)^2[S_{2,0}^{(p_1-1)}]*[K_{\co(k_1\vec{c}+\vec{x}_1)}]
			-(q-1)^2\big[[S_{1,1}], \haTh_{\star,1}\big]_{\sqq^{-2}}*[K_{\co(k_1\vec{c})}]
			\\
			=&(q-1)^2\Big({-\sqq^{-1}[2]_\sqq\haB_{[1,1],-1}}*[K_\de] -\big[\haB_{[1,1],0},\haTh_{\star,1}\big]_{\sqq^{-2}}\Big)*[K_{\co(k_1\vec{c})}]
			\\
			\stackrel{\eqref{eqnSSkk+1}}{=}&(1-q)^2\widehat{\R}(k_1,k_1+1\mid 0;\star,[1,1]),
		\end{align*}
		where the second equality follows by (\ref{Theta+}) and (\ref{Theta-}).
		
		If $k_2>k_1+1$,
		then by (\ref{Theta+}) and (\ref{Theta-}), we have
		\begin{align*}
			&\widehat{\SS}(k_1,k_2\mid 0;\star,[1,1])-q\widehat{\SS} (k_1,k_2-2\mid 0;\star,[1,1])*[K_{\delta}]
			\\
			=&\sqq^{-2}(q-1)^3\cdot \sqq^{-1}\big(\varrho(\haTh_{\star,k_2-k_1}^-)-
			q\varrho(\haTh_{\star,k_2-k_1-2}^-)*[K_{\delta}]\big)*[K_{\co(k_1\vec{c}+\vec{x}_1)}]
			\\
			&-\sqq^{-2}(q-1)^3\big(\haTh_{\star,k_2-k_1}^+-
			q\haTh_{\star,k_2-k_1-2}^+*[K_{\delta}]\big)*[K_{\co(k_1\vec{c})}]\\
			=&\sqq^{-2}(q-1)^3\cdot \sqq^{-1}
			\big(\frac{\sqq}{q-1}\big[\haTh_{\star, k_2-k_1-1}, [S_{2,0}^{(p_1-1)}]\big]_{\sqq^{-2}}\\
			&\qquad\qquad\qquad\quad-\frac{q\sqq}{q-1}\big[[S_{1,1}],\haTh_{\star, k_2-k_1-2}\big]_{\sqq^{-2}}*[K_{\de-\alpha_{11}}]\big)
			*[K_{\co(k_1\vec{c}+\vec{x}_1)}]
			\\
			&-\sqq^{-2}(q-1)^3
			\big(\frac{q}{q-1}\big[[S_{1,1}],\haTh_{\star, k_2-k_1}\big]_{\sqq^{-2}}
			\\
			&\qquad\qquad\qquad\quad-\frac{q}{q-1}\big[\haTh_{\star, k_2-k_1-1}, [S_{2,0}^{(p_1-1)}]\big]_{\sqq^{-2}}*[K_{\alpha_{11}}]\big)
			*[K_{\co(k_1\vec{c})}]\\
			=&(q-1)^2 
			\big(\sqq^{-2}\big[\haTh_{\star, k_2-k_1-1}, [S_{2,0}^{(p_1-1)}]\big]_{\sqq^{-2}}-\big[[S_{1,1}],\haTh_{\star, k_2-k_1-2}\big]_{\sqq^{-2}}*[K_{\de-\alpha_{11}}]\big)
			*[K_{\co(k_1\vec{c}+\vec{x}_1)}]
			\\
			&-(q-1)^2 \big(\big[[S_{1,1}],\haTh_{\star, k_2-k_1}\big]_{\sqq^{-2}}-\big[\haTh_{\star, k_2-k_1-1}, [S_{2,0}^{(p_1-1)}]\big]_{\sqq^{-2}}*[K_{\alpha_{11}}]\big)
			*[K_{\co(k_1\vec{c})}]\\
			=&(q-1)^2 
			\big(\sqq^{-1}[2]_\sqq\big[\haTh_{\star, k_2-k_1-1}, [S_{2,0}^{(p_1-1)}]\big]_{\sqq^{-2}}*
			[K_{\co(k_1\vec{c}+\vec{x}_1)}]
			-\big[[S_{1,1}],\haTh_{\star, k_2-k_1}\big]_{\sqq^{-2}}
			*[K_{\co(k_1\vec{c})}]
			\\
			&-\big[[S_{1,1}],\haTh_{\star, k_2-k_1-2}\big]_{\sqq^{-2}}
			*[K_{\co(k_1+1)\vec{c})}]\big)\\
			=&(q-1)^2\Big(-[2]_\sqq\big[\haTh_{\star, k_2-k_1-1}, \haB_{[1,1],-1}\big]_{\sqq^{-2}}*[K_{\co((k_1+1)\vec{c})}]
			-\big[\haB_{[1,1],0},\haTh_{\star, k_2-k_1}\big]_{\sqq^{-2}}*[K_{\co(k_1\vec{c})}]\\
			&-\big[\haB_{[1,1],0},\haTh_{\star, k_2-k_1-2}\big]_{\sqq^{-2}}*[K_{\co((k_1+1)\vec{c})}]\Big)
			\\
			\stackrel{(\ref{Theta and B k_2>k_1+1})}{=}&(q-1)^2\Big(\widehat{\R}(k_1,k_2\mid 0;\star,[1,1])-q\widehat{\R} (k_1,k_2-2\mid 0;\star,[1,1])*[K_{\delta}]\Big).
		\end{align*}
		Therefore, the desired formula \eqref{SerreO11} follows by induction.
	\end{proof}

	\begin{proposition}
		\label{prop:SerreO11l}
		For any $k_1,k_2,l\in\Z$, we have
		\begin{align}
			\label{SerreO11l}
			\widehat{\SS}(k_1,k_2\mid l;\star,[1,1])=&(1-q)^2\widehat{\R}(k_1,k_2\mid l;\star,[1,1]).
		\end{align}
	\end{proposition}
	
	\begin{proof}
		We only need to consider the case $l\ge0$, since the other case $l<0$ is similar.
		By Lemma \ref{lem:SerreO11}, we do induction on $l$, and assume that \eqref{SerreO11l} holds for $0\le l'\leq l$.
		
		Using \eqref{humbvl} and \eqref{eq:Theta111Ol}, a direct computation shows
		\begin{align}
			\label{eq:THO}
			&\big[\haTh_{[1,1],1},\widehat{\SS}(k_1,k_2\mid l;\star,[1,1])\big]
			\\\notag
			=&-\widehat{\SS}(k_1+1,k_2\mid l;\star,[1,1])+\widehat{\SS}(k_1-1,k_2\mid l;\star,[1,1])*[K_\de]
			\\\notag
			&-\widehat{\SS}(k_1,k_2+1\mid l;\star,[1,1])+\widehat{\SS}(k_1,k_2-1\mid l;\star,[1,1])*[K_\de]
			\\\notag
			&+[2]_\sqq\widehat{\SS}(k_1,k_2\mid l+1;\star,[1,1])
			\\\notag
			=&(1-q)^2\Big(-\widehat{\R}(k_1+1,k_2\mid l;\star,[1,1])+\widehat{\R}(k_1-1,k_2\mid l;\star,[1,1])*[K_\de]
			\\\notag
			&-\widehat{\R}(k_1,k_2+1\mid l;\star,[1,1])+\widehat{\R}(k_1,k_2-1\mid l;\star,[1,1])*[K_\de]\Big)
			\\\notag
			&+[2]_\sqq\widehat{\SS}(k_1,k_2\mid l+1;\star,[1,1])
			\\\notag
			=&[2]_\sqq\widehat{\SS}(k_1,k_2\mid l+1;\star,[1,1]).
		\end{align}
		Here the second equality uses the inductive assumptions, and the third equality uses the following equality obtained from \eqref{eq:haRkk}:
		\begin{align*}
			\widehat{\R}(k_1,k_2\pm 1\mid l;\star,[1,1])=\widehat{\R}(k_1\mp1,k_2\mid l;\star,[1,1])*[K_\de]^{\pm1}.
		\end{align*}
		
		On the other hand, using the same argument as \eqref{haThetaij1starm}, one can obtain that
		$\big[\widehat{\Theta}_{[1,1],1}, \widehat{\Theta}_{\star,m}\big]=0$ for any $m>0$. Together with \eqref{humbvl}, we have
		\begin{align}
			\label{eq:THR}
			&[\haTh_{[1,1],1},\widehat{\R}(k_1,k_2\mid l;\star,[1,1])]
			= [2]_\sqq\widehat{\R}(k_1,k_2\mid l+1;\star,[1,1]).
		\end{align}
		Comparing \eqref{eq:THO} and \eqref{eq:THR}, using the inductive assumption, one can obtain
		$$\widehat{\SS}(k_1,k_2\mid l+1;\star,[1,1])=(1-q)^2\widehat{\R}(k_1,k_2\mid l+1;\star,[1,1]).$$
		The proof is completed.
	\end{proof}

	\subsection{The relations \eqref{humstarbvl} and \eqref{eq:HHstari1}}
	
	\begin{proposition}
		\label{prop:humstarbvl}
		Relation \eqref{humstarbvl} holds in $\iH(\X_\bfk,\varrho)$ for any $m>0$ and $l\in\Z$.
	\end{proposition}
	
	\begin{proof}
		The proof is completely same to \cite[Proposition 8.5]{LR21}, hence omitted here.
	\end{proof}

	\begin{proposition}
		Relation \eqref{eq:HHstari1} holds in $\iH(\X_\bfk,\varrho)$ for any $m,r>0$.
	\end{proposition}
	
	\begin{proof}
		The proof is completely same to \cite[Proposition 8.6]{LR21}, hence omitted here.
	\end{proof}

	\section{Relations between $\star$ and $[i,1]$, II}
	\label{sec:Relationsstari1 II}
	
	In this section, we shall prove the following formulas in $\iH(\X_\bfk,\varrho)$:
	\begin{align}
		\label{humbvlstar}
		&\big[\widehat{H}_{[i,1],m},[\co(l\vec{c})]\big]=-\frac{[m]_\sqq}{m} [\co((l+m)\vec{c})]+\frac{[m]_\sqq}{m} [\co((l-m)\vec{c})]*[K_{m\de}],
		\\
		\label{Serrei1star}
		&\widehat{\SS}(k_1,k_2\mid l;[i,1],\star)=\begin{cases}(1-q)^2\widehat{\R}(k_1,k_2\mid l; [i,1],\star), & \text{ if }\varrho(i)=i,\\
			0,& \text{ if }\varrho(i)\neq i,
		\end{cases}
	\end{align}
	for any $k_1,k_2,l\in\Z$, $m>0$, $1\leq i\leq \bt$.

	In fact, if $\varrho(i)=i$, these formulas can be derived by the same computations as \cite[\S9]{LR21}.
	So we can assume $\varrho(i)\neq i$ in the remaining of this section. 
	Without loss of generality, we assume $i=1$ and $\varrho(i)=2$ throughout this section.


	\subsection{The relation \eqref{humbvlstar}}
	
	In this subsection, we shall prove \eqref{humbvlstar}. Note that twisting with $\vec{c}$ induces an  automorphism on $\iH(\X_\bfk,\varrho)$. By using Lemma \ref{lem:equiv}~(2), it is equivalent to prove that for any $r>0$,
	\begin{align}
		\label{humbvlstar-reform}
		&\big[ \haTh_{[1,1],r},[\co]\big]+\big[\haTh_{[1,1],r-2},[\co]\big]*[K_\de]
		\\
		& \qquad\qquad= \sqq^{-1}\big[\haTh_{[1,1],r-1},[\co(\vec{c})]\big]_{\sqq^2}+\sqq\big[\haTh_{[1,1],r-1},[\co(-\vec{c})]\big]_{\sqq^{-2}}*[K_\de].
		\notag
	\end{align}


	The proof of \eqref{humbvlstar-reform} needs the following Lemmas \ref{Theta 11 times co} and \ref{co times Theta 11}, whose proofs are technical. 
	\begin{lemma}
		\label{Theta 11 times co}
		For any $r>0$, we have
		\begin{align*}
			\haTh_{[1,1],r}\ast [\![\co]\!]=&\frac{\sqq^{r+1}}{q-1}\sum\limits_{|\mu|=r}\bn(\ell(\mu))[\![\co\oplus S_{1,0}^{(\mu)} ]\!]+\sqq^{r-1}\sum\limits_{[N]\in\cm_{r\de}}\bn(\ell(N)-1)[\![\co\oplus N]\!]\\&
			-\sum\limits_{k>0}
			\sum\limits_{|\mu|=r-k}\sqq^{r-2k+1}\cdot \bn(\ell(\mu))[\![\co(k\vec{c})\oplus S_{1,0}^{(\mu)} ]\!]\\
			&+
			\sum\limits_{k>0}\sum\limits_{[N]\in\cm_{1,(r-k)\de}} \sqq^{r-2k-1}(1-q)\cdot \bn(\ell(N)-1)[\![\co(k\vec{c})\oplus N]\!].
		\end{align*}
	\end{lemma}
	
	\begin{proof}
		The proof is the same as \cite[Lemma 9.3]{LR21}, hence omitted here.  
	\end{proof}

	\begin{lemma}
		\label{co times Theta 11}
		For any $r>0$, we have
		\begin{align*}
			[\![\co]\!]\ast \haTh_{[1,1],r}=&\frac{\sqq^{r+1}}{q-1}\sum\limits_{|\mu|={r}}   \bn(\ell(\mu))[\![\co\oplus S_{1,0}^{(\mu)}]\!]+\sqq^{r-1}\sum\limits_{[N]\in\cm_{r\de}}  \bn(\ell(N)-1)[\![\co\oplus N]\!]\\
			&- \sum\limits_{k>0}\sum\limits_{|\mu|={r-k}}  \sqq^{r-2k+1}\cdot \bn(\ell(\mu))[\![\co(-k\vec{c})\oplus S_{1,0}^{(\mu)}]\!]\ast [K_{k\delta}]
			\\
			&+\sum\limits_{k> 0}\sum\limits_{[N]\in\cm_{1,(r-k)\de}}  \sqq^{r-2k-1}(1-q)\cdot {\bf n}(\ell(N)-1)[\![\co(-k\vec{c})\oplus N]\!]\ast [K_{k\delta}].
		\end{align*}
	\end{lemma}

	\begin{proof}
		For any $[M]\in\cm_{1,r\de}$, we have $\langle\widehat{\co}, \widehat{M} \rangle=r$, $\Ext_{\X}^1(\co, M)=0$ and $\Hom_{\X}(\co, \varrho(M))\cong\bfk^r$. For any non-zero morphism $f:\co\rightarrow \varrho(M)$, we have $\Im(f)\cong S_{2,0}^{(kp_1)}$ or $S_{2,0}^{(k  p_1-1)}$ for some $k>0$.
		
		If $\Im(f)=S_{2,0}^{(kp_1)}$, then $\Ker(f)=\co(-k\vec{c})$ and $[\coker(\varrho(f))]\in\cm_{1,(r-k)\de}$; see \cite[Sublemma B.7]{LR21}. For any $[N]\in\cm_{(r-k)\de}$, we have
		\begin{align*}&|\{f:\co\rightarrow \varrho(M)\,|\,\Ker(f)\cong \co(-k\vec{c}), \coker(\varrho(f))\cong N\}|\\=&|\{g:S_{2,0}^{(kp_1)}\rightarrowtail \varrho(M)\,|\, \coker(\varrho(g))\cong N\}|\\
			=&|\{\varrho(g):S_{1,0}^{(kp_1)}\rightarrowtail M\,|\, \coker(\varrho(g))\cong N\}|
			\\
			=&F_{N,S_{1,0}^{(kp_1)}}^{M}\cdot|\Aut(S_{1,0}^{(kp_1)})|.
		\end{align*}
		Similarly, if $\Im(f)=S_{2,0}^{(kp_1-1)}$, then $\Ker(f)=\co(-k\vec{c}+\vec{x}_2)$ and $[\coker(\varrho(f))]\in\cm_{(r-k)\de+\alpha_{11}}$. For any $[N']\in\cm_{(r-k)\de+\alpha_{11}}$,
		\begin{align*}&|\{f:\co\rightarrow \varrho(M)\,|\,\Ker(f)\cong \co(-k\vec{c}+\vec{x}_2), \coker(\varrho(f))\cong N'\}|\\=&|\{g:S_{2,0}^{(kp_1-1)}\rightarrowtail \varrho(M)\,|\, \coker(\varrho(g))\cong N'\}|\\
			=&F_{N',S_{1,0}^{(kp_1-1)}}^{M}\cdot|\Aut(S_{1,0}^{(kp_1-1)})|.
		\end{align*}
		Then by \eqref{Hallmult1}, we have
		\begin{align*}
			&[\co]\ast[M]
			=\sqq^{-r}\sum\limits_{k\geq 0}\sum\limits_{[N]\in\cm_{(r-k)\de}} F_{N,S_{1,0}^{(kp_1)}}^{M}\cdot|\Aut(S_{1,0}^{(kp_1)})|\cdot[\co(-k\vec{c})\oplus N]\ast [K_{k\delta}]\\
			&\quad+\sqq^{-r}\sum\limits_{k\geq 1}\sum\limits_{[N']\in\cm_{(r-k)\de+\alpha_{11}}} F_{N',S_{1,0}^{(kp_1-1)}}^{M}\cdot|\Aut(S_{1,0}^{(kp_1-1)})|\cdot[\co(-k\vec{c}+\vec{x}_2)\oplus N']\ast [K_{k\delta-\alpha_{21}}];
		\end{align*}
		or equivalently, by using the Riedtmann-Peng formula, we have
		\begin{align}\label{formula of co times M}
			&[\![\co]\!]\ast[\![M]\!]
			=\sqq^{-r}\sum\limits_{k\geq 0}\sum\limits_{[N]\in\cm_{(r-k)\de}}q^{r-k}\frac{|\Ext^1(N, S_{1,0}^{(kp_1)})_{M}|}{|\Hom(N, S_{1,0}^{(kp_1)})|}\cdot[\![\co(-k\vec{c})\oplus N]\!]\ast [K_{k\delta}] \\\notag
			&\quad+\sqq^{-r}\sum\limits_{k\geq 1}\sum\limits_{[N']\in\cm_{(r-k)\de+\alpha_{11}}}q^{r-k+1}\frac{|\Ext^1(N', S_{1,0}^{(kp_1-1)})_{M}|}{|\Hom(N', S_{1,0}^{(kp_1-1)})|}\cdot[\![\co(-k\vec{c}+\vec{x}_2)\oplus N']\!]\ast [K_{k\delta-\alpha_{21}}].
		\end{align}
		
		Similar arguments show that
		\begin{align}\label{the formula of co times S lambda}
			&[\![\co]\!]\ast  [\![S_{1,0}^{(\lambda)}]\!]= \sqq^{-r}\sum\limits_{k\geq 0}\sum\limits_{|\mu|=r-k} q^{r-k}\frac{|\Ext^1(S_{1,0}^{(\mu)}, S_{1,0}^{(kp_1)})_{S_{1,0}^{(\lambda)}}|}{|\Hom(S_{1,0}^{(\mu)}, S_{1,0}^{(kp_1)})|}\cdot[\![\co(-k\vec{c})\oplus S_{1,0}^{(\mu)}]\!]\ast [K_{k\delta}].
		\end{align}
		
		Based on \eqref{eq:THeta11r}, \eqref{formula of co times M}--\eqref{the formula of co times S lambda}, the remaining proof is completely the same as \cite[Lemma 9.4]{LR21}, hence omitted here.
	\end{proof}
	
	\begin{proposition}
		\label{prop:humbvlstar-reform}
		Relation \eqref{humbvlstar-reform} holds in $\iH(\X_\bfk,\varrho)$ for any $r>0$.
	\end{proposition}
	
	\begin{proof}
		By using Lemmas \ref{Theta 11 times co} and \ref{co times Theta 11}, (\ref{humbvlstar-reform}) follows by easy (but tedious) cancellations, which we omit here.
	\end{proof}

	\subsection{The relation \eqref{Serrei1star}}
	
	In this subsection, we shall prove the relation \eqref{Serrei1star}. First, let us give a lemma.
	
	\begin{lemma}
		\label{lem:serrel=0}
		The following formula holds in $\iH(\X_\bfk,\varrho)$ for any $l\in\Z$:
		\begin{align*}
			\widehat{\SS}(0,0 |l;[1,1],\star)=0.
		\end{align*}
	\end{lemma}
	
	\begin{proof}
		By \eqref{exchange relation of co and S_11} we obtain
		\begin{align*}
			&\widehat{\SS}(0,0\mid l;[1,1],\star)
			\\
			=&[S_{1,1}]*[S_{1,1}]*[ \co(l\vec{c})]-[2]_\sqq[S_{1,1}]*[ \co(l\vec{c})]*[S_{1,1}]+[\co(l\vec{c})]*[S_{1,1}]*[S_{1,1}]\\
			=&[S_{1,1}]*\big(\sqq^{-1}[\co(l\vec{c})]*[S_{1,1}]+(\sqq-\sqq^{-1}) [\co(l\vec{c}+\vec{x}_1)]\big)-[2]_\sqq[S_{1,1}]*[ \co(l\vec{c})]*[S_{1,1}]
			\\
			&+\big(\sqq[S_{1,1}]*[\co(l\vec{c})] -(q-1)[\co(l\vec{c}+\vec{x}_1)]  \big)*[S_{1,1}]
			\\
			=&(\sqq-\sqq^{-1})\big[[S_{1,1}], [\co(l\vec{c}+\vec{x}_1)]\big]_\sqq\\
			=&0.
		\end{align*}
	\end{proof}
	
	\begin{proposition}
		Relation \eqref{Serrei1star} holds in $\iH(\X_\bfk,\varrho)$ for any $k_1,k_2,l\in\Z$.
	\end{proposition}

	\begin{proof}
		With the help of Proposition \ref{prop:humbvlstar-reform},  \eqref{humbvl} for $\nu=[1,1]$ and $\mu=[1,1]$, we have 
		\begin{align*}
			\big[\haH_{[1,1],k_1}, \widehat{\SS}(0,0\mid l;[1,1],\star) \big]=&\frac{[2k_1]_\sqq}{[k_1]}\widehat{\SS}(k_1,0\mid l;[1,1],\star) -\frac{[k_1]_\sqq}{k_1} \widehat{\SS}(0,0\mid l+k_1;[1,1],\star)
			\\
			&+\frac{[k_1]_\sqq}{k_1} \widehat{\SS}(0,0\mid l-k_1;[1,1],\star)*[K_{k_1\delta}]
		\end{align*}
		for any $k_1\geq1$.
		Then $\widehat{\SS}(k_1,0\mid l;[1,1],\star)=0$ by Lemma \ref{lem:serrel=0}.
		Furthermore, it is easy to get that $\widehat{\SS}(k_1,0\mid l;[1,1],\star)=0$ for any $k_1,l\in\Z$ by computing $\big[\haH_{[2,1],k_1}, \widehat{\SS}(0,0\mid l;[1,1],\star) \big]$.

		Similarly, \begin{align*} \big[\haH_{[1,1],k_2}, &\widehat{\SS}(k_1,0\mid l;[1,1],\star) \big]=\frac{[2k_2]_\sqq}{k_2}\widehat{\SS}(k_1+k_2,0\mid l;[1,1],\star)+\frac{[2k_2]_\sqq}{k_2}\widehat{\SS}(k_1,k_2\mid l;[1,1],\star)\\
			&-\frac{[k_2]_\sqq}{k_2}\widehat{\SS}(k_1,0\mid l+k_2;[1,1],\star)+\frac{[k_2]_\sqq}{k_2}\widehat{\SS}(k_1,0\mid l-k_2;[1,1],\star)*[K_{k_2\de}]
		\end{align*}
		for any $k_2\ge1$. Then we have 
		$\widehat{\SS}(k_1,k_2\mid l;[1,1],\star)=0$ for any $k_1,l\in\Z,k_2\geq0$. 
		By considering $\big[\haH_{[2,1],k_2}, \widehat{\SS}(k_1,0\mid l;[1,1],\star) \big]$, we have
		$\widehat{\SS}(k_1,k_2\mid l;[1,1],\star)=0$ for any $k_1,k_2,l\in\Z$.
	\end{proof}

	

\end{document}